\documentclass[12pt]{amsart}
\usepackage{soul}

\usepackage[OT2,T1]{fontenc}
\newcommand{\mifody}{%
  \renewcommand\rmdefault{wncyr}%
  \renewcommand\sfdefault{wncyss}%
  \renewcommand\encodingdefault{OT2}%
  \normalfont
  \selectfont}
\newcommand{\Sh}{\mathop{\mifody{\textsf{Sh}}}}

\usepackage[vcentermath]{youngtab}

\usepackage{mathabx}

\usepackage{times}
\usepackage{mathptmx}

\usepackage[mathscr]{eucal}
\usepackage{amsthm}
\usepackage{mathtools}
\usepackage{ytableau}

\usepackage{epstopdf}
\usepackage{ifpdf}
\ifpdf
\usepackage{t-angles}
\usepackage[matrix,arrow,curve,pdf]{xy} %%%,frame,
\else
\usepackage{t-angles}
\usepackage[matrix,arrow,curve,dvips]{xy} %%%,frame,
\UsePSspecials{dvips}
\fi

\newcommand{\qDim}{\mathsf{qdim}}

\newcommand{\q}{\mathsf{q}}

\newcommand{\theZ}{\mathsf{Z}}
\newcommand{\theH}{\mathsf{H}}
\newcommand{\theF}{\mathsf{F}}
\newcommand{\theFF}{\widetilde{\mathsf{F}}}
\newcommand{\nofrac}[2]{#1/#2}

\newcommand{\tbl}{\mathfrak{t}}

\newcommand{\Threetableaux}[3]{\Bigl(\,#1\,,\,\;#2\,,\,\;#3\,\Bigr)}

\usepackage{accents}

\advance\textheight\topskip
\allowdisplaybreaks

\newcommand{\shape}[1]{\overline{#1}}
\newcommand{\theD}{\mathsf{D}}
\newcommand{\Mass}{\mathsf{M}}
\newcommand{\Massi}{\widetilde{\mathsf{M}}}
\newcommand{\Ent}{\mathsf{E}}

\newcommand{\up}{\raisebox{2pt}{\mbox{\tiny$\scriptscriptstyle\nearrow$}}}
\newcommand{\down}{\raisebox{2pt}{\mbox{\tiny$\scriptscriptstyle\searrow$}}}

\newcommand{\hupinv}{h^{-1\,\up}}
\newcommand{\hdowninv}{h^{-1\,\down}}
\newcommand{\gup}{g^{\up}}
\newcommand{\gdown}{g^{\down}}
\newcommand{\gupinv}{g^{-1\,\up}}
\newcommand{\gdowninv}{g^{-1\,\down}}

\newcommand{\Ytableau}[1]{\mbox{\scriptsize$\displaystyle\begin{ytableau}#1\end{ytableau}$}}
\newcommand{\Ytableaui}[1]{\raisebox{3pt}{\mbox{\scriptsize$\begin{ytableau}#1\end{ytableau}$}}}
\newcommand{\Ytableauii}[1]{\raisebox{3pt}{\mbox{\tiny$\begin{ytableau}#1\end{ytableau}$}}}

\newcommand{\Chain}{\mathscr{T}}
\newcommand{\ket}[1]{\left|#1\right\rangle}
\newcommand{\ketdemo}[3]{\Bigl|\mbox{\small$(#1)$},\;\mbox{\small$#2\,,\;#3$}\Bigr\rangle}
\newcommand{\dgr}{\lambda}
\newcommand{\Dgr}{\Lambda}

\newcommand{\res}{\mathop{\mathrm{res}}\nolimits}

\newcommand{\bmcat}{{\mathbb{B}_{\bullet\circ}}}

\newcommand{\elli}{\ell'}
\newcommand{\ellii}{\ell}
\newcommand{\dgri}{\dgr'}
\newcommand{\dgrii}{\dgr}
\newcommand{\Dgri}{\Dgr'}
\newcommand{\Dgrii}{\Dgr}
\newcommand{\tbli}{\tbl'}
\newcommand{\tblii}{\tbl}

\newcommand{\Specht}[1]{\mathscr{S}^{#1}}
\newcommand{\SN}[1]{\mathscr{L}^{#1}}

\newcommand{\abmcat}{{\mathbb{A}_{\bullet\circ}}}

\newcommand{\omegaii}{\widetilde{\omega}}

\newcommand{\HyperT}{\mathbf{T}}
\newcommand{\hyp}{\mathsf{t}}
\newcommand{\qwB}{\mathsf{qw}\kern-1.8pt\mathscr{B}\kern-1.5pt}
\newcommand{\Hecke}{\mathscr{H}}

\newcommand{\JMBare}{J}
\newcommand{\Cas}{C}
\newcommand{\Casii}{\widetilde{\Cas}}

\newcommand{\Cons}{\mathscr{A}}
\newcommand{\modV}{\mathscr{V}}
\newcommand{\Jww}{\widetilde{\mathscr{J}}}
\newcommand{\Jwww}{\overline{\mathscr{J}}}
\newcommand{\Jbb}{\mathscr{J}}
\newcommand{\Uu}{\mathscr{U}}

\newcommand{\tensor}{\mathbin{\otimes}}

\newcommand{\bu}{\overline{u}}
\newcommand{\bx}{\overline{y}}

\newcommand{\Star}{\raisebox{-.5pt}{\mbox{\large$\star$}}}
\newcommand{\Circ}{\raisebox{-.5pt}{\mbox{$\circ$}}}
\newcommand{\Ast}{\raisebox{-.5pt}{\mbox{\large$*$}}}

\newcommand{\Bobject}[1]{\object{\mbox{\Large$#1$}}}
\newcommand{\fobject}[1]{\object{\scriptstyle #1}}
\newcommand{\upfobject}[1]{\object{\raisebox{4pt}{\mbox{$\scriptstyle#1$}}}}
\newcommand{\ffobject}[1]{\object{\scriptscriptstyle #1}}

\newcommand{\htilde}{\widetilde{h}}
\newcommand{\htildeii}{\widetilde{\htilde}}

\newcommand{\EE}{\mathscr{E}}

\newcommand{\wt}{\mathsf{wt}}
\newcommand{\weight}{\mathsf{weight}}
\newcommand{\Wt}{\textup{\textsf{Wt}}}

\newcommand{\pos}[1]{\{#1\}}
\newcommand{\Pos}[2]{#1[#2]}

\newcommand{\BareWt}{\kern2pt\overline{\kern-2pt\mathsf{Wt}\kern-2pt}\kern2pt}

\newcommand{\Orb}{\mathsf{Orb}}
\newcommand{\Orbii}{\mathsf{Orb}_{1'2'}}

\newcommand{\acts}{\mathop{\pmb{.}}\nolimits}

\newcommand{\HookDistance}{\mathsf{\Gamma}\!}
\newcommand{\Deltav}{\Delta_{\text{v}}}
\newcommand{\Deltah}{\Delta_{\text{h}}}
\newcommand{\Id}{\text{id}}

\newcommand{\remove}{\mathbin{\boxbackslash}}
\newcommand{\superimpose}{\lozenge}

\newcommand{\Corner}[1]{\mbox{\setul{}{.4pt}
\ul{\,$#1\kern2pt$}\kern-.6pt\rule[-4pt]{.4pt}{2.5ex}}\kern1pt}

\newcommand{\Corners}[1]{\mbox{\setul{}{1pt}\ul{\,$#1\kern2pt$}\kern-.6pt\rule[-4pt]{1.2pt}{2.5ex}}\kern1pt}

\newcommand{\corners}[1]{\mbox{\setul{}{1pt}\ul{$\,\scriptstyle#1\kern1pt$}\kern-.7pt\rule[-4.1pt]{1.2pt}{2.1ex}}}

\newcommand{\Cocorners}[1]{\;\rule[-2pt]{1.2pt}{2.4ex}\kern-.5pt\accentset{\rule{.9em}{1pt}}{{\;#1}}\;}

\newcommand{\cocorners}[1]{\;\rule[-1pt]{1.2pt}{1.75ex}\kern-2pt\accentset{\rule{.8em}{1pt}}{{\;\;\scriptstyle#1\,}}}

\newcommand{\cocornersi}[1]{\;\rule[-1pt]{1.2pt}{1.85ex}\kern-2pt\accentset{\rule{2.5em}{1pt}}{{\;\;\scriptstyle#1\,}}}

\newcommand{\cocornersii}[1]{\;\rule[-1pt]{1.2pt}{1.85ex}\kern-2pt\accentset{\rule{4em}{1pt}}{{\;\;\scriptstyle#1\,}}}

\newcommand{\eqdef}{\mathrel{\mbox{\raisebox{.6pt}{\footnotesize{:}}}{=}}}

\newcommand{\bref}[1]{\textup{\textbf{\ref{#1}}}}

\newcommand{\mfrac}[2]{\raisebox{.3pt}{\mbox{\small$\displaystyle\frac{#1}{#2}$}}}
\newcommand{\ffrac}[2]{\raisebox{.5pt}{\mbox{\footnotesize$\displaystyle\frac{#1}{#2}$}}}
\newcommand{\fffrac}[2]{\raisebox{.9pt}{\mbox{\tiny$\displaystyle\frac{#1}{#2}$}}}

\newcommand{\oC}{\mathbb{C}}
\newcommand{\oZ}{\mathbb{Z}}
\newcommand{\oN}{\mathbb{N}}

\newcommand{\oS}{\mathbb{S}}

\numberwithin{equation}{section}

\makeatletter
\@addtoreset{equation}{section}
\@addtoreset{subsubsection}{section}

\def\@secnumfont{\bfseries}
\def\subsubsection{\@startsection{subsubsection}{3}%
  \z@{.5\linespacing\@plus.7\linespacing}{-.5em}%
  {\normalfont\bfseries}}
\def\paragraph{\@startsection{paragraph}{4}%
  \z@\z@{-\fontdimen2\font}%
  \normalfont\bfseries}
\def\subparagraph{\@startsection{subparagraph}{5}%
  \z@\z@{-\fontdimen2\font}%
  \normalfont\bfseries}

\makeatother

\swapnumbers
\newtheorem{Thm}[subsection]{Theorem}
\newtheorem{thm}[subsubsection]{Theorem}
\newtheorem{Lemma}[subsection]{Lemma}
\newtheorem{lemma}[subsubsection]{Lemma}

\theoremstyle{definition}

\newtheorem{Rem}[subsection]{Remark}%[section]
\newtheorem{rem}[subsubsection]{Remark}
\newtheorem{example}[subsubsection]{Example}

\begin{document}
 \title[Quantum walled Brauer algebra]{Quantum walled Brauer algebra:
  commuting families, Baxterization, and representations}

\author[Semikhatov]{A.\,M.\;Semikhatov}
\author[Tipunin]{I.\,Yu.\;Tipunin}

\address{Lebedev Physics Institute, Moscow 119991,
  Russia\hfill\mbox{}\linebreak \texttt{asemikha@gmail.com}, \
  \texttt{tipunin@gmail.com}}

\begin{abstract}
  For the quantum walled Brauer algebra, we construct its Specht
  modules and (for generic parameters of the algebra) seminormal
  modules.  The latter construction yields the spectrum of a commuting
  family of Jucys--Murphy elements.  We also propose a Baxterization
  prescription; it involves representing the quantum walled Brauer
  algebra in terms of morphisms in a braided monoidal category and
  introducing parameters into these morphisms, which allows
  constructing a ``universal transfer matrix'' that generates
  commuting elements of the algebra.
\end{abstract}

\maketitle
\thispagestyle{empty}

\section{Introduction}
We study the quantum (``quantized'') walled Brauer algebra $\qwB$ and
its representations.\footnote{Speaking of an algebra rather than
  algebras, we mean a particular member of a family $\qwB_{m,n}$ of
  quantum walled Brauer algebras with $m,n\geq 1$; this algebra,
  moreover, depends on parameters, as we discuss below.}  The
classical version of the algebra was introduced in
\cite{[BCLLS],[Koi]} in the context of generalized Schur--Weyl
duality: the algebra was shown to centralize the $g\ell(N)$ action on
``mixed'' tensor products $X^*{}^{\otimes m}\tensor X^{\otimes n}$ of
the natural $g\ell(N)$ representation and its dual; for special
parameter values, the walled Brauer algebra centralizes the action of
$g\ell(M|N)$~\cite{[Serg],[BS]}. \ The structure of the algebra was
explored in~\cite{[CDvDM]}.  The quantum version of the algebra was
introduced in~\cite{[L],[Ha],[KM],[Kos]}, and its role as the
centralizer of $U_q(g\ell_N)$ on the mixed tensor product
$X^*{}^{\otimes m}\tensor X^{\otimes n}$ was elucidated
in~\cite{[DDS-1],[DDS-2]}.  We also note a recent ``super'' extension
of quantum walled Brauer algebras in~\cite{[BGJKW]}.

We here view $\qwB_{m,n}$ following \cite{[Ha]} (also see~\cite{[LR]})
as a diagram algebra, with the diagrams supplied by a braided monoidal
category.\footnote{To be developed in full rigor, this approach would
  require a braided version of Deligne's category \cite{[Deligne]},
  which can apparently be done, but is beyond the scope of this paper
  (see \cite{[CW]} and the references therein).}  We use the diagrams
representing category morphisms to construct two types of commutative
families of $\qwB_{m,n}$ elements:
\begin{itemize}
\item[(i)] a family of ``conservation laws'' following from a
  Baxterization procedure, and

\item[(ii)] a family of Jucys--Murphy elements
  $\JMBare(n)_{2},\dots,\JMBare(n)_{m+n}\in\qwB_{m,n}$ (which we
  diagonalize, as is discussed below).
\end{itemize}

The commuting $\qwB_{m,n}$ elements in item~(i) are called
``conservation laws'' or ``Hamiltonians'' in view of the physical
interpretation of (quantum)\pagebreak[3] walled Brauer algebras as
pertaining to integrable models of statistical mechanics, e.g., the
$t$--$J$ model (see \cite{[F],[LF],[C]} and the references therein).
The commuting elements follow by expanding a ``universal transfer
matrix'' $\Cons_{m, n}(z)$ in the spectral parameter $z\in\oC$
introduced by a trick that generalizes the well-known Baxterization of
Hecke algebras.  For $\qwB_{m,n}$, however, it is applied not only to
the algebra generators but also to morphisms in a braided monoidal
category $\bmcat$ that do not belong to $\qwB_{m,n}$.  Using these, it
is straightforward to show that $[\Cons_{m, n}(z),\Cons_{m, n}(w)]=0$;
at the same time it turns out that despite the occurrence of
``extraneous'' morphisms in the definition,
$\Cons_{m,n}(z)\in\qwB_{m,n}\tensor\oC(z)$, whence a commuting family
of $\qwB_{m,n}$ elements follows by expanding around a suitable value
of the spectral parameter.

By borrowing more from integrable systems of statistical mechanics,
but staying within the $\qwB_{m,n}$ algebra, it would be quite
interesting to diagonalize our conservation laws by Bethe-ansatz
techniques, but we here solve only a more modest diagonalization
problem, the one for a family of Jucys--Murphy elements (which are in
fact the $z\to\infty$ limits of ``monodromy matrices'' constructed
similarly to the transfer matrices).

We recall that Jucys--Murphy elements were originally introduced for
the symmetric group algebra~\cite{[Juc],[Mu-81],[Mu-83]} and were then
discussed for some other diagram algebras \cite{[GL-hecke],[HMR],[Li]}
and in even broader contexts (see, e.g.,
\cite{[Re],[Na],[LR],[Ram-semi]}), up to the generality of cellular
algebras~\cite{[MS],[GG]}. \ As with the conservation laws, we here
define Jucys--Murphy elements $\JMBare(n)_{j}\in\qwB_{m,n}$,
$j=2,\dots,m+n$, in a way that makes their commutativity manifest, but
involves braided-category morphisms not from $\qwB_{m,n}$; again,
their apparent ill-definedness is in fact superficial, and the
$\JMBare(n)_{j}$ can eventually be expressed via relatively explicit
formulas in terms of generators.

For generic values of the algebra parameters, when $\qwB_{m,n}$ is
semisimple, we diagonalize the commuting family of Jucys--Murphy
elements in each irreducible representation.  This amounts to
constructing seminormal $\qwB_{m,n}$ representations.  Seminormal
bases$/$representations have been studied rather extensively for a
number of ``related'' algebras (see, e.g.,
\cite{[We-88],[Na],[Ram-semi],[Ram-skew],[MS],[Eny-semi]}), but
apparently not for~ $\qwB_{m,n}$.  The seminormal basis is made of
triples of Young tableaux of certain shapes, and an essential novelty
compared with the Hecke-algebra case is the ``mobile
elements''---corners of Young tableaux that can change their position.

For special parameter values, $\qwB_{m,n}$ is not semisimple and
seminormal representations may not exist.  By contrast, Specht modules
exist for all parameter values, are generically irreducible, and
become reducible at special parameter values (playing a role somewhat
similar to that played for Lie algebras by Verma modules).  The
$\qwB_{m,n}$ Specht modules have been discussed
in~\cite{[Eny-cellular],[RS1],[RS2]}. \ We construct them rather
explicitly, by extending our ``diagrammatic'' view of the $\qwB$
algebra to representations.  A construction that combines categorial
diagrams with Young tableaux is called the link-state representation
here, to emphasize a similarity (or the authors' prejudices regarding
this similarity) to a link-state construction for the Temperley--Lieb
algebra (see, e.g., \cite{[M],[RsA]}).

This paper is organized as follows.  In Sec.~\ref{sec:category}, we
introduce $\qwB_{m,n}$ as an algebra of tangles of a particular type
satisfying certain relations~\cite{[Ha]}.  This language naturally
suggests a construction for Jucys--Murphy elements.  In
Sec.~\ref{sec:baxter}, by ``introducing spectral parameters into
tangles,'' we obtain a Baxterization (not of the algebra, but of the
``ambient'' category), which allows us to construct commutative
families of algebra elements; the families depend on two parameters in
addition to the $\qwB_{m,n}$ parameters and are also determined by one
out of three expansion points of the universal transfer matrix.  In
Sec.~\ref{sec:specht}, we construct the $\qwB_{m,n}$ Specht modules in
terms of tangles in which the propagating lines (``defects'') end at
the boxes of two standard Young tableaux.  In
Sec.~\ref{sec:seminormal}, we construct the seminormal representations
and find the spectrum of Jucys--Murphy elements.

\section{$\qwB_{m,n}$ from braided monoidal
  categories}\label{sec:category}
In this section, we define the quantum walled Brauer algebra in terms
of diagrams (tangles); these can be thought of as being supplied by a
rigid braided monoidal category $\bmcat$; specifically, the category
objects are given by (ordered) collections of nodes of two sorts
($\bullet$ and~$\circ$) and morphisms are given by tangles on these
nodes.  Quantum walled Brauer algebras $\qwB_{m,n}$ with $m,n\in\oN_0$
are endomorphism algebras of $\bmcat$ objects; taking the
abelianization $\abmcat$ of $\bmcat$ then gives a category whose
simple objects are all simple $\qwB_{m,n}$-modules for $m,n\in\oN_0$.

\subsection{The category $\bmcat$}\label{the-cat}
We fix a braided $\oC$-vector space $X$ (a vector space with a
bilinear map $\psi:X\tensor X\to X\tensor X$ satisfying the
braid$/$Yang--Baxter equation). \ Objects in $\bmcat$ are tensor
products of $X$ and its dual $X^*$. \ We use the respective
notation~$\bullet$ and~$\circ$ for $X^*$ and $X$ and call them the
black and white objects.  Their tensor products are represented simply
as $\bullet{}\,{}\bullet=X^*\otimes X^*$,
$\bullet\circ\bullet=X^*\otimes X\otimes X^*$, and so on.

\subsubsection{Diagram notation for morphisms} We use the standard
pictorial notation,
%% (see, e.g., \cite{[Besp-TMF]},[Kaufman]) for the braiding, evaluation
%% and coevaluation morphisms, and so on.  
representing morphisms as tangles~\cite{[Tu]}
(cf.~\cite{[Ha],[LR]}). \ Characteristic examples of morphisms are
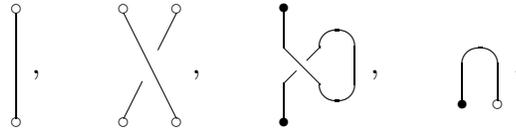
\begin{equation*}
  \begin{tangles}{l}
    \fobject{\circ}\\[-4.5pt]
    \vstr{200}\id\\[-4.5pt]
    \fobject{\circ}
  \end{tangles}\ \ ,\kern30pt 
  \begin{tangles}{l}
    \fobject{\circ}\step[2]\fobject{\circ}\\[-4.5pt]
    \vstr{200}\x\\[-4.5pt]
    \fobject{\circ}\step[2]\fobject{\circ}
  \end{tangles}\ \ ,\kern30pt 
  \begin{tangles}{l}
     \fobject{\bullet}\\[-4.5pt]
     \vstr{67}\hstr{67}\id\step[2]\coev\\
     \vstr{67}\hstr{67}\x\step[2]\id\\
     \vstr{67}\hstr{67}\id\step[2]\ev\\[-4.5pt]
     \fobject{\bullet}
  \end{tangles}\ \ ,\kern30pt
  \begin{tangles}{l}
    \hstr{67}\coev\\
    \hstr{67}\vstr{50}\id\step[2]\id\\[-4.5pt]
    \hstr{67}\fobject{\bullet}\step[2]\fobject{\circ}
  \end{tangles}\ \ ,
\end{equation*}
which are the identity morphism of $X$, the braiding
$X\tensor X\to X\tensor X$, the ribbon morphism on (or the twist of)
$X^*$, and the coevaluation of $X^*$ and $X$. \ All diagrams are
considered modulo isotopy (more precisely, a tangle is an isotopy
class of diagrams), and we assume that there are no triple
intersections of lines.

\subsubsection{Relations in $\mathrm{Mor}(\bmcat)$}\label{cat-relations}
The set of morphisms between two objects in $\bmcat$ is a $\oC$-vector
space of formal linear combinations of tangles modulo the relations
that we now describe and which depend on four complex parameters
$\alpha,\beta,\kappa,\kappa'\in\oC$.

The basic axiom is that the braiding of two white objects satisfies
the Hecke relation
\begin{equation}\label{Hecke-white}
  \begin{tangles}{l}
    \hstr{133}\fobject{\circ}\step\fobject{\circ}\\[-4.5pt]
    \hstr{133}\hx\\
    \hstr{133}\hx\\[-4.5pt]
    \hstr{133}\fobject{\circ}\step\fobject{\circ}
  \end{tangles}\ 
  = -\alpha\beta \ \
  \begin{tangles}{l}
    \fobject{\circ}\step[1.5]\fobject{\circ}\\[-4.5pt]
    \vstr{200}\id\step[1.5]\id\\[-4.5pt]
    \fobject{\circ}\step[1.5]\fobject{\circ}
  \end{tangles}\ \ 
  + (\alpha+\beta) \ \
  \begin{tangles}{l}
    \fobject{\circ}\step[2]\fobject{\circ}\\[-4.5pt]
    \vstr{200}\x\\[-4.5pt]
    \fobject{\circ}\step[2]\fobject{\circ}
  \end{tangles}\ \ .
\end{equation}

We also assume that the category is rigid, with the black and white
objects being dual to each other; this means that there are the
evaluation and coevaluation morphisms \ \ $\begin{tangles}{l}
  \hstr{67}\fobject{\bullet}\step[2]\fobject{\circ}\\[-4.3pt]
  \hstr{67}\vstr{67}\ev
  \end{tangles}$, \ \ $\begin{tangles}{l}
    \hstr{67}\fobject{\circ}\step[2]\fobject{\bullet}\\[-4.3pt]
    \hstr{67}\vstr{67}\ev
  \end{tangles}$, \ \ $\begin{tangles}{l}
    \hstr{67}\vstr{67}\coev\\[-4.5pt]
    \hstr{67}\fobject{\bullet}\step[2]\fobject{\circ}
\end{tangles}$\ \ , and \ \ $\begin{tangles}{l}
  \hstr{67}\vstr{67}\coev\\[-4.5pt]
  \hstr{67}\fobject{\circ}\step[2]\fobject{\bullet}
\end{tangles}$\
\ .  Moreover, the category is assumed to be ribbon, with the ribbon
map on black and white objects defined by two constants $\kappa$ and
$\kappa'$:
\begin{gather}\label{ribbon}
  \begin{tangles}{l}
    \fobject{\bullet}\\[-4.5pt]
    \vstr{67}\hstr{67}\id\step[2]\coev\\
    \vstr{67}\hstr{67}\x\step[2]\id\\
    \vstr{67}\hstr{67}\id\step[2]\ev\\[-4.5pt]
    \fobject{\bullet}
  \end{tangles}\ \
  = \ffrac{1}{\kappa}\ \ 
  \begin{tangles}{l}
    \fobject{\bullet}\\[-4.5pt]
    \vstr{200}\id
    \fobject{\bullet}
  \end{tangles}\ \
  \qquad\text{and}\qquad
  \begin{tangles}{l}
    \fobject{\circ}\\[-4.5pt]
    \vstr{67}\hstr{67}\id\step[2]\coev\\
    \vstr{67}\hstr{67}\x\step[2]\id\\
    \vstr{67}\hstr{67}\id\step[2]\ev\\[-4.5pt]
    \fobject{\circ}
  \end{tangles}\ \
  = \ffrac{1}{\kappa'}\ \ 
  \begin{tangles}{l}
    \fobject{\circ}\\[-4.5pt]
    \vstr{200}\id
    \fobject{\circ}
  \end{tangles}\ \ .
\end{gather}
We note that by duality, these relations imply that
\begin{gather}\label{ribbon1}
  \begin{tangles}{l}
    \hstr{67}\step[4]\fobject{\circ}\\[-4.5pt]
    \vstr{67}\hstr{67}\coev\step[2]\id\\
    \vstr{67}\hstr{67}\id\step[2]\x\\
    \vstr{67}\hstr{67}\ev\step[2]\id\\[-4.5pt]
    \hstr{67}\step[4]\fobject{\circ}
  \end{tangles}\ \
  = \ffrac{1}{\kappa}\ \ 
  \begin{tangles}{l}
    \fobject{\circ}\\[-4.5pt]
    \vstr{200}\id
    \fobject{\circ}
  \end{tangles}\ \
  \qquad\text{and}\qquad
  \begin{tangles}{l}
    \hstr{67}\step[4]\fobject{\bullet}\\[-4.5pt]
    \vstr{67}\hstr{67}\coev\step[2]\id\\
    \vstr{67}\hstr{67}\id\step[2]\x\\
    \vstr{67}\hstr{67}\ev\step[2]\id\\[-4.5pt]
    \hstr{67}\step[4]\fobject{\bullet}
  \end{tangles}\ \
  = \ffrac{1}{\kappa'}\ \ 
  \begin{tangles}{l}
    \fobject{\bullet}\\[-4.5pt]
    \vstr{200}\id
    \fobject{\bullet}
  \end{tangles}\ \ .
\end{gather}
Below, these relations are conveniently used in the form
\begin{gather}\label{kappas}
  \begin{tangles}{l}
    \hstr{67}\vstr{67}\coev\\
    \hstr{67}\vstr{67}\x\\[-4.5pt]
    \hstr{67}\fobject{\bullet}\step[2]\fobject{\circ}
  \end{tangles}\ \ = \kappa \ \ 
  \begin{tangles}{l}
    \hstr{67}\coev\\
    \hstr{67}\vstr{50}\id\step[2]\id\\[-4.5pt]
    \hstr{67}\fobject{\bullet}\step[2]\fobject{\circ}
  \end{tangles}
  \qquad\text{and}\qquad
  \begin{tangles}{l}
    \hstr{67}\vstr{67}\coev\\
    \hstr{67}\vstr{67}\x\\[-4.5pt]
    \hstr{67}\fobject{\circ}\step[2]\fobject{\bullet}
  \end{tangles}\ \ = \kappa' \ \ 
  \begin{tangles}{l}
    \hstr{67}\coev\\
    \hstr{67}\vstr{50}\id\step[2]\id\\[-4.5pt]
    \hstr{67}\fobject{\circ}\step[2]\fobject{\bullet}
  \end{tangles}\ \ .
\end{gather}

Composing \eqref{Hecke-white} with evaluation and coevaluation maps,
we readily obtain the Hecke-algebra relations for two black objects,
the same as for the white ones:
\begin{align}\label{Hecke-black}
  \begin{tangles}{l}
    \hstr{133}\fobject{\bullet}\step\fobject{\bullet}\\[-4.5pt]
    \hstr{133}\hx\\
    \hstr{133}\hx\\[-4.5pt]
    \hstr{133}\fobject{\bullet}\step\fobject{\bullet}
  \end{tangles}\ 
  &= -\alpha\beta \ \
  \begin{tangles}{l}
    \fobject{\bullet}\step[1.5]\fobject{\bullet}\\[-4.5pt]
    \vstr{200}\id\step[1.5]\id\\[-4.5pt]
    \fobject{\bullet}\step[1.5]\fobject{\bullet}
  \end{tangles}\ \ 
  + (\alpha+\beta) \ \
  \begin{tangles}{l}
    \fobject{\bullet}\step[2]\fobject{\bullet}\\[-4.5pt]
    \vstr{200}\x\\[-4.5pt]
    \fobject{\bullet}\step[2]\fobject{\bullet}
  \end{tangles}\ \ .
  \\
  \intertext{Simple manipulations (see \bref{demo-braid} below) also
  show that}
  \label{square-of-bad}
  \begin{tangles}{l}
    \hstr{133}\fobject{\bullet}\step\fobject{\circ}\\[-4.5pt]
    \hstr{133}\hx\\
    \hstr{133}\hx\\[-4.5pt]
    \hstr{133}\fobject{\bullet}\step\fobject{\circ}
  \end{tangles}\ 
  &= -\ffrac{1}{\alpha\beta} \ \
  \begin{tangles}{l}
    \fobject{\bullet}\step[1.5]\fobject{\circ}\\[-4.5pt]
    \vstr{200}\id\step[1.5]\id\\[-4.5pt]
    \fobject{\bullet}\step[1.5]\fobject{\circ}
  \end{tangles}\ \ 
  + \ffrac{\alpha+\beta}{\alpha \beta}\kappa \ \
  \begin{tangles}{l}
    \fobject{\bullet}\step[2]\fobject{\circ}\\[-4.5pt]
    \ev\\
    \coev\\[-4.5pt]
    \fobject{\bullet}\step[2]\fobject{\circ}
  \end{tangles} \ \ .
\end{align}
A similar formula holds with the colors interchanged and $\kappa$
replaced with $\kappa'$.

Compositions  with \ \ $\begin{tangles}{l}
  \hstr{67}\fobject{\circ}\step[2]\fobject{\bullet}\\[-4.3pt]
  \hstr{67}\vstr{67}\ev
\end{tangles}$ \ and \ $\begin{tangles}{l}
  \hstr{67}\fobject{\bullet}\step[2]\fobject{\circ}\\[-4.3pt]
  \hstr{67}\vstr{67}\ev
\end{tangles}$ \ \ yield
\begin{gather}\label{circlebw}
    \begin{tangles}{l}
    \Coev\\[-4.3pt]
    \fobject{\circ}\step[3]\fobject{\bullet}\\
    \Ev
  \end{tangles}\ \
  =\ffrac{\alpha\beta\kappa\kappa' + 1}{\kappa'(\alpha+\beta)}
  \qquad\text{and}\qquad
  \begin{tangles}{l}
    \Coev\\[-4.3pt]
    \fobject{\bullet}\step[3]\fobject{\circ}\\
    \Ev
  \end{tangles}\ \
  =\ffrac{\alpha\beta\kappa\kappa' + 1}{\kappa(\alpha+\beta)}.
\end{gather}

\subsubsection{Braid diagram manipulations}\label{demo-braid}As an
example of derivations with diagrams, we show how
\eqref{square-of-bad} follows from~\eqref{Hecke-black} and the axioms.
Composing~\eqref{Hecke-black} with a ``cup'' gives
\begin{equation*}
  \begin{tangles}{l}
    \hstr{66}\fobject{\circ}\step[2]\fobject{\bullet}\step[2]\fobject{\bullet}\\[-4.3pt]
    \hstr{66}\id\step[2]\x\\
    \hstr{66}\ev\step[2]\id\\[-4.3pt]
    \hstr{66}\step[4]\fobject{\bullet}
  \end{tangles}\ \
  = -\alpha\beta\ \
  \begin{tangles}{l}
    \hstr{66}\fobject{\circ}\step[2]\fobject{\bullet}\step[2]\fobject{\bullet}\\[-4.3pt]
    \hstr{66}\x\step[2]\id\\
    \hstr{66}\id\step[2]\ev\\[-4.3pt]
    \fobject{\bullet}
  \end{tangles}\ \
  + (\alpha+\beta) \ \
  \begin{tangles}{l}
    \hstr{66}\fobject{\circ}\step[2]\fobject{\bullet}\step[2]\fobject{\bullet}\\[-4.3pt]
    \hstr{66}\ev\step[2]\id\\
    \hstr{66}\step[4]\id\\[-4.3pt]
    \hstr{66}\step[4]\fobject{\bullet}
  \end{tangles}\ \ .
\end{equation*}
At the top of this relation, we attach \ \ 
$\begin{tangles}{l}
  \fobject{\bullet}\step[1]\fobject{\circ}\step[1]\fobject{\bullet}\\[-4.3pt]
  \vstr{50}\hx\step[1]\id\\[-4.3pt]
  \fobject{\circ}\step[1]\fobject{\bullet}\step[1]\fobject{\bullet}
\end{tangles}$\ , after which the equality takes the form
\begin{equation*}
  \begin{tangles}{l}
    \hstr{66}\fobject{\bullet}\step[2]\fobject{\circ}\step[2]\fobject{\bullet}\\[-4.3pt]
    \hstr{66}\vstr{150}\id\step[2]\ev\\
    \id\\[-4.3pt]
    \fobject{\bullet}
  \end{tangles}\ \
  = -\alpha\beta\ \
  \begin{tangles}{l}
    \hstr{66}\fobject{\bullet}\step[2]\fobject{\circ}\step[2]\fobject{\bullet}\\[-4.3pt]
    \hstr{66}\x\step[2]\id\\
    \hstr{66}\id\step[2]\xx\\
    \hstr{66}\ev\step[2]\id\\[-4.3pt]
    \hstr{66}\step[4]\fobject{\bullet}
  \end{tangles}\ \
  + (\alpha+\beta)\ \
  \begin{tangles}{l}
    \hstr{66}\fobject{\bullet}\step[2]\fobject{\circ}\step[2]\fobject{\bullet}\\[-4.3pt]
    \hstr{66}\x\step[2]\id\\
    \hstr{66}\ev\step[2]\id\\[-4.3pt]
    \hstr{66}\step[4]\fobject{\bullet}
  \end{tangles}\ \ .
\end{equation*}
But the first relation in~\eqref{ribbon} can be rewritten as\ \ 
$\begin{tangles}{l}
  \fobject{\bullet}\step[1]\fobject{\circ}\\[-4.3pt]
  \hstr{50}\ev
\end{tangles}\ \ = \fffrac{1}{\kappa}\ \ 
\begin{tangles}{l}
  \fobject{\bullet}\step[1]\fobject{\circ}\\[-4.3pt]
  \vstr{66}\vstr{50}\hx\\
  \vstr{66}\hstr{50}\ev
\end{tangles}
$\
, and hence the last term is equal to
$\kappa(\alpha+\beta)\ \ \begin{tangles}{l}
  \fobject{\bullet}\step[1]\fobject{\circ}\step[1]\fobject{\bullet}\\[-4.3pt]
  \vstr{66}\hstr{50}\ev\step[2]\id\\[-4.3pt]
  \step[2]\fobject{\bullet}
\end{tangles}$\
.  Composing the resulting relaton with a ``cap,'' we then obtain
\begin{equation*}
  \begin{tangles}{l}
    \fobject{\bullet}\step[1.5]\fobject{\circ}\\[-4.5pt]
    \vstr{200}\id\step[1.5]\id\\[-4.5pt]
    \fobject{\bullet}\step[1.5]\fobject{\circ}
  \end{tangles}\ 
  = -\alpha\beta \ \
  \begin{tangles}{l}
    \hstr{133}\fobject{\bullet}\step\fobject{\circ}\\[-4.5pt]
    \hstr{133}\hx\\
    \hstr{133}\hx\\[-4.5pt]
    \hstr{133}\fobject{\bullet}\step\fobject{\circ}
  \end{tangles}
  + \kappa(\alpha+\beta) \ \
  \begin{tangles}{l}
    \fobject{\bullet}\step[2]\fobject{\circ}\\[-4.5pt]
    \ev\\
    \coev\\[-4.5pt]
    \fobject{\bullet}\step[2]\fobject{\circ}
  \end{tangles} \ \ ,
\end{equation*}
which is~\eqref{square-of-bad}.

\subsection{The algebra $\qwB_{m,n}$
  \cite{[Ha]}}\label{qwB-from-braided}
The quantum walled Brauer algebra $\qwB_{m,n}$ is the algebra of
endomorphisms of the object (``mixed tensor space'')
\begin{equation*}
  \Chain_{m,n}=(X^*)^{\otimes m}\otimes X^{\otimes n}=
  \underbrace{\bullet\dots\bullet}_m\,
  \underbrace{\circ\dots\circ}_n\;.
\end{equation*}
Each such endomorphism can be represented by a tangle, considered
modulo the above relations, with $m$ black and $n$ white nodes on the
top edge and the same numbers of black and white nodes on the bottom
edge, and with strands connecting a white (black) node with a white
(black) node on the opposite edge or with a black (white) node on the
same edge.  If a strand connects nodes on the same edge, we call it an
arc (top or bottom depending on the edge); there must be the same
number of top and bottom arcs. The strands connecting different edges
are called defects\footnote{Propagating lines in another nomenclature
  in a similar context~\cite{[CDvDM]}.} (black or white depending on
the type of connected nodes).

A vertical \textit{wall} can be imagined to separate
$(X^*)^{\otimes m}=\bullet\dots\bullet$ from
$X^{\otimes n}=\circ\dots\circ$.  Arcs necessarily cross the wall,
while defects do not.\footnote{For each pair of positive integers
  $(m,n)$, the algebra $\qwB_{m,n}$ can also be represented as
  endomorphisms of an object in $\bmcat$ where the $m$ factors $X^*$
  and the $n$ factors $X$ are taken in a different order.  All such
  algebras are isomorphic because of the existence of braiding
  morphisms in the category $\bmcat$. \ For example, in the case
  $m=n,n\pm1$, there is a natural ``physical'' order of factors in the
  tensor product $X\otimes X^*\otimes X\otimes X^*\otimes\dots$, which
  represents a spin chain of alternating atoms of two sorts.}

\subsubsection{Numbering convention}\label{convention}
With the order of factors chosen as
$X^*{}^{\otimes m}\tensor X^{\otimes n}$ in what follows, we adopt the
convention that the nodes are enumerated ``from the wall'' outwards.
We also often use a primed collection of the integers for labeling the
(``black'' ) $X^*$ factors:
\begin{equation*}
  \begin{tangles}{l}
    \hstr{50}\fobject{\dots}\step[2]\fobject{\bullet}\step[2]\fobject{\bullet}\step[2]\fobject{\bullet}\step[2]\fobject{\bullet}\step[2]\fobject{\circ}\step[2]\fobject{\circ}\step[2]\fobject{\circ}\step[2]\fobject{\circ}\step[2]\fobject{\dots}\\
    \hstr{50}\fobject{\dots}\step[2]\fobject{4'}\step[2]\fobject{3'}\step[2]\fobject{2'}\step[2]\fobject{1'}\step[2]\fobject{1}\step[2]\fobject{2}\step[2]\fobject{3}\step[2]\fobject{4}\step[2]\fobject{\dots}
    \end{tangles}
\end{equation*}

\subsubsection{The $\qwB_{m,n}$ generators and relations} The tangles
are multiplied standardly (in accordance with our convention of
reading the diagrams from top down), with
relations~\eqref{Hecke-white}--\eqref{circlebw} applied whenever
needed, and then the algebra $\qwB_{m,n}$ is generated by the tangles
\cite{[DDS-1]}
\begin{align}\label{eq:g-j}
 g_j &=\ \ \begin{tangles}{l}
    \fobject{\bullet}\step[2]\fobject{\dots}\step[2]\fobject{\bullet}\step[2]\fobject{\bullet}\step[2]
     \fobject{\dots}\step[2]\fobject{\bullet}\step[2]\fobject{\circ}\step[2]
     \fobject{\dots}\step[2]\fobject{\circ}\\[-4.3pt]
    \vstr{167}\id\step[2]\step[2]\x\step[2]\step[2]\id\step[2]\id\step[2]\step[2]\id\\[-4.5pt]
    \fobject{\bullet}\step[2]\fobject{\dots}\step[2]\fobject{\bullet}\step[2]\fobject{\bullet}\step[2]
    \fobject{\dots}\step[2]\fobject{\bullet}\step[2]
\fobject{\circ}\step[2]
    \fobject{\dots}\step[2]\fobject{\circ}\\
    \fobject{m'}\step[5]\fobject{j+1'\step[1]j'}\step[5]\fobject{1'}\step[2]\fobject{1}\step[4]\fobject{n}
  \end{tangles}\,\,\,,\qquad j=1,\dots, m-1,
  \\
  \label{eq:EE}
  \EE &= \ \ \
  \begin{tangles}{l}
    \fobject{\bullet}\step[2]\fobject{\dots}\step[2]\fobject{\bullet}\step[2]\fobject{\bullet}\step[2]
    \fobject{\circ}\step[2]\fobject{\circ}\step[2]\fobject{\dots}\step[2]\fobject{\circ}\\[-4.3pt]
    \id\step[2]\step[2]\id\step[2]\ev\step[2]\id\step[2]\step[2]\id\\
    \id\step[2]\step[2]\id\step[2]\coev\step[2]\id\step[2]\step[2]\id\\[-4.5pt]
    \fobject{\bullet}\step[2]\fobject{\dots}\step[2]\fobject{\bullet}\step[2]
    \fobject{\bullet}\step[2]\fobject{\circ}\step[2]\fobject{\circ}\step[2]\fobject{\dots}\step[2]\fobject{\circ}\\
 \fobject{m'}\step[4]\fobject{2'}\step[2]\fobject{1'}\step[2]\fobject{1}\step[2]\fobject{2}\step[4]\fobject{n}
  \end{tangles}
  \\
  \intertext{and}
  \label{eq:h-i}
  h_i&=\ \ \begin{tangles}{l}
    \fobject{\bullet}\step[2]\fobject{\dots}\step[2]\fobject{\bullet}\step[2]\fobject{\circ}\step[2]
     \fobject{\dots}\step[2]\fobject{\circ}\step[2]\fobject{\circ}\step[2]\fobject{\dots}\step[2]\fobject{\circ}\\[-4.3pt]
    \vstr{167}\id\step[2]\step[2]\id\step[2]\id\step[2]\step[2]\x\step[2]\step[2]\id\\[-4.5pt]
    \fobject{\bullet}\step[2]\fobject{\dots}\step[2]\fobject{\bullet}\step[2]\fobject{\circ}\step[2]
    \fobject{\dots}\step[2]\fobject{\circ}\step[2]
\fobject{\circ}\step[2]
    \fobject{\dots}\step[2]\fobject{\circ}\\
    \fobject{m'}\step[4]\fobject{1'}\step[2]\fobject{1}\step[5]\fobject{i\step[1]i+1}\step[5]\fobject{n}
  \end{tangles}\,\,\, ,\qquad i=1,\dots, n-1.
\end{align}

By Eqs.~\eqref{Hecke-white} and~\eqref{Hecke-black}, the $g_j$ and
$h_i$ are standard generators of two commuting Hecke algebras
$\Hecke_{m}(\alpha,\beta)$ and $\Hecke_{n}(\alpha,\beta)$, 
\begin{gather}
  g_j^2=(\alpha+\beta)g_j-\alpha \beta \cdot 1,\qquad
  h_i^2=(\alpha+\beta)h_i-\alpha \beta \cdot 1,\\
  g_j h_i = h_i g_j,
\end{gather}
where $1\leq j\leq m-1$ and $1\leq i\leq n-1$, and, of course, the
Hecke-algebra relations are
\begin{alignat*}{5}
  g_j g_k &= g_k g_j&&\text{for}& |j - k| &\geq 2&
  &\text{and}& g_j g_{j+1} g_j &= g_{j+1} g_{j} g_{j+1},\\
  h_i h_k &= h_k h_i& \quad&\text{for}&\quad |i - k| &\geq 2&\quad
  &\text{and}&\quad h_i h_{i+1} h_i &= h_{i+1} h_{i} h_{i+1}.
\end{alignat*}
The other relations are as follows:
\begin{gather}
  \label{EE^2}
  \EE \EE = \ffrac{\theta +1}{\kappa (\alpha +\beta )}\EE,
  \\
  \label{loops}
  \EE g_1 \EE = \ffrac{1}{\kappa}\EE,\qquad
  \EE h_1 \EE = \ffrac{1}{\kappa}\EE,
  \\
  \EE g_j=g_j \EE,\quad 2\leq j\leq m-1,\qquad
  \EE h_i=h_i \EE,\quad 2\leq i\leq n-1,
  \\
  \label{quintic}
  \EE  g_1  h_1^{-1} \EE  (g_1 - h_1)=0,\qquad
  (g_1 - h_1) \EE  g_1  h_1^{-1}  \EE = 0.
\end{gather}
Here,
\begin{equation*}
  \theta=\alpha\beta\kappa\kappa'.
\end{equation*}

The dimenion of $\qwB_{m,n}$ with generic parameter values is that of
the classical walled Brauer algebra, $(m+n)!$.

\subsubsection{}
Relations \eqref{EE^2}--\eqref{quintic} follow immediately from the
diagram representation of the generators.  In particular, it is
obvious that $\EE$ commutes with $g_{\geq 2}$ and $h_{\geq 2}$.  Next,
relations~\eqref{loops}, which take the form
\begin{gather*}
  \begin{tangles}{l}
    \hstr{67}\fobject{\bullet}\step[2]\fobject{\bullet}\step[2]\fobject{\circ}\\[-4.5pt]
    \vstr{67}\hstr{67}\id\step[2]\ev\\
    \vstr{67}\hstr{67}\id\step[2]\coev\\
    \vstr{67}\hstr{67}\x\step[2]\id\\
    \vstr{67}\hstr{67}\id\step[2]\ev\\
    \vstr{67}\hstr{67}\id\step[2]\coev\\[-4.5pt]
    \hstr{67}\fobject{\bullet}\step[2]\fobject{\bullet}\step[2]\fobject{\circ}
  \end{tangles}\ \
  = \ffrac{1}{\kappa}\ \ \
  \begin{tangles}{l}
    \hstr{67}\fobject{\bullet}\step[2]\fobject{\bullet}\step[2]\fobject{\circ}\\[-4.5pt]
    \vstr{150}\hstr{67}\id\step[2]\ev\\
    \vstr{150}\hstr{67}\id\step[2]\coev\\[-4.5pt]
    \hstr{67}\fobject{\bullet}\step[2]\fobject{\bullet}\step[2]\fobject{\circ}
  \end{tangles}
  \ \ ,
  \qquad\qquad
  \begin{tangles}{l}
    \hstr{67}\fobject{\bullet}\step[2]\fobject{\circ}\step[2]\fobject{\circ}\\[-4.5pt]
    \vstr{67}\hstr{67}\ev\step[2]\id\\
    \vstr{67}\hstr{67}\coev\step[2]\id\\
    \vstr{67}\hstr{67}\id\step[2]\x\\
    \vstr{67}\hstr{67}\ev\step[2]\id\\
    \vstr{67}\hstr{67}\coev\step[2]\id\\[-4.5pt]
    \hstr{67}\fobject{\bullet}\step[2]\fobject{\circ}\step[2]\fobject{\circ}
  \end{tangles}\ \
  = \ffrac{1}{\kappa}\ \ \
  \begin{tangles}{l}
    \hstr{67}\fobject{\bullet}\step[2]\fobject{\circ}\step[2]\fobject{\circ}\\[-4.5pt]
    \vstr{150}\hstr{67}\ev\step[2]\id\\
    \vstr{150}\hstr{67}\coev\step[2]\id\\[-4.5pt]
    \hstr{67}\fobject{\bullet}\step[2]\fobject{\circ}\step[2]\fobject{\circ}
  \end{tangles}\ \ ,
\end{gather*}
are corollaries of the $\kappa$-relation
in~\eqref{ribbon}--\eqref{kappas}. \ As regards
relations~\eqref{quintic}, we note that, graphically,
\begin{gather}\label{Lambda}
  \Omega\eqdef\EE  g_1  h_1^{-1} \EE = \ \ \
  \begin{tangles}{l}
    \hstr{67}\fobject{\bullet}\step[2]\fobject{\bullet}\step[2]\fobject{\circ}\step[2]\fobject{\circ}\\[-4.5pt]
    \hstr{133}\Ev\hstr{67}\ev\\
    \hstr{133}\Coev\hstr{67}\coev\\[-4.5pt]
    \hstr{67}\fobject{\bullet}\step[2]\fobject{\bullet}\step[2]\fobject{\circ}\step[2]\fobject{\circ}
  \end{tangles}
\end{gather}
whence it is obvious that $\Omega g_1$ is the same as $\Omega h_1$ up
to isotopy, and equivalently for $g_1 \Omega$ and $h_1
\Omega$.

The two relations in~\eqref{quintic} can be equivalently rewritten
with $g_1 h_1^{-1}$ replaced by $g_1^{-1}h_1$.

\subsubsection{The algebra parameters}\label{sec:parameters} The
algebra relations involve the parameters $\alpha$, $\beta$, $\kappa$,
and $\theta$, and we sometimes write
$\qwB_{m,n}(\alpha,\beta,\kappa,\theta)$ for the algebra, although two
parameters can be eliminated from the relations by rescaling the
generators.

\subsection{Reduced tangles}\label{red-rules}
\begin{figure}[tb]
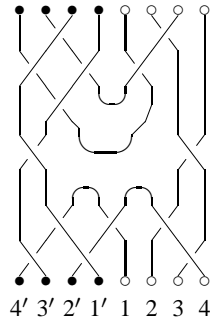

  \begin{equation*}
    \begin{tangles}{l}
      \hstr{50}\fobject{\bullet}\step[2]\fobject{\bullet}\step[2]\fobject{\bullet}\step[2]\fobject{\bullet}\step[2]\fobject{\circ}\step[2]\fobject{\circ}\step[2]\fobject{\circ}\step[2]\fobject{\circ}\\[-4pt] 
      \vstr{67}\hstr{50}\id\step[2]\xx\step[2]\id\step[2]\id\step[2]\xx\step[2]\id\\
      \vstr{67}\hstr{50}\xx\step[2]\xx\step[2]\xx\step[2]\id\step[2]\id\\
      \vstr{67}\hstr{50}\id\step[2]\xx\step[2]\ev\step[1.5]\ddh\step[2]\id\step[2]\id\\
      \vstr{25}\hstr{50}\id\step[2]\id\step[2]\hd\step[1.5]\step[2.5]\dd\step[2.5]\id\step[2]\id\\
      \vstr{66}\hx\step[1.25]\ev\step[1.75]\hx\\
      \vstr{67}\hstr{50}\id\step[2]\id\step[2]\coev\step[2]\coev\step[2]\id\step[2]\id\\
      \vstr{67}\hstr{50}\id\step[2]\x\step[2]\xx\step[2]\x\step[2]\id\\
      \vstr{67}\hstr{50}\x\step[2]\x\step[2]\id\step[2]\id\step[2]\x\\[-4pt]
      \hstr{50}\fobject{\bullet}\step[2]\fobject{\bullet}\step[2]\fobject{\bullet}\step[2]\fobject{\bullet}\step[2]\fobject{\circ}\step[2]\fobject{\circ}\step[2]\fobject{\circ}\step[2]\fobject{\circ}\\
      \hstr{50}\fobject{4'}\step[2]\fobject{3'}\step[2]\fobject{2'}\step[2]\fobject{1'}\step[2]\fobject{1}\step[2]\fobject{2}\step[2]\fobject{3}\step[2]\fobject{4}
    \end{tangles}
  \end{equation*}
  \caption{\small A reduced tangle representing an endomorphism of
    $(X^*)^{\otimes4}\otimes X^{\otimes4}$.  The tangle contains two
    top and two bottom arcs.}
  \label{F:r_g_e_tangle}
\end{figure}

Using relations \eqref{Hecke-white}--\eqref{circlebw}, each tangle can
be rewritten as a linear combination of tangles of special form, which
we call reduced tangles.  A reduced tangle is a tangle in which
\begin{enumerate}
\item no strand crosses itself,
\item every two strands cross at most once (in particular, 
      top arcs do not cross bottom arcs),
\item there are no loops,
\item\label{crossing-item} and the following preference rules hold for
  all other crossings:
  \begin{enumerate}
  \item A strand connected to a black node overcrosses any strand not
    connected to a black node.

  \item A defect or bottom arc connected to a bottom black node with a
    number $a'$ (with convention~\bref{convention}) overcrosses any
    defect or bottom arc connected to a bottom black node with a
    number $a'<b'$.

  \item Any black defect overcrosses any top arc.
    
  \item For top arcs, the same convention is applied as for the bottom
    arcs, based on the number of black ends of the arcs (also counted
    outward from the wall): the arc with a smaller number of the black
    end overcrosses.
    
%%   \item Any top arc overcrosses any white defect.

  \item For two white defects, the one connected to the $a$th bottom
    node overcrosses that connected to the $b$th bottom node if and
    only if $a>b$.
  \end{enumerate}
  \end{enumerate}
  An example of a reduced tangle is shown in
  Fig.~\ref{F:r_g_e_tangle}.\footnote{Somewhat different conventions,
    with the same effect, are used in~\cite{[W-sym]} (the terminology
    also differs in one essential point: our defects and arcs are all
    called arcs there).}

The overcrossing preferences set for reduced tangles can be
equivalently expressed in terms of a \textit{drawing order}.  Imagine
a tangle drawn in ink; lines are drawn one by one, and a new line
breaks each time it meets a line already drawn, which means a new line
undercrosses every old one.  Then, first, lines connected to the
bottom black dots $1'$, $2'$, \dots\ are drawn in this order.  Next,
top arcs are drawn in the ascending order $b'_1<b'_2<\dots$ of their
black ends. Finally, the white defects are drawn in the
\textit{descending} order $w_1>w_2>\dots$ of their bottom ends.

The reduced tangles form a basis in $\qwB_{m,n}$.

Multiplication of reduced tangles is defined in a standard way: for
two reduced tangles $T_1$ and $T_2$, the product $T_1T_2$ is the
tangle obtained by placing $T_1$ under $T_2$, identifying the
nodes on the top edge of $T_1$ with those on the bottom edge of $T_2$
according to their numbers, and then removing the intermediate nodes.
The resultant tangle is not reduced in general and should be rewritten
as a linear combination of reduced tangles using
relations~\eqref{Hecke-white}--\eqref{circlebw}.

\begin{figure}[tb]
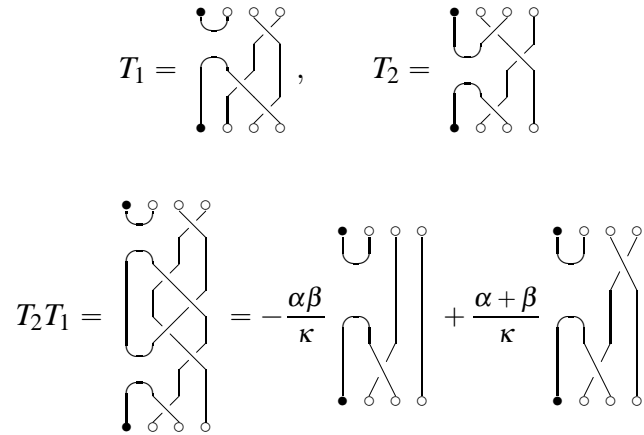

  \begin{gather*}
    T_1=\ \ \begin{tangles}{l}
      \hstr{50}\fobject{\bullet}\step[2]\fobject{\circ}\step[2]\fobject{\circ}\step[2]\fobject{\circ}\\[-4pt]
      \vstr{50}\hstr{50}\ev\step[2]\x\\
      \vstr{50}\hstr{50}\coev\step[2]\id\step[2]\id\\
      \vstr{50}\hstr{50}\id\step[2]\x\step[2]\id\\
      \vstr{50}\hstr{50}\id\step[2]\id\step[2]\x\\[-4pt]
      \hstr{50}\fobject{\bullet}\step[2]\fobject{\circ}\step[2]\fobject{\circ}\step[2]\fobject{\circ}
    \end{tangles}\ \ ,
    \qquad
    T_2=\ \ \begin{tangles}{l}
      \hstr{50}\fobject{\bullet}\step[2]\fobject{\circ}\step[2]\fobject{\circ}\step[2]\fobject{\circ}\\[-4pt]
      \vstr{50}\hstr{50}\id\step[2]\xx\step[2]\id\\
      \vstr{50}\hstr{50}\ev\step[2]\x\\
      \vstr{50}\hstr{50}\coev\step[2]\id\step[2]\id\\
      \vstr{50}\hstr{50}\id\step[2]\x\step[2]\id\\[-4pt]
      \hstr{50}\fobject{\bullet}\step[2]\fobject{\circ}\step[2]\fobject{\circ}\step[2]\fobject{\circ}
    \end{tangles}
    \\[\baselineskip]   
    T_2T_1=\ \
    \begin{tangles}{l}
      \hstr{50}\fobject{\bullet}\step[2]\fobject{\circ}\step[2]\fobject{\circ}\step[2]\fobject{\circ}\\[-4pt]
      \vstr{50}\hstr{50}\ev\step[2]\x\\
      \vstr{50}\hstr{50}\coev\step[2]\id\step[2]\id\\
      \vstr{50}\hstr{50}\id\step[2]\x\step[2]\id\\
      \vstr{50}\hstr{50}\id\step[2]\id\step[2]\x\\
\vstr{50}\hstr{50}\id\step[2]\xx\step[2]\id\\
      \vstr{50}\hstr{50}\ev\step[2]\x\\
      \vstr{50}\hstr{50}\coev\step[2]\id\step[2]\id\\
      \vstr{50}\hstr{50}\id\step[2]\x\step[2]\id\\[-4pt]
      \hstr{50}\fobject{\bullet}\step[2]\fobject{\circ}\step[2]\fobject{\circ}\step[2]\fobject{\circ}
    \end{tangles}\ \ = -\ffrac{\alpha\beta}{\kappa}\ \
    \begin{tangles}{l}
      \hstr{50}\fobject{\bullet}\step[2]\fobject{\circ}\step[2]\fobject{\circ}\step[2]\fobject{\circ}\\[-4pt]
      \hstr{50}\ev\step[2]\id\step[2]\id\\
      \hstr{50}\coev\step[2]\id\step[2]\id\\
      \hstr{50}\id\step[2]\x\step[2]\id\\[-4pt]
      \hstr{50}\fobject{\bullet}\step[2]\fobject{\circ}\step[2]\fobject{\circ}\step[2]\fobject{\circ}
    \end{tangles}\ \
    +\ffrac{\alpha+\beta}{\kappa}\ \
    \begin{tangles}{l}
      \hstr{50}\fobject{\bullet}\step[2]\fobject{\circ}\step[2]\fobject{\circ}\step[2]\fobject{\circ}\\[-4pt]
      \hstr{50}\ev\step[2]\x\\
      \hstr{50}\coev\step[2]\id\step[2]\id\\
      \hstr{50}\id\step[2]\x\step[2]\id\\[-4pt]
      \hstr{50}\fobject{\bullet}\step[2]\fobject{\circ}\step[2]\fobject{\circ}\step[2]\fobject{\circ}
    \end{tangles}
  \end{gather*}
    \caption{\small Reduced tangles $T_1$ and $T_2$, their product as
      a nonreduced tangle, and the product expressed in terms of
      reduced tangles.}
  \label{F:T1T2}
\end{figure}

An example of calculating the product of two reduced tangles is given
in Fig.~\ref{F:T1T2}.

\subsection{Jucys--Murphy elements}\label{sec:JM}
We use the diagram representation to define a family of Jucys--Murphy
elements for $\qwB_{m,n}$: $\JMBare(n)_1=1$, $\JMBare(n)_2$, \dots,
$\JMBare(n)_{m+n}$. \ We first define $\JMBare(n)_2$, \dots,
$\JMBare(n)_{n}$ just as the Jucys--Murphy elements for the ``white''
Hecke algebra $\Hecke_{n}(\alpha,\beta)$ (see~\bref{JM-Hecke}):
\begin{equation}\label{JM-white}
  \JMBare(n)_1=1,\qquad
  \JMBare(n)_{i}= (-\alpha
  \beta)^{-1}h_{n+1-i}\JMBare(n)_{i-1}h_{n+1-i},
  \quad i=2,\dots,n.
\end{equation}
This is equivalent to the following tangle definition of
$\JMBare(n)_{i}$, up to a factor: the $i$th strand from the right
double-braids with all the strands to the right of it.  In particular,
\begin{gather}\label{Jnn}
  \JMBare(n)_{n} = (-\alpha\beta)^{-n+1}\ \
  \begin{tangles}{l}
    \vstr{50}\hstr{50}\fobject{\bullet}\step[2]\fobject{\bullet}\step[2]\fobject{\bullet}\step[2]\fobject{\bullet}\step[2]\fobject{\circ}\step[2]\fobject{\circ}\step[2]\fobject{\circ}\\[-4.5pt]
%%     \vstr{62}\hstr{50}\id\step[2]\id\step[2]\x\step[2]\id\step[2]\id\\
    \vstr{80}\hstr{50}\id\step[2]\id\step[2]\id\step[2]\id\step[2]\x\step[2]\id\\
    \vstr{80}\hstr{50}\id\step[2]\id\step[2]\id\step[2]\id\step[2]\id\step[2]\x\\
    \vstr{80}\hstr{50}\id\step[2]\id\step[2]\id\step[2]\id\step[2]\id\step[2]\x\\
    \vstr{80}\hstr{50}\id\step[2]\id\step[2]\id\step[2]\id\step[2]\x\step[2]\id\\[-4.5pt]
    \hstr{50}\fobject{\bullet}\step[2]\fobject{\bullet}\step[2]\fobject{\bullet}\step[2]\fobject{\bullet}\step[2]\fobject{\circ}\step[2]\fobject{\circ}\step[2]\fobject{\circ}
  \end{tangles}
\end{gather}
(with the total number $n$ of white strands).  This tangle definition
is then extended to all the ``higher'' elements $\JMBare(n)_{n+1}$,
\dots, $\JMBare(n)_{m+n}$: each $\JMBare(n)_{n+j}$ is, up to a factor,
given by double-braiding the $(n+j)$th strand counted from the right
(just the $j$th black strand counted from the right) with all the
strands to the right of it.  In particular,
\begin{gather}\label{Jnnp1}
  \JMBare(n)_{n+1} = (-\alpha\beta)^{n}\ \ 
  \begin{tangles}{l}
    \vstr{50}\hstr{50}\fobject{\bullet}\step[2]\fobject{\bullet}\step[2]\fobject{\bullet}\step[2]\fobject{\bullet}\step[2]\fobject{\circ}\step[2]\fobject{\circ}\step[2]\fobject{\circ}\\[-4.5pt]
    \vstr{62}\hstr{50}\id\step[2]\id\step[2]\id\step[2]\x\step[2]\id\step[2]\id\\
    \vstr{62}\hstr{50}\id\step[2]\id\step[2]\id\step[2]\id\step[2]\x\step[2]\id\\
    \vstr{62}\hstr{50}\id\step[2]\id\step[2]\id\step[2]\id\step[2]\id\step[2]\x\\
    \vstr{62}\hstr{50}\id\step[2]\id\step[2]\id\step[2]\id\step[2]\id\step[2]\x\\
    \vstr{62}\hstr{50}\id\step[2]\id\step[2]\id\step[2]\id\step[2]\x\step[2]\id\\
    \vstr{62}\hstr{50}\id\step[2]\id\step[2]\id\step[2]\x\step[2]\id\step[2]\id\\[-4.5pt]
    \hstr{50}\fobject{\bullet}\step[2]\fobject{\bullet}\step[2]\fobject{\bullet}\step[2]\fobject{\bullet}\step[2]\fobject{\circ}\step[2]\fobject{\circ}\step[2]\fobject{\circ}
  \end{tangles}
\end{gather}
and
\begin{gather}\label{Jnnpmore}
  \JMBare(n)_{n+2} = (-\alpha\beta)^{n-1}\ \ 
  \begin{tangles}{l}
    \vstr{50}\hstr{50}\fobject{\bullet}\step[2]\fobject{\bullet}\step[2]\fobject{\bullet}\step[2]\fobject{\bullet}\step[2]\fobject{\circ}\step[2]\fobject{\circ}\step[2]\fobject{\circ}\\[-4.5pt]
    \vstr{62}\hstr{50}\id\step[2]\id\step[2]\x\step[2]\id\step[2]\id\step[2]\id\\
    \vstr{62}\hstr{50}\id\step[2]\id\step[2]\id\step[2]\x\step[2]\id\step[2]\id\\
    \vstr{62}\hstr{50}\id\step[2]\id\step[2]\id\step[2]\id\step[2]\x\step[2]\id\\
    \vstr{62}\hstr{50}\id\step[2]\id\step[2]\id\step[2]\id\step[2]\id\step[2]\x\\
    \vstr{62}\hstr{50}\id\step[2]\id\step[2]\id\step[2]\id\step[2]\id\step[2]\x\\
    \vstr{62}\hstr{50}\id\step[2]\id\step[2]\id\step[2]\id\step[2]\x\step[2]\id\\
    \vstr{62}\hstr{50}\id\step[2]\id\step[2]\id\step[2]\x\step[2]\id\step[2]\id\\
    \vstr{62}\hstr{50}\id\step[2]\id\step[2]\x\step[2]\id\step[2]\id\step[2]\id\\[-4.5pt]
    \hstr{50}\fobject{\bullet}\step[2]\fobject{\bullet}\step[2]\fobject{\bullet}\step[2]\fobject{\bullet}\step[2]\fobject{\circ}\step[2]\fobject{\circ}\step[2]\fobject{\circ}
  \end{tangles}\ \ ,
  \quad
  \JMBare(n)_{n+3} = (-\alpha\beta)^{n-2}\ \ 
  \begin{tangles}{l}
    \vstr{50}\hstr{50}\fobject{\bullet}\step[2]\fobject{\bullet}\step[2]\fobject{\bullet}\step[2]\fobject{\bullet}\step[2]\fobject{\circ}\step[2]\fobject{\circ}\step[2]\fobject{\circ}\\[-4.5pt]
    \vstr{50}\hstr{50}\id\step[2]\x\step[2]\id\step[2]\id\step[2]\id\step[2]\id\\
    \vstr{50}\hstr{50}\id\step[2]\id\step[2]\x\step[2]\id\step[2]\id\step[2]\id\\
    \vstr{50}\hstr{50}\id\step[2]\id\step[2]\id\step[2]\x\step[2]\id\step[2]\id\\
    \vstr{50}\hstr{50}\id\step[2]\id\step[2]\id\step[2]\id\step[2]\x\step[2]\id\\
    \vstr{50}\hstr{50}\id\step[2]\id\step[2]\id\step[2]\id\step[2]\id\step[2]\x\\
    \vstr{50}\hstr{50}\id\step[2]\id\step[2]\id\step[2]\id\step[2]\id\step[2]\x\\
    \vstr{50}\hstr{50}\id\step[2]\id\step[2]\id\step[2]\id\step[2]\x\step[2]\id\\
    \vstr{50}\hstr{50}\id\step[2]\id\step[2]\id\step[2]\x\step[2]\id\step[2]\id\\
    \vstr{50}\hstr{50}\id\step[2]\id\step[2]\x\step[2]\id\step[2]\id\step[2]\id\\
    \vstr{50}\hstr{50}\id\step[2]\x\step[2]\id\step[2]\id\step[2]\id\step[2]\id\\[-4.5pt]
    \hstr{50}\fobject{\bullet}\step[2]\fobject{\bullet}\step[2]\fobject{\bullet}\step[2]\fobject{\bullet}\step[2]\fobject{\circ}\step[2]\fobject{\circ}\step[2]\fobject{\circ}
  \end{tangles}\ \ ,
\end{gather}
and so on.  This means that
\begin{gather}\label{higherJM}
  \JMBare(n)_{n+j}=(-\alpha\beta)^{-1}\,g_{j-1} \JMBare(n)_{n+j-1}
  g_{j-1},\quad j\geq2,
\end{gather}
which allows expressing all $\JMBare(n)_{n+j}$ with $j\geq 2$ in terms
of $\JMBare(n)_{n+1}$.  For this last, crucially, we must prove
(see~\bref{JJ-lemma} below) that the definition yields an element of
$\qwB_{m,n}$; this is not obvious from~\eqref{Jnn} because the tangle
contains morphisms form the braided tensor category, \ \
$\begin{tangles}{l}
  \fobject{\bullet}\step\fobject{\circ}\\[-4.5pt]
  \vstr{67}\hx\\[-4.5pt]
  \fobject{\circ}\step\fobject{\bullet}
\end{tangles}$\ \ and
$\begin{tangles}{l}
  \fobject{\circ}\step\fobject{\bullet}\\[-4.5pt]
  \vstr{67}\hx\\[-4.5pt]
  \fobject{\bullet}\step\fobject{\circ}
\end{tangles}$\ \ , which are not $\qwB$ elements.

On the other hand, the definition immediately implies that the
$\JMBare(n)_{j}$ with $j\geq 2$ pairwise commute, simply because
(cf.~\cite{[HMR]})
\begin{gather*}
  \begin{tangles}{l}
    \vstr{20}\id\step\id\step\id\step\id\step\id\step\id\\
    \vstr{33}\hx\step\id\step\id\step\id\step\id\\
    \vstr{33}\id\step\hx\step\id\step\id\step\id\\
    \vstr{33}\id\step\id\step\hx\step\id\step\id\\
    \vstr{33}\id\step\id\step\id\step\hx\step\id\\
    \vstr{33}\id\step\id\step\id\step\id\step\hx\\
    \vstr{33}\id\step\id\step\id\step\id\step\hx\\
    \vstr{33}\id\step\id\step\id\step\hx\step\id\\
    \vstr{33}\id\step\id\step\hx\step\id\step\id\\
    \vstr{33}\id\step\hx\step\id\step\id\step\id\\
    \vstr{33}\hx\step\id\step\id\step\id\step\id\\
    \vstr{60}\id\step\id\step\id\step\id\step\id\step\id\\
    \vstr{33}\id\step\hx\step\id\step\id\step\id\\
    \vstr{33}\id\step\id\step\hx\step\id\step\id\\
    \vstr{33}\id\step\id\step\id\step\hx\step\id\\
    \vstr{33}\id\step\id\step\id\step\id\step\hx\\
    \vstr{33}\id\step\id\step\id\step\id\step\hx\\
    \vstr{33}\id\step\id\step\id\step\hx\step\id\\
    \vstr{33}\id\step\id\step\hx\step\id\step\id\\
    \vstr{33}\id\step\hx\step\id\step\id\step\id\\
    \vstr{20}\id\step\id\step\id\step\id\step\id\step\id\\
  \end{tangles}
  \ \ \ = \ \ \
  \begin{tangles}{l}
    \vstr{20}\id\step\id\step\id\step\id\step\id\step\id\\
    \vstr{33}\id\step\hx\step\id\step\id\step\id\\
    \vstr{33}\id\step\id\step\hx\step\id\step\id\\
    \vstr{33}\id\step\id\step\id\step\hx\step\id\\
    \vstr{33}\id\step\id\step\id\step\id\step\hx\\
    \vstr{33}\id\step\id\step\id\step\id\step\hx\\
    \vstr{33}\id\step\id\step\id\step\hx\step\id\\
    \vstr{33}\id\step\id\step\hx\step\id\step\id\\
    \vstr{33}\id\step\hx\step\id\step\id\step\id\\
    \vstr{60}\id\step\id\step\id\step\id\step\id\step\id\\
    \vstr{33}\hx\step\id\step\id\step\id\step\id\\
    \vstr{33}\id\step\hx\step\id\step\id\step\id\\
    \vstr{33}\id\step\id\step\hx\step\id\step\id\\
    \vstr{33}\id\step\id\step\id\step\hx\step\id\\
    \vstr{33}\id\step\id\step\id\step\id\step\hx\\
    \vstr{33}\id\step\id\step\id\step\id\step\hx\\
    \vstr{33}\id\step\id\step\id\step\hx\step\id\\
    \vstr{33}\id\step\id\step\hx\step\id\step\id\\
    \vstr{33}\id\step\hx\step\id\step\id\step\id\\
    \vstr{33}\hx\step\id\step\id\step\id\step\id\\
    \vstr{20}\id\step\id\step\id\step\id\step\id\step\id\\
  \end{tangles}
\end{gather*}
irrespective of the color of the lines, just by virtue of the braid
equation.

\begin{lemma}\label{JJ-lemma}
  All $\JMBare(n)_j$, $j=2,\dots,m+n$, are elements of $\qwB_{m,n}$.
\end{lemma}
This follows by noting that the tangle definition implies the
recursion relations
\begin{align*}%%\label{recursion-n}
  \JMBare(n)_{n+j}&=\JMBare(n-1)_{n+j-1} - \kappa(\alpha+\beta)U(j,n),\quad j\geq 2,\\
  \intertext{where}
%%   \notag
  U(j,n) &=
  (-\alpha \beta)^{n - j}
  \gdown_{j-1,1} \hdowninv_{n-1,1}\,\EE\,\hupinv_{1,n-1} \gup_{1,j-1}.
\end{align*}
We here define the ``contiguous'' products of generators
\begin{gather}\label{contiguous}
  \begin{split}
    \gup_{m,n}&=g_m g_{m+1}\dots g_n% \quad (m\leq n)
    ,
    \\
    \gdown_{m,n}&=g_m g_{m-1}\dots g_n% \quad(m \geq n)
    ,
  \end{split}
  \qquad\qquad
  \begin{split}
    \gupinv_{m,n}&=g_m^{-1} g_{m+1}^{-1}\dots g_n^{-1}% \quad(m\leq n)
    ,
    \\
    \gdowninv_{m,n}&=g_m^{-1} g_{m-1}^{-1}\dots g_n^{-1}% \quad(m \geq n)
  \end{split}
\end{gather}
(no inversion of the order of factors in the second definition), and
similarly for $h_{m,n}$.  

The recursion follows by applying~\eqref{square-of-bad} to the
rightmost double braiding in \eqref{Jnnp1}, \eqref{Jnnpmore}, and so
on.  All the $U(j,n')$ with $n'<n$ arising in applying the recursion
are elements of $\qwB_{m,n}$ by the embeddings of the ``lower'' $\qwB$
algebras.  The initial condition for the recursion,
$\JMBare(1)_{1+j}$, is an element of $\qwB_{m,n}$ as well:
by~\eqref{higherJM}, the problem reduces to $j=1$, but
$\JMBare(1)_{2}$ is just the right-hand side of~\eqref{square-of-bad}
up to a factor:
\begin{gather*}%%\label{initial-c}
  \JMBare(1)_{2} = \,\Id - (\alpha +\beta )\kappa\,\EE.
\end{gather*}
This proves~\bref{JJ-lemma}.

\subsubsection{}\label{J-explicit}
Solving the above recursion relations, we find the higher
Jucys--Murphy elements explicitly:
\begin{gather*}
  \JMBare(n)_{n + j}
  =(-\alpha \beta)^{1 - j}
    \gdown_{j - 1,1}
    \Bigl(\!1 - (\alpha+\beta)\kappa\sum_{i=1}^{n}(-\alpha\beta)^{i-1}
    \hdowninv_{i-1,1}\,\EE\,\hupinv_{1,i-1}
    \Bigr)
    \gup_{1,j - 1}.
\end{gather*}

\subsection{Casimir elements}
We continue using diagram representations to obtain some special
$\qwB$ elements.

\subsubsection{}\label{sec:cas1}
We define a central element $\Cas_{m,n}\in\qwB_{m,n}$, which we call a
Casimir element, as
\begin{equation*}
  \kappa(\alpha+\beta)\Cas_{m,n}= \ \ 
  \begin{tangles}{l}
    \vstr{67}\coev\\[-8pt]
    \fobject{\circ}\step[2]\fobject{\bullet}\step[1.5]\Bobject{\bullet}\step[1.2]\Bobject{\circ}\\[-1pt]
    \hstr{67}\id\step[3]\x\step[2]\id\\
    \hstr{67}\id\step[3]\id\step[2]\x\\
    \hstr{67}\id\step[3]\id\step[2]\x\\
    \hstr{67}\id\step[3]\x\step[2]\id\\[-4.5pt]
    \upfobject{\circ}\step[2]\upfobject{\bullet}\step[1.5]\Bobject{\bullet}\step[1.2]\Bobject{\circ}\\[-4pt]
    \vstr{67}\ev
  \end{tangles}\ \ \
  =\ \ \
  \begin{tangles}{l}
    \vstr{67}\step[2.67]\coev\\[-8pt]
    \Bobject{\bullet}\step[1.33]\Bobject{\circ}\step[1.33]\fobject{\circ}\step[2]\object{\bullet}\\[-1pt]
    \hstr{67}\id\step[2]\xx\step[3]\id\\
    \hstr{67}\xx\step[2]\id\step[3]\id\\
    \hstr{67}\xx\step[2]\id\step[3]\id\\
    \hstr{67}\id\step[2]\xx\step[3]\id\\[-4.5pt]
    \Bobject{\bullet}\step[1.33]\Bobject{\circ}\step[1.33]\upfobject{\circ}\step[2]\upfobject{\bullet}\\[-4pt]
    \vstr{67}\step[2.67]\ev
  \end{tangles}
\end{equation*}
where the two strands with larger endpoints respectively represent the
bunches of $m$ black and $n$ white strands.  That $\Cas_{m,n}$ is
central is obvious from this representation.  That it is an element of
$\qwB_{m,n}$ follows from the easily established recursion relations
\begin{align*}
  \Cas_{j,n} &=
  -\alpha \beta \Cas_{j - 1, n} - \ffrac{1}{\theta}
  (-\alpha \beta)^{j - n} \JMBare(n)_{n + j},\quad j\geq 1,\\
  \Cas_{0, n} &= -\ffrac{1}{\alpha \beta}
  \Cas_{0, n - 1} + \ffrac{1}{\alpha \beta}
  \hdowninv_{n - 1, 1} \hupinv_{1, n - 1},\quad n\geq 2,\\
  \Cas_{0, 1} &= \omega \cdot 1,
\end{align*}
where $\omega$ can be chosen arbitrarily (the above diagram defines
$\omega=\fffrac{1}{\alpha\beta}
-\fffrac{\theta+1}{\theta(\alpha+\beta)^2}$,
but introducing an arbitrary constant does not affect the property of
$\Cas_{m,n}$ to be central).\enlargethispage{\baselineskip}

\subsubsection{}\label{sec:cas1-more}Solving these recursion
relations, we find
\begin{align*}
  \Cas_{m, n} = \Cas_{m,n}(\omega) &= -\ffrac{1}{\theta}
  (-\alpha \beta)^{m - n}
  \sum_{j=n + 1}^{n + m}\JMBare(n)_{j} -
   \sum_{i=1}^{n - 1}(-\alpha \beta)^{m - i}
   \hupinv_{i, n - 1} \hdowninv_{n - 1, i}\\
   &\quad{}+\omega (-\alpha \beta)^{m + 1 - n}\cdot 1.\pagebreak[3]
\end{align*}

\subsubsection{}\label{sec:cas2}
Another central element is obviously given by
\begin{equation*}
  \kappa'(\alpha+\beta)\Casii_{m,n}(\omegaii)= \ \ 
  \begin{tangles}{l}
    \vstr{67}\coev\\[-8pt]
    \fobject{\bullet}\step[2]\fobject{\circ}\step[1.5]\Bobject{\bullet}\step[1.2]\Bobject{\circ}\\[-1pt]
    \hstr{67}\id\step[3]\x\step[2]\id\\
    \hstr{67}\id\step[3]\id\step[2]\x\\
    \hstr{67}\id\step[3]\id\step[2]\x\\
    \hstr{67}\id\step[3]\x\step[2]\id\\[-4.5pt]
    \upfobject{\bullet}\step[2]\upfobject{\circ}\step[1.5]\Bobject{\bullet}\step[1.2]\Bobject{\circ}\\[-4pt]
    \vstr{67}\ev
  \end{tangles}\ \ ,
\end{equation*}
which is again defined as a $\qwB$ element by appropriate recursion
relations, with the boundary condition
$\Casii_{1,0}(\omegaii)=\omegaii\cdot 1$, where $\omegaii$ can also be
chosen arbitrarily.

%% \subsection{A Temperley--Lieb subalgebra}
%% We note that $\qwB_{m,n}$ has subalegbras isomorphic to the
%% Temperley--Lieb algebra $\TL_{2\min(m,n)}$. \ A particular embedding of
%% the Temperley--Lieb algebra in $\qwB_{m,n}$ is provided by the
%% formulas
%% \begin{align*}
%%   e_{1}&=\kappa\sqrt{\ffrac{\alpha\beta}{\theta}}\EE,\\
%%   e_{2j+1}&=\kappa\sqrt{\ffrac{\alpha\beta}{\theta}}
%%             \gdowninv_{j,1}\hdown_{j,1}\EE \hupinv_{1,j}\gup_{1,j},
%%             \qquad j\geq1,\\
%%   e_{2j+2}&=\kappa\sqrt{\ffrac{\alpha\beta}{\theta}}
%%             \gdowninv_{j+1,1}\hdown_{j,1}\EE
%%             \hupinv_{1,j}\gup_{1,j+1},\quad j\geq0. 
%% \end{align*}
%% In other words, $e_{2j}=g_j^{-1} e_{2j-1} g_j$ and
%% $e_{2j+1} = h_j e_{2j} h_j^{-1}$. \ It then follows from the
%% $\qwB_{m,n}$ relations that the $e_i$ satisfy the relations
%% \begin{equation*}
%%   \begin{aligned}
%%     e_ie_j&=e_je_i,\quad |i-j|>1,\\
%%     e_ie_{i\pm1}e_i&=e_i,\\
%%     e_i^2&=\delta e_i,
%%   \end{aligned}
%%   \qquad\qquad
%%   \delta=\ffrac{\theta^{1/2}+\theta^{-1/2}}{(\alpha/\beta)^{1/2}+(\alpha/\beta)^{-1/2}}.
%% \end{equation*}

\section{Baxterization and commuting families}\label{sec:baxter}
It is well known that the relations in the Hecke algebra
$\Hecke_{m}(\alpha,\beta)$ (see Appendix~\ref{app:Hecke} for the
conventions) can be equivalently stated as the ``Yang--Baxter
equations with a spectral parameter,''
\begin{gather}\label{yb-hecke}
  g_i(w)g_{i+1}(z w) g_i(z) = g_{i+1}(z)g_{i}(z w) g_{i+1}(w),
\end{gather}
for
\begin{gather*}
  g_i(z)=g_i + \ffrac{\alpha + \beta}{z - 1}\cdot 1,\qquad z\in\oC.
\end{gather*}
We now propose a similar construction for $\qwB_{m,n}$: namely, we
``Baxterize'' the morphisms  $g= \ \begin{tangles}{l}
  \fobject{\bullet}\step[1]\fobject{\bullet}\\[-4.2pt]
  \vstr{66}\hx\\[-4.2pt]
  \fobject{\bullet}\step[1]\fobject{\bullet}
\end{tangles}$\ and
 $h=\ \begin{tangles}{l}
  \fobject{\circ}\step[1]\fobject{\circ}\\[-4.2pt]
  \vstr{66}\hx\\[-4.2pt]
  \fobject{\circ}\step[1]\fobject{\circ}
\end{tangles}$\  , as well as
\ \ $\begin{tangles}{l}
  \fobject{\bullet}\step[1]\fobject{\circ}\\[-4.2pt]
  \vstr{66}\hx\\[-4.2pt]
  \fobject{\circ}\step[1]\fobject{\bullet}
\end{tangles}$\ \ 
and
\ \ $\begin{tangles}{l}
  \fobject{\circ}\step[1]\fobject{\bullet}\\[-4.2pt]
  \vstr{66}\hx\\[-4.2pt]
  \fobject{\bullet}\step[1]\fobject{\circ}
\end{tangles}$\
, even though the last two are not elements of $\qwB$\,; but using
them allows constructing commuting families of $\qwB$ elements modeled
on conservation laws in integrable systems of statistical mechanics.

\subsection{Baxterized morphisms}
In the language of diagrams, we write the above $g_i(z)$ as
\begin{gather*}
  g(z) = \ \begin{tangles}{l}
    \fobject{\bullet}\step\fobject{\bullet}\\[-4.5pt]
    \hx\object{\raisebox{8pt}{\footnotesize$z$}}\\[-4.5pt]
    \fobject{\bullet}\step\fobject{\bullet}
  \end{tangles}\ \ \ 
  = \ \ 
  \begin{tangles}{l}
    \fobject{\bullet}\step\fobject{\bullet}\\[-4.5pt]
    \hx\\[-4.5pt]
    \fobject{\bullet}\step\fobject{\bullet}
  \end{tangles}\ \
  + \ffrac{\alpha + \beta}{z - 1}\ \ 
  \begin{tangles}{l}
    \fobject{\bullet}\step\fobject{\bullet}\\[-4.5pt]
    \id\step\id\\[-4.5pt]
    \fobject{\bullet}\step\fobject{\bullet}
  \end{tangles}\ \ .
\end{gather*}
The Yang--Baxter equation with a spectral parameter,
Eq.~\eqref{yb-hecke}, is then standardly expressed as
\begin{gather*}
  \begin{tangles}{l}
    \hstr{133}\fobject{\bullet}\step\fobject{\bullet}\step\fobject{\bullet}\\[-4.5pt]
    \hstr{133}\hx\object{\raisebox{10pt}{$_z$}}\step\id\\
    \hstr{133}\id\step\hx\object{\raisebox{10pt}{\ \ $_{z w}$}}\\
    \hstr{133}\hx\object{\raisebox{10pt}{$_w$}}\step\id\\[-4.5pt]
    \hstr{133}\fobject{\bullet}\step\fobject{\bullet}\step\fobject{\bullet}
  \end{tangles}\quad\
  = \ \
  \begin{tangles}{l}
    \hstr{133}\fobject{\bullet}\step\fobject{\bullet}\step\fobject{\bullet}\\[-4.5pt]
    \hstr{133}\id\step\hx\object{\raisebox{10pt}{$_w$}}\\
    \hstr{133}\hx\object{\raisebox{10pt}{\ \ $_{z w}$}}\step\id\\
    \hstr{133}\id\step\hx\object{\raisebox{10pt}{$_z$}}\\[-4.5pt]
    \hstr{133}\fobject{\bullet}\step\fobject{\bullet}\step\fobject{\bullet}
  \end{tangles}\ \ ,
\end{gather*}
The same of course holds for white lines, with
\begin{gather*}
  h(z) = \ \begin{tangles}{l}
    \fobject{\circ}\step\fobject{\circ}\\[-4.5pt]
    \hx\object{\raisebox{8pt}{\footnotesize$z$}}\\[-4.5pt]
    \fobject{\circ}\step\fobject{\circ}
  \end{tangles}\ \ \ 
  = \ \ 
  \begin{tangles}{l}
    \fobject{\circ}\step\fobject{\circ}\\[-4.5pt]
    \hx\\[-4.5pt]
    \fobject{\circ}\step\fobject{\circ}
  \end{tangles}\ \
  + \ffrac{\alpha + \beta}{z - 1}\ \ 
  \begin{tangles}{l}
    \fobject{\circ}\step\fobject{\circ}\\[-4.5pt]
    \id\step\id\\[-4.5pt]
    \fobject{\circ}\step\fobject{\circ}
  \end{tangles}\ \ .
\end{gather*}

In the setting of the braided category in~\bref{the-cat}, we
now extend these definitions to the ``mixed'' cases.  The following
lemma is easy to show by direct verification.
\begin{lemma}
  Let
  \begin{gather*}
    \begin{tangles}{l}
      \fobject{\bullet}\step\fobject{\circ}\\[-4.5pt]
      \hx\object{\raisebox{8pt}{\footnotesize$z$}}\\[-4.5pt]
      \fobject{\circ}\step\fobject{\bullet}
    \end{tangles}\ \
    = \ \ 
    \begin{tangles}{l}
      \fobject{\bullet}\step\fobject{\circ}\\[-4.5pt]
      \hx\\[-4.5pt]
      \fobject{\circ}\step\fobject{\bullet}
    \end{tangles}\ \
    + \ffrac{\alpha +\beta}{\alpha\beta}\ffrac{u}{z-u}\ \ 
    \begin{tangles}{l}
      \fobject{\bullet}\step[1]\fobject{\circ}\\[-4.5pt]
      \vstr{50}\hstr{50}\ev\\
      \vstr{50}\hstr{50}\coev\\[-4.5pt]
      \fobject{\circ}\step[1]\fobject{\bullet}
    \end{tangles}\ \ ,
    \qquad
    \begin{tangles}{l}
      \fobject{\circ}\step\fobject{\bullet}\\[-4.5pt]
      \hx\object{\raisebox{8pt}{\footnotesize$z$}}\\[-4.5pt]
      \fobject{\bullet}\step\fobject{\circ}
    \end{tangles}\ \
    = \ \ 
    \begin{tangles}{l}
      \fobject{\circ}\step\fobject{\bullet}\\[-4.5pt]
      \hx\\[-4.5pt]
      \fobject{\bullet}\step\fobject{\circ}
    \end{tangles}\ \
    + \ffrac{\alpha +\beta}{\alpha\beta}\ffrac{v}{z-v}\ \ 
    \begin{tangles}{l}
      \fobject{\circ}\step\fobject{\bullet}\\[-4.5pt]
      \vstr{50}\hstr{50}\ev\\
      \vstr{50}\hstr{50}\coev\\[-4.5pt]
      \fobject{\bullet}\step\fobject{\circ}
    \end{tangles}
  \end{gather*}
  with parameters $u$ and $v$.  Then all ``mixed'' Yang--Baxter
  equations with spectral parameter \textup{(}the equations involving
  black and white lines, such as
  \begin{gather*}
    \begin{tangles}{l}
      \hstr{133}\fobject{\bullet}\step\fobject{\circ}\step\fobject{\bullet}\\[-4.5pt]
      \hstr{133}\hx\object{\raisebox{10pt}{$_z$}}\step\id\\
      \hstr{133}\id\step\hx\object{\raisebox{10pt}{\ \ $_{z w}$}}\\
      \hstr{133}\hx\object{\raisebox{10pt}{$_w$}}\step\id\\[-4.5pt]
      \hstr{133}\fobject{\bullet}\step\fobject{\circ}\step\fobject{\bullet}
    \end{tangles}\quad\
    = \ \
    \begin{tangles}{l}
      \hstr{133}\fobject{\bullet}\step\fobject{\circ}\step\fobject{\bullet}\\[-4.5pt]
      \hstr{133}\id\step\hx\object{\raisebox{10pt}{$_w$}}\\
      \hstr{133}\hx\object{\raisebox{10pt}{\ \ $_{z w}$}}\step\id\\
      \hstr{133}\id\step\hx\object{\raisebox{10pt}{$_z$}}\\[-4.5pt]
      \hstr{133}\fobject{\bullet}\step\fobject{\circ}\step\fobject{\bullet}
    \end{tangles}
  \end{gather*}
  and others\textup{)} hold if
  \begin{gather}\label{uv}
    u v = -\theta.
  \end{gather}
\end{lemma}
We assume the last equation to hold in what follows.

\subsubsection{}\label{identity-zw}We note that 
\begin{gather*}
  \begin{tangles}{l}
    \hstr{133}\fobject{\bullet}\step\fobject{\circ}\\[-4.5pt]
    \hstr{133}\hx\object{\raisebox{10pt}{$_z$}}\\
    \hstr{133}\hx\object{\raisebox{10pt}{$_w$}}\\[-4.5pt]
    \hstr{133}\fobject{\bullet}\step\fobject{\circ}
  \end{tangles}\ \
  = -\ffrac{1}{\alpha\beta}\ \ \
  \begin{tangles}{l}
    \hstr{133}\fobject{\bullet}\step\fobject{\circ}\\[-4.5pt]
    \hstr{133}\vstr{150}\id\step\id\\[-4.5pt]
    \hstr{133}\fobject{\bullet}\step\fobject{\circ}
  \end{tangles}\ \
  + \ffrac{\kappa(\alpha +\beta)}{\alpha\beta}
  \ffrac{z w - 1}{(z-u)(w-v)}\ \
  \begin{tangles}{l}
    \hstr{133}\fobject{\bullet}\step\fobject{\circ}\\[-4.5pt]
    \vstr{67}\hstr{67}\ev\\
    \vstr{67}\hstr{67}\coev\\[-4.5pt]
    \hstr{133}\fobject{\bullet}\step\fobject{\circ}
  \end{tangles}\ \ .
\end{gather*}
Hence, first, the left-hand side is an element of
$\qwB_{m,n}(z,w)=\qwB_{m,n}\tensor\oC(z,w)$ and, second, setting
$w=1/z$ yields a pair of morphisms that are essentially inverse to
each other.

\subsection{``Universal transfer matrix'' and conservation laws}
\subsubsection{}\label{A-def}
Following~\bref{sec:cas1} in spirit, but using the Baxterized
operations introduced above, we define
\begin{equation*}
  \kappa(\alpha+\beta)\Cons_{m,n}(z,w) = \ \ 
  \begin{tangles}{l}
    \vstr{67}\coev\\[-8pt]
    \fobject{\circ}\step[2]\fobject{\bullet}\step[1.5]\Bobject{\bullet}\step[1.2]\Bobject{\circ}\\[-1pt]
    \hstr{67}\id\step[3]\x\object{\raisebox{10pt}{$_z$}}\step[2]\id\\
    \hstr{67}\id\step[3]\id\step[2]\x\object{\raisebox{10pt}{$_z$}}\\
    \hstr{67}\id\step[3]\id\step[2]\x\object{\raisebox{10pt}{$_w$}}\\
    \hstr{67}\id\step[3]\x\object{\raisebox{10pt}{$_w$}}\step[2]\id\\[-4.5pt]
    \upfobject{\circ}\step[2]\upfobject{\bullet}\step[1.5]\Bobject{\bullet}\step[1.2]\Bobject{\circ}\\[-4pt]
    \vstr{67}\ev
  \end{tangles}\ \ .
\end{equation*}
The black and white strands with larger endpoints respectively
represent $m$ black and $n$ white strands.  We emphasize that the
diagram contains \textit{three} types of Baxterized operations:
traveling from top down along the right-hand part of the loop, we
first encounter \ \ $\begin{tangles}{l}
  \fobject{\bullet}\step\fobject{\bullet}\\[-5pt]
  \vstr{67}\hx\object{\raisebox{7pt}{$_z$}}\\[-4.5pt]
  \fobject{\bullet}\step\fobject{\bullet}
  \end{tangles}$\ \
repeated $m$ times, then 
 \ \ $\begin{tangles}{l}
  \fobject{\bullet}\step\fobject{\circ}\\[-5pt]
  \vstr{67}\hx\object{\raisebox{7pt}{$_z$}}\\[-4.5pt]
  \fobject{\circ}\step\fobject{\bullet}
\end{tangles}$\ \ repeated $n$ times, and then \ \ $\begin{tangles}{l}
  \fobject{\circ}\step\fobject{\bullet}\\[-5pt]
  \vstr{67}\hx\object{\raisebox{7pt}{$_w$}}\\[-4.5pt]
  \fobject{\bullet}\step\fobject{\circ}
\end{tangles}$\
\ \ repeated also $n$ times; these are the three different types.  The
\ \ $\begin{tangles}{l}
  \fobject{\bullet}\step\fobject{\bullet}\\[-5pt]
  \vstr{67}\hx\object{\raisebox{7pt}{$_w$}}\\[-4.5pt]
  \fobject{\bullet}\step\fobject{\bullet}
\end{tangles}$\
\ \ at the bottom are of the same type as those at the top, only with
a different argument.
%% \footnote{In particular, the maximum-order pole of $\Cons_{m,n}(z,w)$
%%   is $(z-1)^{-m}(z-u)^{-n}(w-v)^{-n}(w-1)^{-m}$.}

\begin{lemma}\label{sec:recursion}
  $\Cons_{m,n}(z,w)$ is an element of $\qwB_{m,n}(z,w)$.
\end{lemma}
This follows from a system of recursion relations for the
$\Cons_{i, j}(z, w)$, which we now derive.  First, it is easy to
establish the identity
\begin{equation}%%\label{G-identity}
  \begin{tangles}{l}
    \hstr{67}\vstr{67}\coev\\[-4.5pt]
    \hstr{67}\fobject{\circ}\step[2]\fobject{\bullet}\step[1.5]\fobject{\bullet}\step[1.2]\fobject{\circ}\\[-1pt]
    \hstr{45}\vstr{67}\id\step[3]\x\object{\raisebox{7pt}{$_z$}}\step[2]\id\\
    \hstr{45}\vstr{45}\id\step[3]\id\step[2]\ev\\
    \hstr{45}\vstr{45}\id\step[3]\id\step[1.97]\coev\\[-.1pt]
    \hstr{45}\vstr{67}\id\step[3]\x\object{\raisebox{7pt}{\,$_w$}}\step[2]\id\\[-4.5pt]
    \hstr{67}\fobject{\circ}\step[2]\fobject{\bullet}\step[1.5]\fobject{\bullet}\step[1.2]\fobject{\circ}\\[-1pt]
    \hstr{67}\vstr{67}\ev
  \end{tangles}\ \ 
  =
  -\alpha \beta
  \ \ \ 
  \begin{tangles}{l}
    \hstr{80}\fobject{\bullet}\step[2]\fobject{\circ}\\[-4.5pt]
    \hstr{80}\vstr{230}\id\step[2]\id\\[-4.5pt]
    \hstr{80}\fobject{\bullet}\step[2]\fobject{\circ}
  \end{tangles}\ \
  {}+{}
  \ffrac{\alpha\beta\kappa(\alpha + \beta )}{\theta}
  \ffrac{z w + \theta}{(z-1)(w-1)}\ \
  \begin{tangles}{l}
    \hstr{80}\fobject{\bullet}\step[2]\fobject{\circ}\\[-4.5pt]
    \hstr{80}\vstr{120}\ev\\
    \hstr{80}\vstr{120}\coev\\[-4.5pt]
    \hstr{80}\fobject{\bullet}\step[2]\fobject{\circ}
  \end{tangles}\ \ .
\end{equation}
It immediately implies that
\begin{align}\label{A-rec-2}
  \Cons_{m, n}(z, w) &= -\alpha \beta \Cons_{m - 1, n}(z, w)
  \\
  \notag
  &\quad{}+ \ffrac{\alpha \beta}{\theta}
  \ffrac{z w + \theta}{(z - 1) (w - 1)}
  (-\alpha \beta)^{m - n - 1} \Jbb_{m, n}(z, w),\qquad m\geq 1,
\end{align}
where
\begin{gather}\label{Jbb}
  \Jbb_{m, n}(z, w) = \ \ \
  \begin{tangles}{l}
    \fobject{\bullet}\step[1.5]\Bobject{\bullet}\step[1.2]\Bobject{\circ}\\[-5.5pt]
    \hstr{67}\x\object{\raisebox{10pt}{$_z$}}\step[2]\id\\
    \hstr{67}\id\step[2]\x\object{\raisebox{10pt}{$_z$}}\\
    \hstr{67}\id\step[2]\x\object{\raisebox{10pt}{$_w$}}\\
    \hstr{67}\x\object{\raisebox{10pt}{$_w$}}\step[2]\id\\[-4.5pt]
    \upfobject{\bullet}\step[1.5]\Bobject{\bullet}\step[1.2]\Bobject{\circ}
  \end{tangles}
\end{gather}
with the total of $m$ black strands ($m-1$ of which are represented by
the blob) and $n$ white strands.  Assuming that the
$\Jbb_{m, n}(z, w)$ are known (to which we return momentarily), we use
relations~\eqref{A-rec-2} repeatedly until we encounter
$\Cons_{0, n}(z, w)$ in the right-hand side; then the identity
in~\bref{identity-zw} yields further relations
\begin{equation}\label{A-rec-3}
  \Cons_{0, n}(z, w) = -\ffrac{1}{\alpha \beta} \Cons_{0, n - 1}(z, w) + 
  \ffrac{1}{\alpha \beta} \ffrac{z w - 1}{(z - u) (w - v)}
  \Jww_{0, n}(z, w),
\end{equation}
where, somewhat more generally than we actually need, we define
\begin{gather}\label{Jww}
  \Jww_{m, n}(z, w) = \ \ 
  \begin{tangles}{l}
    \vstr{67}\coev\\[-8pt]
    \fobject{\circ}\step[2]\fobject{\bullet}\step[1.5]\Bobject{\bullet}\step[1.2]\Bobject{\circ}\step[1.33]\fobject{\circ}\\[-1pt]
    \hstr{67}\id\step[3]\x\object{\raisebox{10pt}{$_z$}}\step[2]\id\step[2]\id\\
    \hstr{67}\id\step[3]\id\step[2]\x\object{\raisebox{10pt}{$_z$}}\step[2]\id\\
    \vstr{67}\hstr{67}\id\step[3]\id\step[2]\id\step[2]\ev\\
    \vstr{67}\hstr{67}\id\step[3]\id\step[2]\id\step[2]\coev\\
    \hstr{67}\id\step[3]\id\step[2]\x\object{\raisebox{10pt}{$_w$}}\step[2]\id\\
    \hstr{67}\id\step[3]\x\object{\raisebox{10pt}{$_w$}}\step[2]\id\step[2]\id\\[-4.5pt]
    \upfobject{\circ}\step[2]\upfobject{\bullet}\step[1.5]\Bobject{\bullet}\step[1.2]\Bobject{\circ}\step[1.33]\upfobject{\circ}\\[-4pt]
    \vstr{67}\ev
  \end{tangles}\ \ 
       =  \htildeii_{n-1}(w)
  \Jww_{m, n-1}(z, w)
  \htilde_{n-1}(z),\qquad n\geq 1,
\end{gather}
where, first, the big circles represent $m$ black and $n-1$ white
strands and, second, we define
\begin{align*}
  \htilde(z)
   &\eqdef \ \ \
   \begin{tangles}{l}
     \hstr{67}\vstr{67}\coev\\[-4.5pt]
     \hstr{67}\fobject{\circ}\step[2]\fobject{\bullet}\step[2]\fobject{\circ}\step[2]\fobject{\circ}\\[-1pt]
     \hstr{67}\id\step[2]\x\object{\raisebox{10pt}{$_z$}}\step[2]\id\\[-4.5pt]
     \hstr{67}\fobject{\circ}\step[2]\fobject{\circ}\step[2]\fobject{\bullet}\step[2]\fobject{\circ}\\
     \hstr{67}\vstr{67}\step[4]\ev
   \end{tangles}\ \
   = \ \ 
   \begin{tangles}{l}
     \hstr{67}\fobject{\circ}\step[2]\fobject{\circ}\\[-4.5pt]
     \hstr{67}\xx\\[-4.5pt]
     \hstr{67}\fobject{\circ}\step[2]\fobject{\circ}
   \end{tangles}\ \
  +
  \ffrac{(\alpha +\beta) u}{\alpha\beta}
  \ffrac{1}{z-u}\ \ 
   \begin{tangles}{l}
     \hstr{67}\fobject{\circ}\step[2]\fobject{\circ}\\[-4.5pt]
     \hstr{67}\id\step[2]\id\\[-4.5pt]
     \hstr{67}\fobject{\circ}\step[2]\fobject{\circ}
   \end{tangles}\ ,
  \\[-6pt]
%%   \intertext{and}
   \htildeii(z)
   &\eqdef \ \ \
   \begin{tangles}{l}
     \hstr{67}\vstr{67}\step[4]\coev\\[-4.5pt]
     \hstr{67}\fobject{\circ}\step[2]\fobject{\circ}\step[2]\fobject{\bullet}\step[2]\fobject{\circ}\\     
     \hstr{67}\id\step[2]\x\object{\raisebox{10pt}{$_z$}}\step[2]\id\\[-4.5pt]
\hstr{67}\fobject{\circ}\step[2]\fobject{\bullet}\step[2]\fobject{\circ}\step[2]\fobject{\circ}\\[-1pt]
     \hstr{67}\vstr{67}\ev
  \end{tangles}\ \ 
  = \ \ 
   \begin{tangles}{l}
     \hstr{67}\fobject{\circ}\step[2]\fobject{\circ}\\[-4.5pt]
     \hstr{67}\xx\\[-4.5pt]
     \hstr{67}\fobject{\circ}\step[2]\fobject{\circ}
   \end{tangles}\ \
  +
  \ffrac{(\alpha +\beta)v}{\alpha \beta}
  \ffrac{1}{z-v}\ \ 
   \begin{tangles}{l}
     \hstr{67}\fobject{\circ}\step[2]\fobject{\circ}\\[-4.5pt]
     \hstr{67}\id\step[2]\id\\[-4.5pt]
     \hstr{67}\fobject{\circ}\step[2]\fobject{\circ}
   \end{tangles}\ \ ,
\end{align*}
with the property that
\begin{equation}\label{htildehtilde2}
  \htilde(z)\htildeii(w)=
  -\ffrac{(z\alpha + u\beta) (u\alpha + z\beta)}{(\alpha\beta)^2(z - u)^2}
  \cdot 1.
\end{equation}
We need only $\Jww_{0, n}(z, w)$ in~\eqref{A-rec-3}. Because
$\Jww_{0, 1}(z, w) = 1$, we conclude that
\begin{gather*}
  \Jww_{0, n}(z, w) =
  \bigl(\htildeii_{n-1}(w)\dots
  \htildeii_{1}(w)\bigr)
  \bigl(\htilde_{1}(z)\dots
  \htilde_{n-1}(z)\bigr),
\end{gather*}
which is an element of $\qwB_{0,n}(z,w)\subset\qwB_{m,n}(z,w)$.

We return to elements~\eqref{Jbb}, which appear in recursion
relations~\eqref{A-rec-2}. \ Clearly,
\begin{equation}\label{A-rec-4}
  \Jbb_{m, n}(z, w) = -\ffrac{1}{\alpha \beta}
  g_{m - 1}(w) \Jbb_{m - 1, n}(z, w)  g_{m - 1}(z),\qquad m\geq 2.
\end{equation}
It therefore remains to calculate $\Jbb_{1, n}(z, w)$. \ Once again by
the identity in~\bref{identity-zw}, we obtain recursion relations
\begin{equation}\label{A-rec-5}
  \Jbb_{1, n}(z, w) =
  \Jbb_{1, n - 1}(z, w)
  + \kappa(\alpha + \beta)(-\alpha \beta)^{n -1}
  \ffrac{1 - z w}{(z - u) (w - v)}
  \Uu_{1, n}(z, w),
\end{equation}
where
\begin{gather*}
  \Uu_{1, n}(z, w) =
  \htildeii_{n - 1}(w)\dots
  \htildeii_{1}(w)\,\EE\,
  \htilde_{1}(z)\dots
  \htilde_{n - 1}(z).
\end{gather*}
The initial condition is $\Jbb_{1,0}(z,w)=1$.

This finishes the proof because recursion relations~\eqref{A-rec-2},
\eqref{A-rec-3}, \eqref{A-rec-4}, and \eqref{A-rec-5} define
$\Cons_{m,n}(z,w)$ as an element of $\qwB_{m,n}(z,w)$.

\begin{rem}\label{rem:limits}
  We note the limits
  \begin{align*}
    \lim_{z,w\to\infty}\Jbb_{m,n}(z,w)&=\JMBare(n)_{m+n}
    \\
    \intertext{with a Jucys--Murphy element in the right-hand side, and}
    \lim_{z,w\to\infty}\Cons_{m,n}(z,w)&=\Cas_{m,n},
  \end{align*}
  where $\Cas_{m,n}$ is defined in~\bref{sec:cas1} (and $\omega$
  chosen as indicated there).
%%   \begin{gather*}
%%     \omega=\ffrac{\alpha^2 \theta
%%       +\alpha  \beta  \theta -\alpha  \beta +\beta^2
%%       \theta}{\alpha  \beta  \theta  (\alpha +\beta)^2}.
%%   \end{gather*}
\end{rem}

\subsubsection{``Universal'' transfer matrix}For a fixed $\rho\in\oC$,
we set
\begin{equation*}
  \Cons_{m,n}(z) = \Cons_{m,n}(z,\rho z).\pagebreak[3]
\end{equation*}\pagebreak[3]%
When $\qwB_{m,n}$ acts on a particular lattice model, this object is a
transfer matrix (and $\Jbb_{m, n}(z, \rho z)$ in~\eqref{Jbb}, the
monodromy matrix). By extension, we call $\Cons_{m,n}(z)$
the transfer matrix or, to emphasize its independence from a
particular lattice model, the universal transfer matrix.

We remind the reader that Eq.~\eqref{uv} is assumed everywhere.  In
particular, $\Cons_{m,n}(z)$ depends on $u$ in addition to $\rho$ and
the parameters of the algebra; we assume all these parameters
temporarily fixed.

\begin{Thm}$\Cons_{m,n}(z)$ is a generating function for a commutative
  family of elements of $\qwB_{m,n}$:
  \begin{equation*}
    \Cons_{m,n}(z) \Cons_{m,n}(w) - \Cons_{m,n}(w) \Cons_{m,n}(z) = 0.
  \end{equation*}
\end{Thm}
The proof is by the (generally standard, but here quite lengthy) use
of the ``train argument'' \cite{[Fad]}\,---\,the Yang--Baxter equation
with the spectral parameter.

\subsection{Solving the recursion relations}
Similarly to~\bref{J-explicit}, we can solve the above recursion
relations to find a relatively explicit expression for
$\Jbb_{m, n}(z, w)$ and $\Cons_{m, n}(z, w)$.
\begin{lemma}\label{Jbb-explicit}
  We have
  \begin{multline*}
    \Jbb_{m, n}(z, w) = (-\alpha \beta)^{-m + 1}
    g_{m - 1}(w)\dots
    g_{1}(w)
    \\
    \Bigl(\!1 + \ffrac{(1 - z w) (\alpha + \beta) \kappa}{(z - u) (w - v)}
    \sum_{s=1}^{n}(-\alpha \beta)^{s - 1}
    \htildeii_{s - 1}(w)
    \dots\htildeii_{1}(w)
    \EE
    \htilde_{1}(z)
    \dots
    \htilde_{s - 1}(z)\Bigr)\\
    g_{1}(z)\dots
    g_{m - 1}(z).
  \end{multline*}
\end{lemma}
We note that
\begin{gather*}
  \Jbb_{m, n}(z, \fffrac{1}{z})=
  \ffrac{(z \alpha + \beta)^{m - 1} (\alpha + z \beta)^{m - 1}}{
    (z - 1)^{2 m - 2} (\alpha \beta)^{m - 1}}\cdot 1.
\end{gather*}
%% Also, the recursion relations in~\bref{sec:recursion} can be solved
%% similarly to how the relations in~\bref{sec:cas1} were solved:
\begin{lemma}\label{A-explicit}
  The universal transfer matrix is given by
  \begin{multline*}
    \Cons_{m, n}(z, w) = \xi_1(z,w)\cdot 1
    - \ffrac{z w + \theta}{\theta(z - 1) (w - 1)}
    (-\alpha \beta)^{m - n} \sum_{j=1}^{m}\Jbb_{j, n}(z, w)
    \\
    + \ffrac{1 - z w}{(z - u) (w - v)}
    \sum_{i=1}^{n - 1}(-\alpha \beta)^{m - i}
    \Jwww_{i, n - 1}(z, w),
  \end{multline*}
  where $\xi_1(z,w)$ is a rational function of $z$ and $w$, and
  \begin{gather*}
    \Jwww_{i, n-1}(z, w) %%%%!!!
    = \bigl(\htilde_{i}(z)\dots
    \htilde_{n-1}(z)\bigr)
    \bigl(\htildeii_{n-1}(w)
    \dots \htildeii_{i}(w)\bigr).
  \end{gather*}
\end{lemma}

\subsection{Expanding the transfer matrix}
It follows that $\Jwww_{i, n - 1}(z, -\theta/z)$ is proportional to
the identity,\footnote{By~\eqref{htildehtilde2}, \
  $\displaystyle\Jwww_{i, n - 1}(z, -\fffrac{\theta}{z}) = (-1)^{n -
    i} \fffrac{(z \alpha + u \beta)^{n - i} (u \alpha + z \beta)^{n -
      i}}{ (z - u)^{2 (n - i)} (\alpha \beta)^{2 (n - i)}}\cdot
  1$.}\pagebreak[3]
and hence the second sum in the formula for $\Cons_{m, n}(z, w)$ has
the form $\xi_2(z,w)\cdot 1 + (z w + \theta)A_2(z,w)$. \ Combining
this with the structure of $\Jbb_{m, n}(z, w)$ in~\bref{Jbb-explicit},
we conclude that
\begin{gather*}
  \Cons_{m, n}(z, w) = \xi(z,w)\cdot 1 - (z w + \theta)(1 - z w)
  A_{m,n}(z,w),
\end{gather*}
with a rational function $\xi(z,w)$ and with some
$A_{m,n}(z,w)\in\qwB_{m,n}(z,w)$ (regular at $w = 1/z$ and
$w = -\theta/z$). \ This formula suggests two natural points around
which the transfer matrix can be expanded to produce commuting
``conservation laws'' (``Hamiltonians'').\footnote{A third possibility
  is, with $w=\rho z$, to take $z\to\infty$.  Then the zeroth-degree
  term $H^{(0)}_{m,n}$ in the expansion
  $\Cons_{m, n}( z, \rho z) = H^{(0)}_{m,n} + \frac{1}{z} H^{(1)}_{m,
    n} + \dots$
  is not proportional to the identity, but is \textit{central}
  (see~\bref{rem:limits}), and the first Hamiltonian may have to be
  defined as $(H^{(0)}_{m,n})^{-1}\,H^{(1)}_{m, n}$.} We comment on
the expansion around one of these; to avoid square roots in the
formulas, it is convenient to define the transfer matrix as a function
of a single spectral parameter by setting $w=\tau^2 z$.
%% as $\Cons_{m, n}(z) =\Cons_{m, n}( z, \tau^2 z)$.
Then
\begin{equation*}
  \Cons_{m, n}( z, \tau^2 z) =
  (\dots)\cdot 1 + \bigl(z - \fffrac{1}{\tau}\bigr) 2 H^{(1)}_{m, n}(\tau)
  + \bigl(z - \fffrac{1}{\tau}\bigr)^2 H^{(2)}_{m, n}(\tau) + \dots
\end{equation*}
with a commutative family of elements $H^{(j)}_{m, n}(\tau)$,
$j\geq 1$.  It then follows from the formulas in~\bref{A-explicit}
and~\bref{Jbb-explicit} that the first Hamiltonian is
\begin{align*}
  H^{(1)}_{m, n}(\tau) &=
  - \ffrac{\tau^2}{(1 - \tau u) (\tau - v)} 
  \sum_{i=1}^{n - 1}  (-\alpha \beta)^{m - i}
  \Jwww_{i, n - 1}(\fffrac{1}{\tau}, \tau)
  \\
  &\quad{}
  -\ffrac{(\alpha + \beta) (\theta + 1)\tau^3}{\theta (\tau - 1)^2
    (1 - \tau u) (\tau - v)}
  \kappa \sum_{i=1}^{m}\sum_{j=1}^{n} (-\alpha \beta)^{m - n - i + j}
  \EE_{i, j}(\fffrac{1}{\tau}, \tau)
  \\
  &\quad{}-
  \sum_{k=2}^{m}\sum_{j=1}^{k - 1}(-1)^{j - k + 1} (-\alpha \beta)^{m
    - k + 1 - n} \ffrac{\tau^3
    (\beta + \alpha \tau)^{k - 1 - j} (\alpha + \beta
    \tau)^{k - 1 - j}}{\theta (\tau - 1)^{2 (1 - j + k)}}
  \\
  &\qquad\qquad{}\times g_{k - 1}(\tau)\dots g_{k - j+1}(\tau)
  g_{k-j}(\fffrac{1}{\tau})\dots g_{k-1}(\fffrac{1}{\tau}),
\end{align*}
where 
\begin{gather*}
  \EE_{i, j}(z,w)=
  \htildeii_{j-1}(w)\dots\htildeii_{1}(w)
  g_{i-1}(w)\dots g_{1}(w)
  \EE
  g_{1}(z)\dots g_{i-1}(z)
  \htilde_{1}(z)\dots\htilde_{j-1}(z).
\end{gather*}

The $H^{(a)}_{m, n}(\tau)$, $a\geq 1$, depend on $\tau$ and the chosen
parameter $u$, in addition to the parameters of the algebra
(see~\bref{sec:parameters}); the last formula applies in the case
where $\tau\neq 1/u$, $\tau\neq -\theta/u$, and $\tau\neq 1$.

\section{$\qwB$ Specht modules}\label{sec:specht}
We now extend the diagram presentation for $\qwB_{m,n}$ to its
modules.  For generic values of the parameters
$\alpha,\beta,\kappa,\theta$, the simple $\qwB_{m,n}$ modules are
labeled by pairs of Young diagrams $(\dgri,\dgrii)$ such that
$m-|\dgri|=n-|\dgrii|$~\cite{[Eny-cellular],[RS1]}.  We now construct
$\qwB_{m,n}$ Specht modules, which are the irreducible representations
at generic parameters, but exist for all parameter values (and duly
become reducible).\footnote{In terms of a systematically categorial
  treatment, the modules must follow by passing from the category
  $\bmcat$, which is not abelian, to its abelianization~$\abmcat$
  \cite{[Maclane]}.  The objects of $\abmcat$ are functors from
  $\bmcat$ to the category of $\oC$-vector spaces, and the morphisms
  are natural transformations between functors. The category $\bmcat$
  is then identified with a subcategory in $\abmcat$ by the Yoneda
  functor.  The monoidal structure in $\abmcat$ can be introduced
  following~\cite{[BDay]}. The category $\abmcat$ also admits an
  explicit description in terms of the $\qwB_{m,n}$ representation
  categories for all $m,n\in\oN_0$, as described in
  \cite[Ch.~10]{[Deligne]} and \cite{[CW]}.  For generic values of the
  parameters, Specht modules give all simple objects of $\abmcat$.} (A
somewhat implicit description of Specht modules as subquotients of the
regular $\qwB_{m,n}$ module was given in~\cite{[Eny-cellular]}.)

These modules are realized here as link-state representations,
somewhat analogous to those for the Temperley--Lieb algebra (see,
e.g.,~\cite{[RsA]} and the references therein).  Compared to the
Temperley--Lieb case, where the construction is in terms of
nonintersecting arcs only, we here have intersecting arcs and defect
lines, as well as Young tableaux; the Young tableaux turn out to be
``targets'' for the defect lines.  Informally, the construction of
link states for $\qwB_{m,n}$ can be summarized as follows: these are
tangles made of (bottom) arcs and of defect lines that end at Young
tableaux and serve to ``propagate'' the action of Hecke subalgebras to
the tableaux (actually, Specht modules over the Hecke
algebras);\footnote{This has evident similarities with the
  construction in~\cite{[CDvDM]}, where Specht modules of the walled
  Brauer algebra were constructed as
  $\Delta_{m,n}(\lambda',\lambda)\cong
  I_{m,n}^f\otimes_{\Sigma_{m-f,n-f}}(S^{\lambda'}\boxtimes
  S^{\lambda})$,
  where $\Sigma_{m-f,n-f}$ is the product of two symmetric groups,
  $S^{\lambda'}$ and $S^{\lambda}$ are their Specht modules, and
  $I_{m,n}^f$ is a space of ``configurations of arcs.''}  the rules
derived from those in~\bref{cat-relations} apply to disentangling arcs
from each other and defects from arcs.

\subsection{Link states}
We fix $m$, $n$, and two Young diagrams $\dgri$ and $\dgrii$ (which
can be empty),\footnote{For example, the fundamental objects $X^*$ and
  $X$ are particular cases of Specht modules:
  $X^*=\Specht{1,0,\Yboxdim{5pt}\yng(1)\,,\emptyset}$ and
  $X=\Specht{0,1,\emptyset,\Yboxdim{5pt}\yng(1)}$.} with
$f\eqdef m-|\dgri|=n-|\dgrii|\geq 0$, and describe a basis in the
corresponding Specht module $\Specht{m,n,\dgri,\dgrii}$.

\subsubsection{}\label{sec:link-rep}
Basis vectors in $\Specht{m,n,\dgri,\dgrii}$ are \textit{link states}
\begin{equation*}
  \ket{\elli\xrightarrow{\chi}\ellii,\;\tbli,\;\tblii},
\end{equation*}
where $\elli=(a'_1,\dots,a'_{f})\subset(1,\dots,m)$,
$\ellii=(a_1,\dots,a_{f})\subset(1,\dots,n)$, and $\chi$ is a
bijection, and $\tbli$ and $\tblii$ are standard Young tableaux built
on respective diagrams $\dgri$ and $\dgrii$.  It is convenient to say
that $\dgri$, $\tbli$, etc., are \textit{black}, and $\dgrii$,
$\tblii$, etc., are \textit{white}.

A link state uniquely determines a tangle of a special form as follows
(see Fig.~\ref{F:tab1_tangle} for an example):
\begin{figure}[tb]
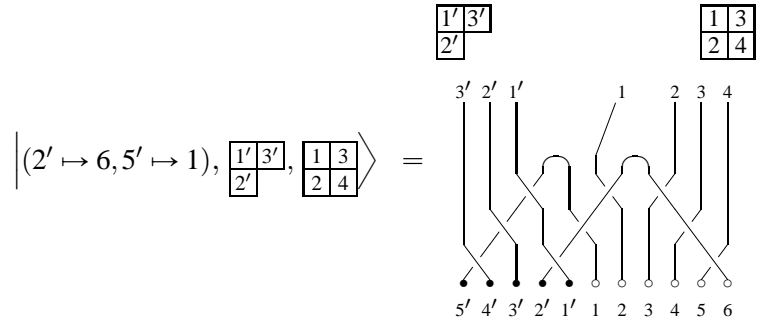
\ytableausetup{smalltableaux}
  \begin{equation*}
    \ketdemo{2'\mapsto 6,5'\mapsto 1}{\begin{ytableau}
          1' & 3'  \\
          2' 
        \end{ytableau}}
        {\begin{ytableau}
          1 & 3  \\
          2 & 4 
        \end{ytableau}}\ \ =\ \ \ \
    \begin{tangles}{l}
      \object{\protect\begin{ytableau}
          1' & 3'  \\
          2' 
        \end{ytableau}}\step[10]
      \object{\protect\begin{ytableau}
          1 & 3  \\
          2 & 4 
        \end{ytableau}}\\
%%       \hstr{50}\fobject{3'}\step[2]\fobject{2'}\step[2]\fobject{1'}\step[8]\fobject{1}\step[4]\fobject{2}\step[2]\fobject{3}\step[2]\fobject{4}\\[2pt]
%%       \vstr{67}\hstr{50}\id\step[1.75]\id\step[1.75]\id\step[6.5]\dd\step[2.5]\step[1.8]\id\step[2]\id\step[2.1]\id\\
%%       \vstr{67}\hstr{50}\id\step[1.75]\id\step[1.75]\id\step[1.9]\coev\step[2.1]\hdd\step[1.4]\coev\step[1.9]\id\step[2]\id\step[2.1]\id\\
%%       \vstr{67}\hstr{50}\id\step[1.75]\id\step[1.75]\x\step[1.9]\id\step[2]\xx\step[2]\x\step[1.9]\id\step[2.1]\id\\
%%       \vstr{67}\hstr{50}\id\step[1.75]\x\step[1.75]\id\step[1.85]\xx\step[2]\id\step[2.1]\id\step[2]\x\step[1.95]\id\\
%%       \vstr{67}\hstr{50}\x\step[1.75]\id\step[1.75]\x\step[1.8]\id\step[2.05]\id\step[2.1]\id\step[2.1]\id\step[2]\x\\[-4pt]
%%       \hstr{50}\ffobject{\bullet}\step[2]\ffobject{\bullet}\step[1.75]\ffobject{\bullet}\step[1.75]\ffobject{\bullet}\step[2]\ffobject{\bullet}\step[2]\ffobject{\circ}\step[2]\ffobject{\circ}\step[2]\ffobject{\circ}\step[2]\ffobject{\circ}\step[2]\ffobject{\circ}\step[2]\ffobject{\circ}\\
%%       \hstr{50}\ffobject{5'}\step[2]\ffobject{4'}\step[2]\ffobject{3'}\step[2]\ffobject{2'}\step[2]\ffobject{1'}\step[2]\ffobject{1}\step[2]\ffobject{2}\step[2]\ffobject{3}\step[2]\ffobject{4}\step[2]\ffobject{5}\step[2]\ffobject{6}
      \hstr{50}\ffobject{3'}\step[2]\ffobject{2'}\step[2]\ffobject{1'}\step[8]\ffobject{1}\step[4]\ffobject{2}\step[2]\ffobject{3}\step[2]\ffobject{4}\\[2pt]
      \hstr{50}\vstr{67}\id\step[2]\id\step[2]\id\step[6.5]\dd\step[2.5]\step[2]\id\step[2]\id\step[2]\id\\
      \hstr{50}\vstr{67}\id\step[2]\id\step[2]\id\step[2]\coev\step[2]\hdd\step[1.5]\coev\step[2]\id\step[2]\id\step[2]\id\\
      \hstr{50}\vstr{67}\id\step[2]\id\step[2]\x\step[2]\id\step[2]\xx\step[2]\x\step[2]\id\step[2]\id\\
      \hstr{50}\vstr{67}\id\step[2]\x\step[2]\id\step[2]\xx\step[2]\id\step[2]\id\step[2]\x\step[2]\id\\
      \hstr{50}\vstr{67}\x\step[2]\id\step[2]\x\step[2]\id\step[2]\id\step[2]\id\step[2]\id\step[2]\x\\[-4pt]
      \hstr{50}\ffobject{\bullet}\step[2]\ffobject{\bullet}\step[2]\ffobject{\bullet}\step[2]\ffobject{\bullet}\step[2]\ffobject{\bullet}\step[2]\ffobject{\circ}\step[2]\ffobject{\circ}\step[2]\ffobject{\circ}\step[2]\ffobject{\circ}\step[2]\ffobject{\circ}\step[2]\ffobject{\circ}\\
      \hstr{50}\ffobject{5'}\step[2]\ffobject{4'}\step[2]\ffobject{3'}\step[2]\ffobject{2'}\step[2]\ffobject{1'}\step[2]\ffobject{1}\step[2]\ffobject{2}\step[2]\ffobject{3}\step[2]\ffobject{4}\step[2]\ffobject{5}\step[2]\ffobject{6}
    \end{tangles}
  \end{equation*}\ytableausetup{nosmalltableaux}%
  \caption{\small A link state for $m=5$, $n=6$, and the two standard
    Young tableaux as indicated.  There are $f=2$ arcs and hence three
    black and four white defects.  Numbers at the upper ends of the
    defects show the tableau entries into which the defects are
    mapped.  The maps $(1'\mapsto 1', 3'\mapsto 2', 4'\mapsto 3')$ and
    $(2\mapsto 1, 3\mapsto 2, 4\mapsto 3, 5\mapsto 4)$ from
    bottom-edge nodes to the entries of the tableaux are monotonic,
    and hence the defects do not cross.}
  \label{F:tab1_tangle}
\end{figure}
\begin{enumerate}
  
\item The tangle has $m$ black and $n$ white nodes in the bottom row.

\item The sets $\elli$ and $\ellii$ contain the numbers of
  (respectively black and white) bottom-row nodes supporting arcs, and
  the bijection identifies pairs of nodes connected by arcs.

\item The remaining $m-f$ black nodes are mapped into the entries of
  $\tbli$ strictly monotonically (and hence bijectively); we can say
  that \textit{defect lines}, whose upper ends are associated with
  boxes of $\tbli$, do the mapping (see Fig.~\ref{F:tab1_tangle}).

  Similarly, the remaining $n-f$ white nodes are connected to defect
  lines that determine their monotonic map into the entries of
  $\tblii$.

\item the crossing preference rules (inherited from those
  in~\bref{red-rules}) apply:
  \begin{enumerate}
    
  \item black defects overcross arcs;
    
  \item an arc attached to a black node $b$ overcrosses any arc
    attached to a black node $c>b$;

  \item arcs overcross white defects; and
  \end{enumerate}
  
\item defect lines do not cross.

\end{enumerate}

\subsubsection{$\qwB_{m,n}$ action on links states}
The action of $\qwB_{m,n}$ on link states is given by attaching the
tangle corresponding to an element of $\qwB_{m,n}$ to the bottom of
the link-state tangle and forgetting the intermediate nodes.  The
resultant tangle is not necessarily a link state because it is not
reduced.  Each such nonreduced tangle can be rewritten as a linear
combination of link states using relations
\eqref{Hecke-white}--\eqref{circlebw} and the additional convention
that a free arc connecting two defects vanishes,
\begin{gather}\label{additional}
  \begin{tangles}{l}
    \fobject{m-f\ \ \ }\step[1]\object{\dots}\step[1]
    \fobject{k'}\step[1]\object{\dots}\step[1]\fobject{1'}\step[2]\fobject{1}
    \step[1]\object{\dots}\step[1]\fobject{j}\step[1]\object{\dots}\step[1]\fobject{\
      \ \ \ n-f}\\[-4.5pt]
    \hstr{200}\step[1]\Ev\\
%    \hstr{133}\Coev\hstr{67}\coev\\[-4.5pt]
%\hstr{67}\fobject{\bullet}\step[2]\fobject{\bullet}\step[2]\fobject{\circ}\step[2]\fobject{\circ}
  \end{tangles}\quad =0
   \qquad\qquad\text{(for free arcs only).}\kern-80pt
\end{gather}
A free arc is one that is not linked with any other arc. (A linked arc
can always be unlinked by using
relations~\eqref{Hecke-white}--\eqref{circlebw}.)

We now apply the above rules to describe the action with the
$\qwB_{m,n}$ generators $g_j$, $\EE$, and $h_i$ in more detail.  After
attaching the corresponding tangle in \eqref{eq:g-j}--\eqref{eq:h-i}
to the bottom of the link state, the following reduction steps are
needed to produce a linear combination of link states.

\begin{description}\addtolength{\itemsep}{6pt}
\item[Acting with $\EE$]
  \begin{enumerate}
  \item if node $1'$ is connected to an arc and node $1$ is connected
    to a defect, apply \eqref{ribbon1} to the resulting configuration
    \ \ $\begin{tangles}{l}
      \vstr{50}\hstr{50}\step[2]\coev\\
      \vstr{50}\hstr{50}\xx\step[2]\d\\[-4pt]
      \vstr{50}\hstr{50}\fobject{\bullet}\step[2]\fobject{\circ}\step[3]\fobject{\circ}\\
      \vstr{50}\hstr{50}\ev
    \end{tangles}$\;;
    
  \item if node $1'$ is connected to a defect and node $1$ is
    connected to an arc, apply \eqref{ribbon} to the resulting
    configuration \ \ $\begin{tangles}{l}
      \vstr{50}\hstr{50}\step[1]\coev\\
      \vstr{50}\hstr{50}\dd\step[2]\xx\\[-4pt]
      \vstr{50}\hstr{50}\fobject{\bullet}\step[3]\fobject{\bullet}\step[2]\fobject{\circ}\\
      \vstr{50}\hstr{50}\step[3]\ev
    \end{tangles}$\;;

  \item if both nodes $1'$ and $1$ are connected to arcs, apply
    \eqref{kappas} to the resulting configuration \ \
    $\begin{tangles}{l}
      \vstr{50}\hstr{50}\step[1]\coev\step[2]\coev\\
      \vstr{50}\hstr{50}\dd\step[2]\xx\step[2]\d\\[-4pt]
      \vstr{50}\hstr{50}\fobject{\bullet}\step[3]\fobject{\bullet}\step[2]\fobject{\circ}\step[3]\fobject{\circ}\\
      \vstr{50}\hstr{50}\step[3]\ev
    \end{tangles}$\;;

  \item if both nodes $1$ and $1'$ are connected to defects, the
    result is zero if the tangle has no arcs, and is evaluated using
    the algebra relations for disentangling linked arcs and then
    applying~\eqref{additional} to a free arc (see
    Fig.~\ref{fig:explain}).
  \end{enumerate}
  \begin{figure}[tb]
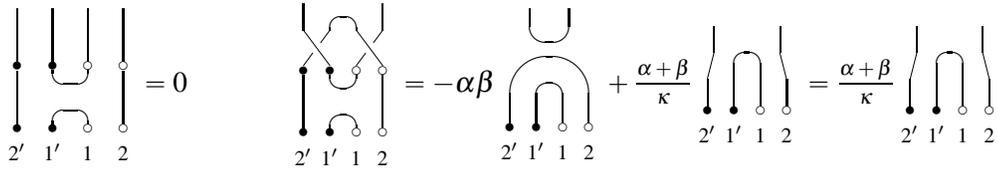

    \centering\small
    \begin{equation*}
      \begin{tangles}{l}        
        \hstr{67}\id\step[2]\id\step[2]\id\step[2]\id\\[-4pt]
        \hstr{67}\fobject{\bullet}\step[2]\fobject{\bullet}\step[2]\fobject{\circ}\step[2]\fobject{\circ}\\
        \hstr{67}\vstr{50}\id\step[2]\ev\step[2]\id\\
        \hstr{67}\vstr{50}\id\step[2]\coev\step[2]\id\\[-4pt]
        \hstr{67}\fobject{\bullet}\step[2]\fobject{\bullet}\step[2]\fobject{\circ}\step[2]\fobject{\circ}\\
        \hstr{67}\fobject{2'}\step[2]\fobject{1'}\step[2]\fobject{1}\step[2]\fobject{2}
      \end{tangles}\ \ = 0\qquad\qquad
      \begin{tangles}{l}
        \vstr{50}\hstr{50}\id\step[2]\coev\step[2]\id\\
        \vstr{67}\hstr{50}\x\step[2]\x\\[-4pt]
        \hstr{50}\fobject{\bullet}\step[2]\fobject{\bullet}\step[2]\fobject{\circ}\step[2]\fobject{\circ}\\
         \vstr{50}\hstr{50}\id\step[2]\ev\step[2]\id\\
         \vstr{50}\hstr{50}\id\step[2]\coev\step[2]\id\\[-4pt]
         \hstr{50}\fobject{\bullet}\step[2]\fobject{\bullet}\step[2]\fobject{\circ}\step[2]\fobject{\circ}\\
         \hstr{50}\fobject{2'}\step[2]\fobject{1'}\step[2]\fobject{1}\step[2]\fobject{2}
       \end{tangles}\ \ =
       -\alpha\beta\ \
       \begin{tangles}{l}
         \hstr{50}\vstr{25}\step[1.5]\id\step[3]\id\\
         \vstr{50}\hstr{50}\step[1.5]\Ev\step[2]\\[8pt]
         \vstr{50}\hstr{50}\COEV\coev\\
         \vstr{50}\hstr{50}\id\step[2]\id\step[2]\id\step[2]\id\\[-4pt]
         \hstr{50}\fobject{\bullet}\step[2]\fobject{\bullet}\step[2]\fobject{\circ}\step[2]\fobject{\circ}\\
         \hstr{50}\fobject{2'}\step[2]\fobject{1'}\step[2]\fobject{1}\step[2]\fobject{2}
       \end{tangles}\ \
       +\fffrac{\alpha+\beta}{\kappa}\ \
       \begin{tangles}{l}
         \hstr{50}\ddh\step[1.5]\coev\step[1.5]\dh\\
         \vstr{50}\hstr{50}\id\step[2]\id\step[2]\id\step[2]\id\\[-4pt]
         \hstr{50}\fobject{\bullet}\step[2]\fobject{\bullet}\step[2]\fobject{\circ}\step[2]\fobject{\circ}\\
         \hstr{50}\fobject{2'}\step[2]\fobject{1'}\step[2]\fobject{1}\step[2]\fobject{2}
       \end{tangles}\ \ =
       \fffrac{\alpha+\beta}{\kappa}\ \
       \begin{tangles}{l}
         \hstr{50}\ddh\step[1.5]\coev\step[1.5]\dh\\
         \vstr{50}\hstr{50}\id\step[2]\id\step[2]\id\step[2]\id\\[-4pt]
         \hstr{50}\fobject{\bullet}\step[2]\fobject{\bullet}\step[2]\fobject{\circ}\step[2]\fobject{\circ}\\
         \hstr{50}\fobject{2'}\step[2]\fobject{1'}\step[2]\fobject{1}\step[2]\fobject{2}
       \end{tangles}
     \end{equation*}
     \caption{\small The result of applying $\EE$ to a link state with
       no arcs is zero, as shown on the left, because in this case the
       ``upper arc'' of $\EE$ is free (not linked with any arc).  But
       whenever the original link state contains arc(s), the new arc
       coming from the action of $\EE$ is necessarily linked, with an
       example shown in the middle; disentangling the arc produces a
       number of terms, in one of which a free arc connects two
       defects, which gives zero.}\label{fig:explain}
  \end{figure}

\item[Acting with $g_{j'}$]
  \begin{enumerate}\addtolength{\itemsep}{4pt}
  \item if black node $j'$ is connected to a defect and black node
    $(j+1)'$ is connected to an arc, then use \eqref{Hecke-black} in the
    resulting configuration \ \ $\begin{tangles}{l}
      \vstr{50}\hstr{50}\step[2]\coev\\
      \vstr{50}\hstr{50}\x\step[2]\d\\[-4pt]
      \vstr{50}\hstr{50}\fobject{\bullet}\step[2]\fobject{\bullet}\step[3]\fobject{\circ}\\
      \vstr{50}\hstr{50}\x
    \end{tangles}$\;.
    [If node $j'$ is connected to an arc and node $(j+1)'$ is connected
    to a defect, we already have a link state.]
      
  \item if both nodes $j'$ and $(j+1)'$ are connected to arcs and the
    arcs intersect, use \eqref{Hecke-black} in \ \ $\begin{tangles}{l}
      \vstr{50}\hstr{50}\coev\step[2]\coev\\
      \vstr{50}\hstr{50}\id\step[2]\xx\step[2]\d\\[-4pt]
      \vstr{50}\hstr{50}\fobject{\bullet}\step[2]\fobject{\bullet}\step[2]\fobject{\circ}\step[3]\fobject{\circ}\\
      \vstr{50}\hstr{50}\x \end{tangles}$\;. [If
      both nodes $j'$ and $(j+1)'$ are connected to arcs and the arcs do
      not intersect, we already have a link state.]
    
    \item if both nodes $j'$ and $(j+1)'$ are connected to defects,
      then the other ends of the defects are attached to boxes of a
      standard Young tableau $\tbli$ that contain numbers $k'$ and
      $(k+1)'$; then set
      \begin{gather*}
        g_{j'}\ket{\elli\xrightarrow{\chi}\ellii,\;\tbli,\;\tblii}
        =\ket{\elli\xrightarrow{\chi}\ellii,\;g_{k'}\acts\tbli,\;\tblii},
      \end{gather*}
      where in the right-hand side the action of $g_{k'}$ is that on a
      Specht module of $\Hecke_{m}(\alpha,\beta)$ (see
      \bref{sec:ytr}.)
  \end{enumerate}
\item[Acting with $h_i$]
  \begin{enumerate}
  \item if node $i$ is connected to a defect and node $i+1$ is
    connected to an arc, use \eqref{Hecke-white}. [If node $i$ is
    connected to an arc and node $i+1$ is connected to a defect, we
    already have a link state.]

  \item if both nodes $i$ and $i+1$ are connected to arcs and the arcs
    intersect, use \eqref{Hecke-white}. [If both nodes $i$ and $i+1$
    are connected to arcs and the arcs do not intersect, we already
    have a link state.]

  \item if both nodes $i$ and $i+1$ are connected to defects, then the
    other ends of the defects are attached to boxes containing numbers
    $k$ and $k+1$ of a standard Young tableau $\tblii$;
    then set
    \begin{gather*}
      h_{i}\ket{\elli\xrightarrow{\chi}\ellii,\;\tbli,\;\tblii}
      =\ket{\elli\xrightarrow{\chi}\ellii,\;\tbli,\;h_{k}\acts\tblii},
    \end{gather*}
    where in the right-hand side $h_k$ acts on a Specht module of
    $\Hecke_n(\alpha,\beta)$ (see \bref{sec:ytr}.)
  \end{enumerate}
\end{description}

\subsubsection{Examples}We act with the $\qwB_{5,6}$ generators (or
with their inverse, depending on which gives simpler expressions) on
the link state \ytableausetup{smalltableaux}%
$\mathscr{X}=\ket{(2'\mapsto 6,5'\mapsto 1),\
  \mbox{\footnotesize$\begin{ytableau}
      1' & 3'  \\
      2'
    \end{ytableau},
    \
    \begin{ytableau}
      1 & 3  \\
      2 & 4 
    \end{ytableau}$}}$ shown in Fig.~\ref{F:tab1_tangle}:
\begin{align*}
  g^{-1}_4\mathscr{X}&=\ketdemo{2'\to 6,4'\to 1}{\begin{ytableau}1'&3'\\ 2'\end{ytableau}}{
  \begin{ytableau}
    1 & 3 \\
    2 & 4 \\
  \end{ytableau}
  }
  \\
  \intertext{(for $g_4$ , the first case of ``Acting with $g_{j'}$''
  should be used here; equivalently, $g_4^{-1}$ ``disentangles''
  the configuration),}
  g^{-1}_3\mathscr{X}&=\ketdemo{2'\to 6,5'\to 1}{\begin{ytableau}1'&2'\\ 3'\end{ytableau}}{
  \begin{ytableau}
    1 & 3 \\
    2 & 4 \\
  \end{ytableau}
  }
  \\
  \intertext{(the action propagates to the Young tableau),}
  g_2\mathscr{X}&=\ketdemo{3'\to 6,5'\to 1}{\begin{ytableau}1'&3'\\ 2'\end{ytableau}}{
  \begin{ytableau}
    1 & 3 \\
    2 & 4 \\
  \end{ytableau}
  }
  \\
  \intertext{(a link state is obtained directly),}
  g^{-1}_1\mathscr{X}&=\ketdemo{1'\to 6,5'\to 1}{\begin{ytableau}1'&3'\\ 2'\end{ytableau}}{
  \begin{ytableau}
    1 & 3 \\
    2 & 4 \\
  \end{ytableau}
  }
  \\
  \intertext{(the same pattern as for $g_4$), and}
  \EE\mathscr{X}&=\ffrac{\alpha^2 \beta^3}{\kappa }
                  \ketdemo{1'\to 1,2'\to
                  6}{\begin{ytableau}1'&2'\\ 3'\end{ytableau}}{
  \begin{ytableau}
    1 & 3 \\
    2 & 4 \\
  \end{ytableau}
  }-\ffrac{\alpha^2 \beta^2}{\kappa }\ketdemo{1'\to
  1,2'\to 6}{\begin{ytableau}1'&3'\\ 2'\end{ytableau}}{
  \begin{ytableau}
    1 & 3 \\
    2 & 4 \\
  \end{ytableau}
  }\\
  &\quad{}+\ffrac{\alpha +\beta }{\kappa }\ketdemo{1'\to
  1,5'\to 6}{\begin{ytableau}1'&3'\\ 2'\end{ytableau}}{
  \begin{ytableau}
    1 & 3 \\
    2 & 4 \\
  \end{ytableau}
  },
\end{align*}
where it is worth visualizing the two tangle configurations involved
here:
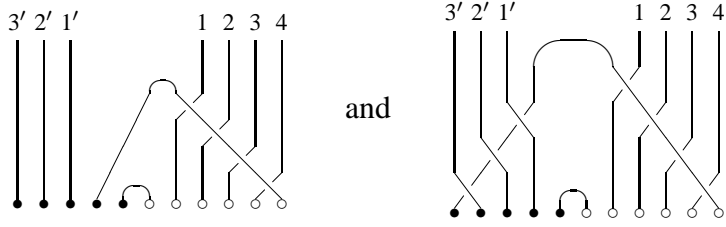
\begin{equation*}
  \begin{tangles}{l}
    \hstr{50}\fobject{3'}\step[2]\fobject{2'}\step[2]\fobject{1'}\step[10]\fobject{1}\step[2]\fobject{2}\step[2]\fobject{3}\step[2]\fobject{4}\\[2pt]
    \vstr{50}\hstr{50}\id\step[2]\id\step[2]\id\step[10]\id\step[2]\id\step[2]\id\step[2]\id\\
    \vstr{50}\hstr{50}\id\step[2]\id\step[2]\id\step[6]\coev\step[2]\id\step[2]\id\step[2]\id\step[2]\id\\
    \vstr{50}\hstr{50}\id\step[2]\id\step[2]\id\step[5]\dd\step[2]\x\step[2]\id\step[2]\id\step[2]\id\\
    \vstr{50}\hstr{50}\id\step[2]\id\step[2]\id\step[4]\dd\step[3]\id\step[2]\x\step[2]\id\step[2]\id\\
    \vstr{50}\hstr{50}\id\step[2]\id\step[2]\id\step[3]\dd\step[4]\id\step[2]\id\step[2]\x\step[2]\id\\
    \vstr{50}\hstr{50}\id\step[2]\id\step[2]\id\step[2]\dd\step[1]\coev\step[2]\id\step[2]\id\step[2]\id\step[2]\x\\[-4pt]
    \hstr{50}\fobject{\bullet}\step[2]\fobject{\bullet}\step[2]\fobject{\bullet}\step[2]\fobject{\bullet}\step[2]\fobject{\bullet}\step[2]\fobject{\circ}\step[2]\fobject{\circ}\step[2]\fobject{\circ}\step[2]\fobject{\circ}\step[2]\fobject{\circ}\step[2]\fobject{\circ}
  \end{tangles}
  \qquad\text{and}\qquad
  \begin{tangles}{l}
    \hstr{50}\fobject{3'}\step[2]\fobject{2'}\step[2]\fobject{1'}\step[10]\fobject{1}\step[2]\fobject{2}\step[2]\fobject{3}\step[2]\fobject{4}\\[2pt]
    \vstr{67}\id\step[1]\id\step[1]\id\step[1]\Coev\step[3]\id\step[1]\id\step[1]\id\step[1]\id\\
    \vstr{67}\hstr{50}\id\step[2]\id\step[2]\id\step[2]\id\step[6]\x\step[2]\id\step[2]\id\step[2]\id\\
    \vstr{67}\hstr{50}\id\step[2]\id\step[2]\x\step[6]\id\step[2]\x\step[2]\id\step[2]\id\\
    \vstr{67}\hstr{50}\id\step[2]\x\step[2]\id\step[6]\id\step[2]\id\step[2]\x\step[2]\id\\
    \vstr{67}\hstr{50}\x\step[2]\id\step[2]\id\step[2]\coev\step[2]\id\step[2]\id\step[2]\id\step[2]\x\\[-4pt]
    \hstr{50}\fobject{\bullet}\step[2]\fobject{\bullet}\step[2]\fobject{\bullet}\step[2]\fobject{\bullet}\step[2]\fobject{\bullet}\step[2]\fobject{\circ}\step[2]\fobject{\circ}\step[2]\fobject{\circ}\step[2]\fobject{\circ}\step[2]\fobject{\circ}\step[2]\fobject{\circ}
  \end{tangles}
\end{equation*}
Next,
\begin{align*}
  h_1\mathscr{X}&=\ketdemo{2'\to 6,5'\to 2}{\begin{ytableau}1'&3'\\ 2'\end{ytableau}}{
  \begin{ytableau}
    1 & 3 \\
    2 & 4 \\
  \end{ytableau}
  },
  \\
  h_2\mathscr{X}&=\beta  \ketdemo{2'\to 6,5'\to 1}{\begin{ytableau}1'&3'\\ 2'\end{ytableau}}{
  \begin{ytableau}
    1 & 3 \\
    2 & 4 \\
  \end{ytableau}
  }-\beta^2 \ketdemo{2'\to 6,5'\to
  1}{\begin{ytableau}1'&3'\\ 2'\end{ytableau}}{
  \begin{ytableau}
    1 & 2 \\
    3 & 4 \\
  \end{ytableau}
  },
  \\
  h^{-1}_3\mathscr{X}&=\ketdemo{2'\to 6,5'\to 1}{\begin{ytableau}1'&3'\\ 2'\end{ytableau}}{
  \begin{ytableau}
    1 & 2 \\
    3 & 4 \\
  \end{ytableau}
  },
  \\
  h_4\mathscr{X}&=\beta  \ketdemo{2'\to 6,5'\to 1}{\begin{ytableau}1'&3'\\ 2'\end{ytableau}}{
  \begin{ytableau}
    1 & 3 \\
    2 & 4 \\
  \end{ytableau}
  }-\beta^2 \ketdemo{2'\to 6,5'\to
  1}{\begin{ytableau}1'&3'\\ 2'\end{ytableau}}{
  \begin{ytableau}
    1 & 2 \\
    3 & 4 \\
  \end{ytableau}
  },
  \\
  h^{-1}_5\mathscr{X}&=\ketdemo{2'\to 5,5'\to 1}{\begin{ytableau}1'&3'\\ 2'\end{ytableau}}{
  \begin{ytableau}
    1 & 3 \\
    2 & 4 \\
  \end{ytableau}
  }.
\end{align*}%
\ytableausetup{nosmalltableaux}%

\subsection{Casimir actions}
The following lemma is an application of the link-state construction.
\begin{lemma}
  On the Specht module $\Specht{m,n,\dgri,\dgrii}$, the eigenvalue of
  the Casimir element $\Cas_{m, n}(1/(\alpha\beta))$
  \textup{(}see~\bref{sec:cas1}\textup{)} is
  \begin{gather*}
    c_{m,n}=
    -(-\alpha \beta)^{m - n}
    \Bigl(\ffrac{Z_1(\dgri)}{\theta} + Z_2(\dgrii)\Bigr),
  \end{gather*}
  where
  \begin{gather*}
    Z_1(\dgr)=\sum_{\Box\,\in\dgr}\left(-\fffrac{\beta}{\alpha}\right)^{\HookDistance(\Box)},
    \qquad
    Z_2(\dgr)=\sum_{\Box\,\in\dgr}\left(-\fffrac{\beta}{\alpha}\right)^{-\HookDistance(\Box)}
  \end{gather*}
  are the diagram contents, whose definitions differ by a minus sign
  in the exponent; $\HookDistance(\Box)$ is defined
  in~\bref{defs:HookD}.

  Similarly, the eigenvalue of $\Casii_{m, n}(1/(\alpha\beta))$
  \textup{(}see~\bref{sec:cas2}\textup{)} is
  \begin{gather*}
    \widetilde{c}_{m,n}=
    -(-\alpha\beta)^{n - m} 
    \Bigl(\ffrac{Z_1(\dgrii)}{\theta} + Z_2(\dgri)\Bigr).
  \end{gather*}
\end{lemma}

This lemma implies a necessary condition for Specht modules to belong
to the same linkage class in the characteristic degenerate cases.
\begin{lemma}
  Let $\theta=-\left(-\frac{\beta}{\alpha}\right)^r$ with
  $r\in\oZ%%\setminus\{0\}
  $.
  For a chosen $j$ such that $1\leq j\leq f=m-|\dgri|$, let $\Dgri$
  and $\Dgrii$ be Young diagrams obtained by adding $j$ boxes to the
  respective diagrams $\dgri$ and $\dgrii$.  Then a necessary
  condition for $\Specht{m,n,\dgri,\dgrii}$ and
  $\Specht{m,n,\Dgri,\Dgrii}$ to be in the same linkage class is
  that
  \begin{gather*}
    \sum_{i=1}^{j}
    \left(-\ffrac{\beta}{\alpha}\right)^{\HookDistance(\Box'_{i})-r}
    =\sum_{i=1}^{j}
    \left(-\ffrac{\beta}{\alpha}\right)^{-\HookDistance(\Box_{i})},
  \end{gather*}
  where $\Box'_{i}$ and $\Box_{i}$ are respectively the boxes of
  $\Dgri/\dgri$ and $\Dgrii/\dgrii$.
\end{lemma}

\section{Seminormal representations and the spectrum of\\
  Jucys--Murphy elements}\label{sec:seminormal}
In this section, we diagonalize the Jucys--Murphy elements of
$\qwB_{m,n}(\alpha,\beta,\kappa,\theta)$ at generic values of the
parameters.  For this, we construct seminormal $\qwB_{m,n}$
representations (seminormal bases in $\qwB_{m,n}$ Specht modules,
which coincide with the irreducible representations because the
algebra is semisimple at generic parameters).

As in the preceding section, we temporarily fix $m$ and $n$, and two
Young diagrams $\dgri$ and $\dgrii$ such that $m-|\dgri| = n -
|\dgrii|\mathbin{{=}{:}}f\geq 0$ (either diagram, or both, can be
empty).

The seminormal representation $\SN{m,n,\dgri,\dgrii}$ has a seminormal
basis of \textit{standard triples}.  These are introduced
in~\bref{sec:hypertableaux}; for each standard triple, we define its
weight in~\bref{sec:weights}, which eventually turns out to be the set
of eigenvalues of the Jucys--Murphy elements; the crucial part of the
construction is the $\qwB_{m,n}$ action on standard triples, which is
defined in~\bref{semi-action}.  The proof of
Theorem~\bref{Thm:action}, which asserts the $\qwB_{m,n}$ action, is
given in~\bref{sec:proof}.  In~\bref{sec:scalar-p}, we define an
invariant scalar product on seminormal representations.
In~\bref{Thm:diagonalize}, we finally prove the diagonalization of
Jucys--Murphy elements.

\subsection{Standard triples and mobile
  elements}\label{sec:hypertableaux}
\subsubsection{Standard triples}
Let $\widetilde{\dgr}$ be a Young diagram obtained by adding $f$ boxes
to~$\dgrii$.  A standard triple $(B,S,W)$ (essentially as
in~\cite{[SW]}, whence the term is borrowed) is a filling of
$(\dgri,\widetilde{\dgr}/\dgrii,\widetilde{\dgr})$ in accordance with
the following rules.
\begin{enumerate}
\item The disjoint union of $\dgri$ and $\widetilde{\dgr}/\dgrii$ is
  filled with $1'$, \dots, $m'$ such that $\dgri$ is made into a
  standard tableau~$B$, and $\widetilde{\dgr}/\dgrii$ into an
  antistandard skew-shape tableau $S$ (with the entries strictly
  decreasing along rows and along columns).

\item $\widetilde{\dgr}$ is filled with $1$, \dots, $n$ into a
  standard tableau $W$.
\end{enumerate}

We let $\HyperT=\HyperT_{m,n}(\dgri, \dgrii)$ be the set of all
standard triples $(B,S,W)$ constructed this way.
%% (in particular, with $n-|\dgrii|=|S|$ and $\overline{B}=\dgri$, where
%% the bar denotes the shape, without a filling).
We repeat that, for a fixed $\dgrii$, the \textit{shapes} of $S$ and
$W$ (denoted by bars) are related as $\shape{S}=\shape{W}/\dgrii$,
which we also express as $\shape{W}=\dgr\sqcup\shape{S}$.

We refer to the three tableaux in a standard triple $\hyp=(B,S,W)$ as
the \textit{black}, \textit{skew}, and \textit{white} ones for the
obvious reason of shape combined with the ``color scheme'' used
throughout this paper.

\subsubsection{Example}\label{hyper-example}
If $m=6$, $n=6$,  \ytableausetup{smalltableaux}%
$\dgri = \Ytableaui{
  {}\\
  {} }$\ , and $\dgrii = \Ytableau{ {}&{} }$\ , with $f=4$,
\ytableausetup{nosmalltableaux}%
a possible standard triple is
\begin{align*}
  (B,S,W) &=\Biggl(\
  \Ytableaui{
    1'\\
    4'
  }\ ,
  \qquad
  \Ytableaui{
    \none&\none&5'\\
    6'&3'&2'
  }\ ,
  \qquad
  \Ytableaui{
    1&2&4\\
    3&5&6
  }\ \Biggr).
  \\
\intertext{Another example is}
  (B,S,W)&=\Biggl(\
  \Ytableaui{
    2'\\
    5'
  }\ ,
  \qquad
  \Ytableaui{
    \none&\none&4'\\
    6'&3'\\
    1'
  }\ ,
  \qquad
  \Ytableaui{
    1&3&4\\
    2&5\\
    6
  }\ \Biggr).
\end{align*}

\subsubsection{Mobile elements}\label{sec:mobile}
Given a standard triple $(B,S,W)$, we write $S\superimpose W$ for $S$
superimposed on~$W$. \ If $1'$ from $S$ and the largest number in $W$
(which is $n$) then occur in the same box, we say that this box is a
\textit{mobile element} (otherwise $S\superimpose W$ contains no
mobile element).

In the first example in~\bref{hyper-example}, there is no mobile
element.  In the second example, this is the box containing $6_{1'}$
in
\begin{equation*}
  S\superimpose W = \ \ 
  \Ytableaui{
    1&3&4_{4'}\\
    2_{6'}&5_{3'}\\
    6_{1'}
  }
\end{equation*}

The fundamental property of mobile elements is that they can
travel.\footnote{A crucial property in a different
  context~\cite{[Shapiro]}.}  By moving a mobile element we mean
detaching it from the rest of $S\superimpose W$ and reattaching to the
remainder in a new position, $(S,W)\to(S_1,W_1)$; the shapes are here
related as $\shape{W_1}=\dgrii\sqcup\shape{S_1}$, and the resulting
$W_1$ is a standard tableau and $S_1$ is a skew antistandard tableau
built on a skew shape.

The \textit{orbit} of a mobile element is the set of all positions
into which it can be moved (including the original position).  The
orbit of $6_{1'}$ in the above example is shown with stars:
\begin{equation*}
  \Ytableaui{
    1&3&4_{4'}&\none[\Star]\\
    2_{6'}&5_{3'}&\none[\Star]\\
    \none[\Star]
  }
\end{equation*}
(thus, in addition to the above $S\superimpose W$, the reader should imagine two more, $S_1\superimpose W_1$ and $S_2\superimpose W_2$, with the $6_{1'}$ box moved to the other two positions indicated by stars).

We can speak of the orbit in terms of just the $S$ or just the $W$
part of the triple, because for any $(S_1,W_1)$ from the orbit of
$(S,W)$, \,$W_1$ is uniquely reconstructed from $S_1$, and vice versa
(of course, for a fixed $\lambda$, which is assumed).  In many cases
in what follows, the existence of a mobile element is assumed; it then
suffices to specify how the box containing $1'$ travels, and we
therefore write \hbox{$\Orb_{1'}(\dgr\sqcup S)$} for the orbit.

\subsection{Weights}\label{sec:weights} For a standard triple
$\hyp=(B,S,W)$, we define its weight $\Wt(\hyp)\in\oC^{m+n-1}$ as
follows.  The first $n-1$ components of this $(m+n-1)$-component
vector are the weight of $W$ viewed as a basis element in a seminormal
representation of the Hecke algebra $\Hecke_{n}(\alpha,\beta)$
(see~\bref{weights-f-white}):
\begin{align*}
  \Wt(B,S,W)_{i}
  &=\wt(W)_i,\qquad 1\leq i\leq n-1.\\
  \intertext{The remaining components are determined by how $(B,S)$ is
  filled with the numbers $1',\dots,m'$:}
  \Wt(B,S,W)_{n-1+i'}
  &=
    \begin{cases}
      \left(-\ffrac{\beta}{\alpha}\right)^{\HookDistance(\Pos{B}{i'})},
      &i'\in B,\\
      -\theta\left(-\ffrac{\beta}{\alpha}\right)^{-\HookDistance(\Pos{S}{i'})},
      &i'\in S,
    \end{cases}
        \qquad 1\leq i'\leq m,
\end{align*}
where, to recall, $\Pos{\tbl}{a}$ is the position (coordinates on the
plane) of a number $a$ in a Young tableau $\tbl$
(see~\bref{positions}) and $\HookDistance(\cdot)$ is defined
in~\bref{defs:HookD}.

\subsubsection{Example}\label{weight-examples} To continue with the
examples in~\bref{hyper-example}, the corresponding weights are
\begin{equation*}
  \Bigl(-\ffrac{\alpha}{\beta},\,-\ffrac{\beta}{\alpha},\,
  \ffrac{\alpha^2}{\beta^2},\,1,\,-\ffrac{\alpha}{\beta};\;
  1,\,\ffrac{\beta\theta}{\alpha},\,-\theta,\,
  -\ffrac{\beta}{\alpha},\,-\ffrac{\beta^2 \theta}{\alpha^2},\,
  \ffrac{\alpha\theta}{\beta}
  \Bigr)
\end{equation*}
and
\begin{equation*}
  \Bigl(-\ffrac{\beta }{\alpha },\,
  -\ffrac{\alpha }{\beta},\,
  \ffrac{\alpha^2}{\beta^2},\,
  1,\,
  \ffrac{\beta^2}{\alpha^2};\;
  -\ffrac{\alpha^2 \theta }{\beta^2},\,
  1,\,
  -\theta,\,
  -\ffrac{\beta^2 \theta }{\alpha^2},\,
  -\ffrac{\beta }{\alpha },\,
  \ffrac{\alpha  \theta}{\beta }\Bigr),
\end{equation*}
with the first $n-1$ components, which are $\wt(W)$, separated by a
semicolon for clarity.

\subsection{The $\qwB_{m,n}$ action on a seminormal
  basis}\label{semi-action}
The standard triples $\HyperT_{m,n}(\dgri,\dgrii)$ are a basis in the
seminormal representation $\SN{m,n,\dgri,\dgrii}$.  The action of
$\qwB_{m,n}$ generators on $\HyperT_{m,n}(\dgri,\dgrii)$ is defined
in~\bref{h-action}, \bref{g-action}, and~\bref{E-action} below.  

For $(B,S,W)\in\HyperT_{m,n}(\dgri,\dgrii)$, we let
$(\mu_1,\dots,\mu_{m+n-1})=\Wt(B,S,W)$.

\subsubsection{}\label{h-action}
For $(B,S,W)\in\HyperT_{m,n}(\dgri,\dgrii)$, we define the action of
$h_i$, $1\leq i\leq n-1$, in accordance
with~\eqref{h-acts-first}:
\begin{gather}
  h_i\acts (B,S,W) = (B,S,h_i\acts W)\nonumber
  \\
  \intertext{(with the obvious linearity assumed here and in what
    follows), where}
  \label{h-acts}
  h_{n-i}\acts W=
  \begin{cases}
    \alpha W,& \text{$i$ and $i+1$ are in the same row},\\
    \beta W,& \text{$i$ and $i+1$ are in the same column},\\
    \ffrac{\alpha + \beta}{1 - \ffrac{\mu_{i-1}}{\mu_i}}W
    +\eta\bigl(\ffrac{\mu_{i-1}}{\mu_i}\bigr)W_{(i,i+1)}& \text{otherwise}.
  \end{cases}
\end{gather}

\subsubsection{}\label{g-action}
The action of $g_j$ on a standard triple $(B,S,W)$ is
\begin{gather*}
  g_j\acts (B,S,W) = (g_j\acts(B,S), W),\quad
  1\leq j\leq m-1,
\end{gather*}\pagebreak[3]%
where
\begin{equation}\label{g-acts}
  g_j\acts(B,S)=
  \begin{cases}
    \alpha\cdot(B,S),&\text{$j$ and $j+1$ are in the same row of $B$ or $S$},\\
    \beta\cdot(B,S),&\text{$j$ and $j+1$ are in the same column of $B$ or $S$},\\
    \ffrac{\alpha + \beta}{1 - \fffrac{\mu_{n + j - 1}}{\mu_{n + j}}}
    \cdot(B,S)
    + \zeta\bigl(\ffrac{\mu_{n + j - 1}}{\mu_{n + j}}\bigr)(B,S)_{(j,j+1)},\kern-130pt&
    \kern130pt\text{otherwise},
  \end{cases}
\end{equation}
where $(B,S)_{(j,j+1)}$ is obtained from $(B,S)$ by transposing $j$
and $j+1$.  The weight components in the third line depend on $B$ and
$S$, but are independent of $W$.  Similarly to~\eqref{eta-inv}, the
function $\zeta$ involved there is such that
\begin{equation}\label{zeta-inv}
  \zeta(x)\zeta\bigl(\ffrac{1}{x}\bigr)=
 -\alpha \beta \,\theF(x)
\end{equation}
(see~\eqref{theF} for $\theF(x)$). \ We in addition require that
\begin{equation}\label{etazeta-consistency}
  \mfrac{\eta(x) \zeta(x) \eta(y) \zeta(y)}{\eta(x y) \zeta(x y)} = 
  -\alpha\beta\,\mfrac{\theF(x) \theF(y)}{\theF(x y)},
\end{equation}
a condition needed for consistency, as we see in what follows.

A ``totally symmetric'' choice satisfying all the relations for $\eta$
and $\zeta$, Eqs.~\eqref{eta-inv}, \eqref{zeta-inv},
and~\eqref{etazeta-consistency}, is
\begin{equation}\label{fully-symmetric}
  \zeta(x)=\eta(x)=\sqrt{-\alpha\beta\theF(x)}.
\end{equation}
Some formulas in what follows are essentially simplified with this
choice, but keeping general $\eta$ and $\zeta$ subject to
Eqs.~\eqref{eta-inv}, \eqref{zeta-inv},
and~\eqref{etazeta-consistency} is quite useful in the proofs, because
it ``explains'' the occurrence of different terms.\footnote{The
  general solution of~\eqref{etazeta-consistency} for
  $\eta(x) \zeta(x)$ is $\eta(x) \zeta(x)= -\alpha\beta x^Q \theF(x)$
  with any $Q$.}

%% \begin{rem}
%%   For a weight $\mu=\Wt(B,S,W)$, exactly $f=m-|\dgri|$ of the
%%   components $\mu_{k\geq n}$ (and of all the $\mu_{k\geq 1}$) contain
%%   $\theta$ as a factor.  If both $\mu_{n + j - 1}$ and $\mu_{n + j}$
%%   contain or both do not contain $\theta$, the ratio in the third line
%%   of~\eqref{g-acts} can be reexpressed in the well-known form
%%   involving the hook distance between the boxes containing $j$ and
%%   $j+1$.  If only one of the two weight components contains $\theta$,
%%   then the boxes with $j$ and $j+1$ belong to disjoint pieces in
%%   $(B, S)$.  (In that case, transposing $j$ and $j+1$ necessarily
%%   yields $(B, S)_{(j,j+1)}\mathbin{{=}{:}}(B_1,S_1)$ with $B_1$ a
%%   standard tableau and $S_1$ an antistandard filling of the skew
%%   shape.)
%% \end{rem}

\subsubsection{}\label{E-action}
The action of $\EE$ on a standard triple $(B, S, W)$ involves the
mobile element (see~\bref{sec:mobile}).  For brevity, we write
$S_1\in\Orb_{1'}(\dgr\sqcup S)$ instead of $(S_1,W_1)\in\Orb(S,W)$ for
elements in the orbit of the mobile element, with the understanding
that (for a given $\lambda$) each $S_1$ from the orbit uniquely
defines an appropriate~$W_1$. 

For a Young diagram $\dgr$, we let $\Corners{\dgr}$ denote its corners
(removable boxes):
\begin{equation*}
  \Corners{\dgr}=\{\Box\in\dgr\mid\Box\text{ is removable}\}.
\end{equation*}
The definition naturally extends to Young tableaux.  Also, if a
(skew-shaped) tableau $S$ contains $1'$, then we let $S\remove 1'$
denote the result of removing the box containing~$1'$.

Then
\begin{equation}\label{E-acts}
  \EE \acts (B, S, W) =
  \begin{cases}
    0,&\text{no mobile element in } S\superimpose W,\\
    \displaystyle
    \sum_{S_1\in\Orb_{1'}(\dgrii\sqcup S)} c_{S,S_1} \cdot(B, S_1, W_1)
    &\text{otherwise},
  \end{cases}
\end{equation}
where
\begin{gather}\nonumber
  c_{S,S_1}=\ffrac{1}{\kappa(\alpha+\beta)}\,c^{(1)}_S\,c^{(2)}_{S_1},
  \\
  \label{cSS1}
  c^{(1)}_S=
  1 + \ffrac{\theta}{\bigl(-\fffrac{\beta}{\alpha}\bigr)^{\HookDistance(\Pos{S}{1'})}},
  \qquad
  c^{(2)}_{S}=
  \frac{\displaystyle
    \prod_{\delta\in\corners{\dgr\sqcup S\,\remove\,1'\,}}
    \Bigl(\!1-\bigl(-\fffrac{\beta}{\alpha}\bigr)^{
      \HookDistance(\delta,\,\Pos{S}{1'})}\Bigr)}{\displaystyle
    \prod_{\substack{S_1\in\Orb_{1'}(\dgr\sqcup S)\\
        S_1\neq S}}
    \Bigl(\!1-\bigl(-\fffrac{\beta}{\alpha}\bigr)^{\HookDistance(\Pos{S_1}{1'},\, \Pos{S}{1'})}\Bigr)
  }.
\end{gather}
The coefficient $c^{(1)}_S$ depends on the position $\Pos{S}{1'}$ of
the mobile element---in fact, via $\HookDistance$, on the diagonal in
which the mobile element is located.  In the ratio in $c^{(2)}_S$, the
denominator is the product of
$1-\bigl(-\fffrac{\beta}{\alpha}\bigr)^{k}$ taken over the orbit,
where each $k$ is the hook distance from the mobile element in
$S\superimpose W$ to its position elsewhere in the orbit.  The
numerator is slightly more elaborate: for each corner $\delta$ of
$\dgr\sqcup S\remove 1'$, we calculate its hook distance $k$ from the
mobile element in $S$, and take the product of
$1-\bigl(-\fffrac{\beta}{\alpha}\bigr)^{k}$ over all such~$k$.  Apart
from the dependence on the corner containing $1'$, \ $c^{(2)}_S$
depends only on the shape $\dgr\sqcup\shape{S}$.

To simplify the formulas in what follows, we often write
\begin{equation*}
  \q=-\ffrac{\beta}{\alpha}.
\end{equation*}

\subsubsection{Examples}\label{action-examples}
We consider the mobile element marked with a star in the diagram
\begin{equation*}
  (\dgr\sqcup s,\Star)=\Ytableaui{
    {}&{}&{}&{}\\
    {}&{}\\
    {}&\Star\\
  }
\end{equation*}
(it is inessential how the blank boxes are divided between $\dgrii$
and the shape $s=\shape{S}$).  Then the positions of the mobile
element in the orbit, shown with stars, and the corners $\delta=\Circ$
are
\begin{equation*}
  \Ytableaui{
    {}&{}&{}&\Circ&\none[\Star]\\
    {}&\Circ&\none[\Star]\\
    \Circ&\none[\Star]\\
    \none[\Star]
  }
\end{equation*}
For each of the four positions of the mobile element attached to the
$(4,2,1)$ remainder, we construct a skew tableau $S_1$. \ The
respective coefficients $c^{(2)}_{S_1}$ for the stars ordered from
left to right are then given by
\begin{multline*}%%%\label{four-fractions}
  \ffrac{(1 - \q^{-1}) (1 - \q^{-3}) (1 - \q^{-6})}{
    (1 - \q^{-2}) (1 - \q^{-4}) (1 - \q^{-7})},
  \quad
  \ffrac{(1 - \q) (1 - \q^{-1}) (1 - \q^{-4})}{
    (1 - \q^{2}) (1 - \q^{-2}) (1 - \q^{-5})},
  \\
  \ffrac{(1 - \q^{3}) (1 - \q) (1 - \q^{-2})}{
    (1 - \q^{4}) (1 - \q^{2}) (1 - \q^{-3})},
  \quad
  \ffrac{(1 - \q^{6}) (1 - \q^{4}) (1 - \q)}
  {(1 - \q^{7})(1 - \q^{5}) (1 - \q^{3})}.
\end{multline*}

\medskip

We next illustrate the entire formula for the $\EE$ action.  Taking
$m = 2$ and $n = 4$, we consider the $16$-dimensional seminormal
representation defined by the pair of Young diagrams
\ytableausetup{smalltableaux}%
$ (\dgr',\dgr)= \bigl(\Ytableaui{{} }, \Ytableaui{
  {}&\\
  {} }\bigr) $.\ytableausetup{nosmalltableaux}
\ An instance of the $\EE$ action on a standard triple from
$\HyperT_{2,4}(\dgri, \dgrii)$ is
\begin{multline*}
  \EE\acts\Threetableaux{
    \Ytableaui{
      2' \\
    }
  }{
    \Ytableaui{
      \none[\cdot] & \none[\cdot] & 1' \\
      \none[\cdot] \\
    }
  }{
    \Ytableaui{
      1 & 3 & 4 \\
      2 \\
    }}
  =\ffrac{1 + \q^2 \theta}{\kappa(\alpha +\beta)}
  \biggl(\ffrac{
    (1 - \q^{-3}) (1 - \q^{-1})
  }{
    (1-\q^{-4}) (1- \q^{-2})}\Threetableaux{
    \Ytableaui{
      2' \\
    }
  }{
    \Ytableaui{
      \none[\cdot] & \none[\cdot] \\
      \none[\cdot] \\
      1' \\
    }
  }{
    \Ytableaui{
      1 & 3 \\
      2 \\
      4 \\
    }
  }
  \\
  {}+\ffrac{
    (1 - \q)
    (1 - \q^3)
  }{
    (1- \q^2)
    (1-\q^4)}\Threetableaux{
    \Ytableaui{
      2' \\
    }
  }{
    \Ytableaui{
      \none[\cdot] & \none[\cdot] & 1' \\
      \none[\cdot] \\
    }
  }{
    \Ytableaui{
      1 & 3 & 4 \\
      2 \\
    }
  }
  +
  \ffrac{
    (1 - \q^{-1})
    (1 - \q)
  }{
    (1- \q^{-2}) 
    (1- \q^2)}
  \Threetableaux{
    \Ytableaui{
      2' \\
    }
  }{
    \Ytableaui{
      \none[\cdot] & \none[\cdot] \\
      \none[\cdot] & 1' \\
    }
  }{
    \Ytableaui{
      1 & 3 \\
      2 & 4 \\
    }
  }\!\!\!\biggr),
\end{multline*}
where the dots are for $\dgrii$, showing how the ``skew shape,'' which
here consists of a single box, is to be completed to a suitable Young
diagram.

\begin{Thm}\label{Thm:action}
  The formulas in \bref{h-action}, \bref{g-action},
  and~\bref{E-action} define a $\qwB_{m,n}$ representation.
\end{Thm}

\subsection{Proof}\label{sec:proof}
We need to prove the ``genuinely $\qwB$'' relations, i.e., those
involving $\EE$, Eqs.~\eqref{EE^2}--\eqref{quintic}; the relations in
the Hecke subalgebras are standard; that the $g_j$ commute with $h_i$
is also obvious.

We begin with relations~\eqref{EE^2} and~\eqref{loops}, and first
rewrite them in a suitable form showing that they depend only on a
Young diagram and a chosen corner (actually, on the diagram obtained
by removing that corner).

\subsubsection{}\label{EE-schematic}
As regards the idempotent property~\eqref{EE^2}, we calculate
\begin{equation*}
   \EE \acts \EE \acts (B, S, W) = 
   \ffrac{1}{\kappa(\alpha+\beta)}\rho_{\dgr\sqcup S} \EE \acts (B, S, W),
\end{equation*}
where, directly from the definitions,
\begin{gather*}
  \rho_{\dgr\sqcup S} =
  \kern-6pt\sum_{S_1\in\Orb_{1'}(\dgr\sqcup S)}\kern-6pt
  c^{(2)}_{S_1}c^{(1)}_{S_1}
  = \kern-6pt\sum_{S_1\in\Orb_{1'}(\dgr\sqcup S)}
  \mfrac{\displaystyle \prod_{\delta\in\corners{\dgr\sqcup
        S_1\remove 1'\,}}
    \Bigl(\!1-\q^{\HookDistance(\delta,\,\Pos{S_1}{1'})}\Bigr)}{\displaystyle
    \prod_{\substack{S_2\in\Orb_{1'}(\dgr\sqcup S_1)\\ S_2\neq
        S_1}}\kern-3pt \Bigl(\!1-\q^{\HookDistance(\Pos{S_2}{1'},\,
      \Pos{S_1}{1'})}\Bigr) } \Bigl(\!1 +
  \mfrac{\theta}{\q^{\HookDistance(\Pos{S_1}{1'})}}\Bigr).
\end{gather*}
Relation~\eqref{EE^2} is equivalent to
\begin{gather*}
  \rho_{\dgr\sqcup S} = \theta + 1,
\end{gather*}
which involves two identities---for the terms proportional to and free
of $\theta$---which are equivalent to each other and which can be
reformulated as follows.  For a given Young diagram $\Dgr$, we let
$\Cocorners{\Dgr}$ denote the (positions of) boxes addable to it.
Then \eqref{EE^2} holds for the action of $\EE$ defined
in~\bref{E-action} if and only if
\begin{gather}\label{id1}
  \sum_{\Star\,\in\,\cocorners{\Dgr}}
  \mfrac{
    \displaystyle
    \prod_{\Circ\in\corners{\Dgr\,}}\Bigl(\!1-\q^{\HookDistance(\Star,\,\Circ)}\Bigr)}{
    \displaystyle
    \prod_{\Star'\neq\Star}\Bigl(\!1-\q^{\HookDistance(\Star,\,\Star')}\Bigr)}
  =1
\end{gather}
for any Young diagram $\Dgr$.  Here, $\Circ$ ranges over corners of
the diagram and $\Star$ ranges over the boxes addable to it.

We next reduce relations \eqref{loops} to a similar identity.  We
select the relation $\EE g_1 \EE = \fffrac{1}{\kappa}\EE$ for
definiteness.  It is nontrivial in the seminormal representation only
when applied to a standard triple containing a mobile element; but the
nondiagonal part of the action of $g_1$ destroys this mobile element
($1'$ is no longer in the same box with $n$ in $S\superimpose W$);
hence, only the diagonal part of the $g_1$ action in the third line
in~\eqref{g-acts} makes a contribution.  Then the coefficient in front
of the first term in the formula for the $g_1$ action,
\begin{equation*}
  \mfrac{\alpha + \beta}{1 - \mu_{n}\mu_{n + 1}^{-1}}
  =\mfrac{\alpha + \beta}{
    1 - \q^{\HookDistance(\Pos{S}{2'},\,\Pos{S}{1'})}},
\end{equation*}
involves the hook distance between the boxes containing $1'$ and $2'$.
But $2'$ can stand only in a corner of $\dgr\sqcup S\remove 1'$;
thus, $\EE g_1 \EE\acts (B, S, W)$ turns out to be proportional to
$\EE\acts (B, S, W)$, with the coefficient similar to
$\rho_{\dgr\sqcup S}$ above, but with one of the factors
corresponding to $\Corners{\dgr\sqcup S_1\remove 1'}$ in the
numerator missing.  It thus follows that relations \eqref{loops} (each
of them, as is easy to see) hold in the seminormal representation if
and only if the identity
\begin{gather}\label{id2}
  \sum_{\Star\,\in\,\cocorners{\Dgr}}
  \mfrac{
    \displaystyle
    \prod_{\Circ\in\corners{\Dgr\,}\;,\;\; {\Circ\neq\Ast}}
    \Bigl(\!1-\q^{\HookDistance(\Star,\,\Circ)}\Bigr)}{
    \displaystyle
    \prod_{\Star'\neq\Star}\Bigl(\!1-\q^{\HookDistance(\Star,\,\Star')}\Bigr)}
  =1,
\end{gather}
holds for any Young diagram $\Dgr$ with a fixed corner
$\Ast\in\Corners{\Dgr}$.

Identities~\eqref{id1} and~\eqref{id2} are particular cases of
two-variate identities established for any Young diagram in
Appendix~\ref{app:2-variate}. \ This proves~\eqref{EE^2}
and~\eqref{loops} for the $\qwB$ action on standard triples.

\subsubsection{}
It remains to establish quintic identities~\eqref{quintic} for the
$\qwB$ action on standard triples.  For this, we calculate how the
element $\Omega=\EE g_1 h_1^{-1} \EE$ involved in these identities
acts on the seminormal basis.  First of all, $\EE\acts(B,S,W)$ is
nonzero if and only if a corner of $S\superimpose W$ is a mobile
element $\Corner{n_{1'}}$. \ Next, $h^{-1}_1$ acts depending on the
relative position of $n-1$ and $n$ in $W$ and $g_1$ acts depending on
the relative position of $2'$ and $1'$ in $S$ (it is obvious that if
$2'\in B$, then the action of $\Omega$ gives zero).  We list the
possible configurations of the relevant corners of $S\superimpose W$.

\begin{enumerate}
\item\label{c1}The superimposed tableau $S\superimpose W$ has
  $\Corner{n_{1'}}$, $\Corner{(n-1)_{?'}}$, $\Corner{?_{2'}}$ as
  (some of) its corners (where the unprimed question mark is an
  integer less than $n-1$ and the primed one is an integer greater
  than $2$);
  
\item\label{c4}boxes near a corner of $S\superimpose W$ are configured
  as \renewcommand{\tabcolsep}{0pt}%
  {\footnotesize\begin{tabular}{c|c|}
    &\;$?_{2'}$\\
    \hline
    $(n-1)_{?'}$\,&\,$n_{1'}$\\\hline
  \end{tabular}}\,; or the configuration with $(n-1)_{?'}$ and
  $?_{2'}$ transposed is realized.

\item\label{c2}$S\superimpose W$ has a row ending with
  $\Corner{(n-1)_{?'}\mid n_{1'}}$ and a corner $\Corner{?_{2'}}$;
  or a ``transposed'' case is realized: a column ends with $n_{1'}$
  just below $(n-1)_{?'}$ and another corner contains $2'$;

\item\label{c3}$S\superimpose W$ has a row ending with
  $\Corner{?_{2'}\mid n_{1'}}$ and a corner $\Corner{(n-1)_{?'}}$; or the
  case where the row is transposed into a column is realized;
  
\item\label{c5}$S\superimpose W$ has corners $\Corner{n_{1'}}$,
  $\Corner{(n-1)_{2'}}$;
  
\item\label{c6}$S\superimpose W$ has a row ending with
  $\Corner{(n-1)_{2'}\mid n_{1'}}$; or a transposed case is realized;
  
\end{enumerate}

For $i=1,\dots,6$, we say that a tableau is in class $i$ if its
configuration of relevant corners is described in item $i$ of the
above list.  We see immediately that the action of $\EE$ on a tableau
in class~\ref{c1} produces tableaux from classes \ref{c1}, \ref{c2},
and \ref{c3} and, if the original tableau has corners arranged as
\renewcommand{\tabcolsep}{0pt}%
{\footnotesize\begin{tabular}{c|c|} &\,$?_{2'}$\,\\\hline
    $(n-1)_{?'}$\\\cline{1-1}
\end{tabular}}\,,
also a tableau from class~\ref{c4}. \ We write this as
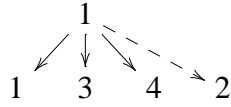
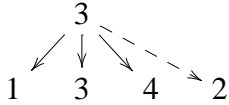
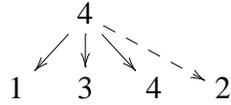
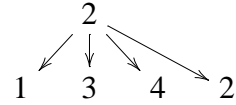
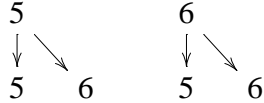
\begin{equation*}
  \xymatrix@R=12pt@C=12pt{
    &\protect\ref{c1}\ar[dl]\ar[d]\ar[dr]\ar@{{}{--}{>}}[drr]&\\
    \protect\ref{c1}&\protect\ref{c2}&\protect\ref{c3}&\protect\ref{c4}&
    }
\end{equation*}
We also have
\begin{equation*}
  \xymatrix@R=12pt@C=12pt{
    &\protect\ref{c2}\ar[dl]\ar[d]\ar[dr]\ar@{{}{--}{>}}[drr]&\\
    \protect\ref{c1}&\protect\ref{c2}&\protect\ref{c3}&\protect\ref{c4}&
  }\qquad
  \xymatrix@R=12pt@C=12pt{
    &\protect\ref{c3}\ar[dl]\ar[d]\ar[dr]\ar@{{}{--}{>}}[drr]&\\
    \protect\ref{c1}&\protect\ref{c2}&\protect\ref{c3}&\protect\ref{c4}&
  }\qquad
  \xymatrix@R=12pt@C=12pt{
    &\protect\ref{c4}\ar[dl]\ar[d]\ar[dr]\ar[drr]&\\
    \protect\ref{c1}&\protect\ref{c2}&\protect\ref{c3}&\protect\ref{c4}&
    }
\end{equation*}
and, similarly,
\begin{equation*}
  \xymatrix@R=12pt@C=12pt{
    \protect\ref{c5}\ar[d]\ar[dr]&\\
    \protect\ref{c5}&\protect\ref{c6}
  }\qquad
  \xymatrix@R=12pt@C=12pt{
    \protect\ref{c6}\ar[d]\ar[dr]&\\
    \protect\ref{c5}&\protect\ref{c6}
  }
\end{equation*}

We next see how $h^{-1}_1 g_1$ acts on the tableaux from each class.
We recall that when the relevant action is in accordance with the
third line in~\eqref{g-acts}, we have
\begin{align}\label{g-zeta}
  g_{1}\acts (B, S, W) &= \mfrac{\alpha + \beta}{1 -
    \fffrac{\mu_{n}}{\mu_{n + 1}}} (B, S, W) +
  \zeta\bigl(\ffrac{\mu_{n}}{\mu_{n + 1}}\bigr) ((B, S)_{(1',2')}, W)
  \\
  \intertext{and, similarly (from the third line in~\eqref{h-acts}),}
  \label{h-eta-inv}
  h^{-1}_{1}\acts (B, S, W)
  &=
    \ffrac{\alpha + \beta}{
    \alpha\beta\bigl(1 - \fffrac{\mu_{n-1}}{\mu_{n-2}}\bigr)}(B, S, W)
    + \ffrac{1}{\alpha \beta}\eta\bigl(\ffrac{\mu_{n-2}}{\mu_{n-1}}\bigr)
    (B, S, W_{(n-1,n)}).
\end{align}
where $(\mu_1,\,\dots,\mu_{m+n-1})=\Wt(B, S, W)$ as usual.  We use
these formulas in the case where both $1'$ and $2'$ are in $S$; then
not only the ratio $\fffrac{\mu_{n-1}}{\mu_{n-2}}$ but also
$\fffrac{\mu_{n}}{\mu_{n + 1}}$ is expressible in terms of hook
distance:
\begin{equation*}
  \ffrac{\mu_{n-1}}{\mu_{n-2}}=
  \q^{\HookDistance(n,\, n - 1)},
  \qquad
  \ffrac{\mu_{n}}{\mu_{n + 1}}=
  \q^{\HookDistance(2',\,1')}.
\end{equation*}
where we write
$\HookDistance(n, n - 1)=\HookDistance(\Pos{W}{n}, \Pos{W}{n - 1})$
and $\HookDistance(2',1')=\HookDistance(\Pos{S}{2'},\Pos{S}{1'})$ for
brevity.  

The following facts are easy to verify directly from the definitions.
\begin{enumerate}
\item In class~\ref{c1}, the action of $h^{-1}_1 g_1$, in addition to
  the original tableau, produces tableaux with $(2',1')$ or$/$and
  $(n-1,n)$ transposed; none of these has a mobile element, and hence
  all ``new'' triples vanish under the subsequent action of $\EE$.
  Thus, effectively (``inside $\Omega$''), $h^{-1}_1 g_1$ acts by the
  corresponding eigenvalue.

\item In class~\ref{c4}, both $h^{-1}_1$ and $g_1$ act literally by
  eigenvalues.

\item In class~\ref{c2}, $h^{-1}_1$ acts by an eigenvalue, and hence
  the configurations produced additionally under the action of
  $h^{-1}_1 g_1$ are those with $\Corner{(n-1)_{?'}\mid n_{2'}}$,
  $\Corner{?_{1'}}$, which are annihilated by $\EE$ (and similarly in
  the transposed case, which we do not mention explicitly any more).
  Effectively, therefore, $h^{-1}_1 g_1$ acts again by an eigenvalue.
  
\item In class~\ref{c3}, $g_1$ acts by an eigenvalue, and the
  configurations produced additionally by the action of $h^{-1}_1 g_1$
  are those with $\Corner{?_{2'}\mid(n-1)_{1'}}$, $\Corner{n_{?'}}$,
  and are also annihilated by $\EE$; hence, $h^{-1}_1 g_1$ effectively
  acts by an eigenvalue.

\item In class~\ref{c5}, the ``new'' corners resulting from the action
  of $h^{-1}_1 g_1$ are $\Corner{(n-1)_{1'}}$, $\Corner{n_{2'}}$, and
  such tableaux are annihilated by~$\EE$; but the ``old'' corners
  $\Corner{(n-1)_{2'}}$ and $\Corner{n_{1'}}$ can now occur in both
  the original and transposed positions, and we readily find
\begin{multline*}
  h^{-1}_1 g_1\Bigr|_{\ref{c5}}:(B,S,W)\mapsto
  \mfrac{(1-\q)(1-\q^{-1})}{
    (1-\q^{\HookDistance(2',\,1')})
    (1 - \q^{\HookDistance(n,\,n-1)})}\,
  (B,S,W)
  \\
  -\ffrac{1}{\alpha\beta}\eta\bigl(\q^{\HookDistance(2',\,1')}\bigr)
  \zeta\bigl(\q^{-\HookDistance(n,\,n-1)}\bigr)
  (B,S_{(1',2')},W_{(n,n-1)}).
\end{multline*}

\item In class~\ref{c6}, both $h^{-1}_1$ and $g_1$ act literally by
  eigenvalues, and
\begin{equation*}
  h^{-1}_1 g_1\Bigr|_{\ref{c6}}:(B,S,W)\mapsto(B,S,W).
\end{equation*}
\end{enumerate}
\begin{lemma}\label{lemma:1234-zero}
  We have
  \begin{align*}
    \Omega\,\Bigr|_{\ref{c1},\ref{c2},\ref{c3},\ref{c4}}=0.
  \end{align*}
  Moreover, $\Omega\acts(B,S,W)\neq 0$ if and only if
  $S\superimpose W$ has corners containing $n_{1'}$ and $(n-1)_{2'}$.
\end{lemma}

\subsubsection{Proof of~\bref{lemma:1234-zero}}
The key observation is that the eigenvalues in classes \ref{c1},
\ref{c2}, \ref{c3}, and~\ref{c4} are expressed the same:
\begin{equation}\label{omega}
  \begin{gathered}
    h^{-1}_1 g_1\Bigr|_{\ref{c1}, \ref{c2}, \ref{c3}, \ref{c4}}:
    (B,S,W)\mapsto
    \omega_{(S,W)}\,
    (B,S,W),
    \\
    \omega_{(S,W)} =
    \mfrac{(1-\q)(1-\q^{-1})}{(1-\q^{\HookDistance(2',\,1')})
      (1 - \q^{\HookDistance(n,\,n-1)})}.
  \end{gathered}
\end{equation}
(In class~\ref{c4}, this does reduce to the formula found directly,
$h^{-1}_1 g_1\Bigr|_{\ref{c4}}:(B,S,W)\mapsto
-\q^{\HookDistance(1',2')}\cdot\linebreak[0](B,S,W)$, because
$\HookDistance(1',2')=\HookDistance(n-1,n)=\pm 1$ in that case.) \ We
then calculate
\begin{align*}
  \Omega \acts (B, S, W) &= 
  \sum_{S_1\in\Orb_{1'}(\dgr\sqcup S)} c_{S,S_1}\; \EE h_1^{-1}g_1\acts(B, S_1, W_1)
  \\
  &=\sum_{S_1\in\Orb_{1'}(\dgr\sqcup S)}
    \sum_{S_2\in\Orb_{1'}(\dgr\sqcup S_1)}
    \omega_{(S_1,W_1)} c_{S,S_1}c_{S_1,S_2}\;  (B, S_2, W_2)
\end{align*}
(where, of course, $\Orb_{1'}(\dgr\sqcup S_1)=\Orb_{1'}(\dgr\sqcup
S)$) and recall the factored structure of the coefficient
in~\eqref{cSS1}. \ The right-hand side of the last formula vanishes
because
\begin{gather*}
  \sum_{S_1\in\Orb_{1'}(\dgr\sqcup S)}
  \omega_{(S_1,W_1)} c^{(2)}_{S_1}c^{(1)}_{S_1} = 0,
\end{gather*}
which is yet another identity from the class established in
Appendix~\ref{app:2-variate}; compared with~\eqref{id1}, \textit{two}
factors are missing in each denominator---those canceled by the two
factors in~\eqref{omega} involving $\HookDistance(2',1')$ and
$\HookDistance(n,n-1)$, which, up to a sign, are the distances from
the mobile element to two corners of $\dgr\sqcup S\remove 1'$.

\subsubsection{}With the vanishing of $\Omega$ on classes \ref{c1},
\ref{c2}, \ref{c3}, and \ref{c4} established in the lemma, it remains
to calculate $\Omega \acts (B, S, W)$ in the cases where
$S\superimpose W$ is in classes~\ref{c5} and~\ref{c6}, i.e., contains
boxes (corners) $n_{1'}$ and $(n-1)_{2'}$.\footnote{This is only
  possible in $\SN{m,n,\dgri,\dgrii}$ with $m-|\dgri|\geq 2$.  In
  representations with $m-|\dgri|<2$, $\Omega$ acts by zero.} We note
that in terms of the weight
$(\mu_1,\dots,\linebreak[0]\mu_{m+n-1})=\Wt(B,S,W)$, this condition is
equivalently stated as
\begin{equation*}
  \left\{
    \begin{aligned}
      \mu_{n-1}\mu_n&=-\theta,\\ \mu_{n-2}\mu_{n+1}&=-\theta.
    \end{aligned}
  \right.
\end{equation*}
In particular,
\begin{equation*}
  \ffrac{\mu_{n-1}}{\mu_{n-2}}=\ffrac{\mu_{n+1}}{\mu_n},
\end{equation*}
which we extensively use in what follows.

For a standard triple $(B,S,W)$ in class~\ref{c5} or~\ref{c6}, we now
see that $\Omega\acts(B,S,W)$ is a sum over \textit{the $1'2'$-orbit
  of $(B,S,W)$}\,---\,all tableaux obtained by detaching the boxes
containing $1'$ and $2'$ from $\dgrii\sqcup S$ and reattaching them so
as to obtain an antistandard tableau:
\begin{equation}\label{sum-z}
  \Omega\acts(B,S,W)=\!\!
  \!\!\sum_{S_1\in\Orbii(\dgr\sqcup S)}
  \ffrac{1}{\kappa^2}\,z^{\phantom{y}}_{S,S_1}(B,S_1,W_1),
\end{equation}
with the coefficients $z^{\phantom{y}}_{S,S_1}$ to be found
in~\bref{lemma:Lambda} below; here and hereafter, we let the
$1'2'$-orbit be denoted by $\Orbii(\dgr\sqcup S)$, once again with the
understanding that each $W_1$ in the sum is uniquely determined by
$\dgr\sqcup S_1$ just because the position of $n$ in $W_1$ coincides
with the position of $1'$ in $S_1$ and the position of $n-1$ coincides
with the position of $2'$.

We need more notation to proceed.

\subsubsection{Notation}
We continue operating in terms of not the superimposed tableaux
$S\superimpose W$ (now assumed to contain the boxes $n_{1'}$ and
$(n-1)_{2'}$) but the antistandard tableaux $S$ such that $1',2'\in
S$; it is understood that every rearrangement of $S$ is
followed by the corresponding rearrangement of $W$, such that $n$
travels together with~$1'$ and $(n-1)$ together with $2'$.
\begin{enumerate}
\item We let $\Corners{\dgr\sqcup S\remove1'2'}$ denote corners of
  $\dgr\sqcup S\remove1'\remove2'$, the tableau with boxes $1'$ and
  $2'$ removed.
  
\item We also let $\Corners{\dgr\sqcup S\remove1'}_{\neq2'}$ denote
  the corners of $\dgr\sqcup S\remove1'$ except the box
  containing~$2'$.  In two examples below, the elements of
  $\Corners{\dgr\sqcup S\remove1'}_{\neq2'}$ are shown in color:
  \begin{equation*}
    \Ytableauii{{}& {}& {}& *(green){}\\ {}& {}& 2'& 1'\\
      *(green){}}
    \qquad\qquad
    \Ytableauii{{}& {}& {}& {}\\ {}& {}& {}& 2'\\ {}& *(green){}& 1'}
  \end{equation*}

\item We define $S\remove 2'$, a skew antistandard tableau with the
  box $2'$ removed. The definition is obvious if $S$ contains two
  corners $\Corner{1'}$ and $\Corner{2'}$; otherwise\,---\,if $S$
  contains\ytableausetup{smalltableaux} $\Ytableau{2'&1'}$ or
  $\Ytableau{2'\\1'}$\ytableausetup{nosmalltableaux}\,---\,we let
  $S\remove 2'$ be the skew tableau with the box $2'$ removed
  \textit{and $1'$ taking its place}.
  For example,
  \begin{equation}\label{example-S}
    \dgr\sqcup S=\Ytableauii{{}& {}& {}& {}\\ {}& 2'& 1'}
    \quad\Longrightarrow\quad
    \dgr\sqcup S\remove 2'= \Ytableauii{{}& {}& {}& {}\\ {}& 1'}\;.
  \end{equation}

\item We let $\Orb_{2'}(\dgr\sqcup S\remove2')$ denote all possible
  ways to attach $2'$ to $\dgr\sqcup S\remove 2'$ so as to obtain
  an antistandard tableau.

  For $\dgr\sqcup S$ in~\eqref{example-S}, for example,
  $\Orb_{2'}(\dgr\sqcup S\remove2')$ consists of two elements marked
  with $\Ast$:
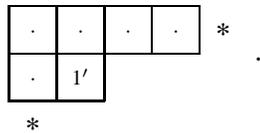
  \begin{equation*}
    \Ytableauii{\cdot& \cdot& \cdot& \cdot&\none[\Ast]\\ \cdot&
      1'
      \\\none[\Ast]}\;.
  \end{equation*}
\end{enumerate}

We also illustrate the definition of %%the $1'2'$-orbit
$\Orbii(\dgr\sqcup S)$ given just above~\eqref{sum-z}. \ For
\begin{equation*}
  \dgr\sqcup S=
  \Ytableauii{
    \none[\cdot] & \none[\cdot] & 2' \\
    6' & 5'  \\
    1'   \\
  },
\end{equation*}
its $1'2'$-orbit is this tableau itself together with five more:
\begin{equation*}
  \Ytableauii{
    \none[\cdot] & \none[\cdot] & 1' \\
    6' & 5'  \\
    2'   \\
  }\,,\qquad \Ytableauii{
    \none[\cdot] & \none[\cdot] & 2' & 1' \\
    6' & 5'   \\
  }\,,\qquad \Ytableauii{
    \none[\cdot] & \none[\cdot] & 2' \\
    6' & 5' & 1' \\
  }\,,
  \qquad
  \Ytableauii{
    \none[\cdot] & \none[\cdot] \\
    6' & 5' \\
    2' & 1' \\
  }\,,\qquad\Ytableauii{
    \none[\cdot] & \none[\cdot] \\
    6' & 5' \\
    2'  \\
    1'  \\
  }\,.
\end{equation*}

\subsubsection{}\label{sec:Lambda}
After these preparations, tracking the occurrence of the different
factors in~\eqref{sum-z} gives the following lemma, where we continue
using $\q=-\fffrac{\beta}{\alpha}$ for
brevity.
\begin{lemma}\label{lemma:Lambda}
  We have
  \begin{equation*}
    z^{\phantom{y}}_{S,T}=
    \begin{cases}
      \displaystyle
      -\ffrac{X_{S} Y_{T}}{\alpha\beta},&\Pos{S}{2'}=\Pos{T}{2'},\\
      \displaystyle
      \ffrac{X_{S} Y_{T}}{(\alpha\beta)^2}      
        \mfrac{\eta(\q^{\HookDistance(\Pos{S}{2'},\,\Pos{T}{2'})})
          \zeta(\q^{\HookDistance(\Pos{S}{2'},\,\Pos{T}{2'})})
        }{\theF(\q^{\HookDistance(\Pos{S}{2'},\,\Pos{T}{2'})})}
      &\text{otherwise},
    \end{cases}
  \end{equation*}
  where $\theF$ is defined in \eqref{theF} and the coefficients are
  \begin{gather*}
    X_{S} =
    \Bigl(\!1 + \ffrac{\theta}{
      \q^{\HookDistance(\Pos{S}{1'})}}\Bigr)
    \Bigl(\!1 + \ffrac{\theta}{
      \q^{\HookDistance(\Pos{S}{2'})}}\Bigr),
    \\
    Y_{S} = 
    \frac{C_{S}\q^{\HookDistance(\Pos{S}{2'},\, \Pos{S}{1'})}}{
      (1 - \q)(1 - \q^{-1})}\,
    \mfrac{\displaystyle
      \prod_{\delta\in\corners{\dgr\sqcup S\,\remove\,1'}_{\neq2'}}
      \Bigl(\!1-\q^{\HookDistance(\delta,\,\Pos{S}{1'})}\!\Bigr)
    }{\displaystyle
      \!\!\prod_{\substack{S_1\in\Orb_{1'}(\dgr\sqcup S)\\S_1\neq S}}\!\!
      \Bigl(\!1-\q^{\HookDistance(\Pos{S_1}{1'},\,\Pos{S}{1'})}\!\Bigr)
    }
    \mfrac{\displaystyle
      \prod_{\delta\in\corners{\dgr\sqcup S\,\remove\,1'2'}}
      \Bigl(\!1-\q^{\HookDistance(\delta,\,\Pos{S}{2'})}\!\Bigr)}{
      \displaystyle      
      \prod_{\substack{S_1\in\Orb_{2'}(\dgr\sqcup S\,\remove\,2')\\S_1\neq S}}
      \!\!\Bigl(\!1-\q^{\HookDistance(\Pos{S_1}{2'},\,\Pos{S}{2'})}\!\Bigr)},
  \end{gather*}
  and
  \begin{equation*}
    C_{S} =
    \begin{cases}\ytableausetup{smalltableaux}
      1 - \q^{-\HookDistance(\Pos{S}{2'},\,
        \Pos{S}{1'})},&S \ \ \text{contains}\ \ \Ytableaui{2'&1'} \ \
      \text{or}\ \ \Ytableaui{2'\\1'},
      \\
      1,&\text{otherwise}.
    \end{cases}\ytableausetup{nosmalltableaux}%
  \end{equation*}
\end{lemma}
We write $X_S$, $Y_S$, etc., instead of the more rigorous
$X_{\dgr\sqcup S}$, $Y_{\dgr\sqcup S}$, etc.

\subsubsection{Examples}For $\dgr\sqcup S$ in \eqref{example-S},
\begin{equation*}
  Y_S=-\mfrac{\q^2 (1- \q^{-1})}{(1- \q)^2}
  \mfrac{1- \q^{-2}}{(1- \q^{-3})
    (1- \q^3)}
  \mfrac{(1- \q^{-3}) (1- \q)}{(1- \q^{-4})
    (1- \q^2)}.
\end{equation*}
Another example is
\begin{multline*}
  \dgr\sqcup S
  =\Ytableaui{{}& {}& {}& {}\\ {}& {} & 2'\\
     {}\\ 1'}
   \\  \Longrightarrow\ 
   Y_S=
   %% \mfrac{-\q^{-3}}{(1- \q)^2}
   \mfrac{\q^{-4}}{(1- \q)(1- \q^{-1})}
   \mfrac{(1- \q^{-6})
   (1- \q^{-1})}{(1- \q^{-7})
   (1- \q^{-5})
   (1- \q^{-2})}
 \mfrac{(1- \q^{-2})
   (1- \q) 
   (1- \q^3)}{(1- \q^{-3}) (1- \q^2)}.
\end{multline*}

\subsubsection{The end of the proof of~\bref{Thm:action}}
With all the instances where $\Omega\acts(B,S,W)$ is nonzero
calculated in~\bref{lemma:1234-zero} and~\bref{lemma:Lambda},
we can finally calculate $(h_1 - g_1)\Omega$ on any vector of the
seminormal basis.  The diagonal parts of the action of $g_1$
in~\eqref{g-zeta} and of $h_1$ in
\begin{equation*}
  h_{1}\acts (B, S, W) =
  \mfrac{\alpha + \beta}{1 - \fffrac{\mu_{n-2}}{\mu_{n-1}}} (B, S, W)
  + \eta\bigl(\ffrac{\mu_{n-2}}{\mu_{n-1}}\bigr)(B, S, W_{(n-1,n)})
\end{equation*}
cancel when $g_1 - h_1$ is applied to each term in~\eqref{sum-z}; in
particular, all tableaux containing
\ytableausetup{smalltableaux}$\Ytableau{2'\\1'}$, or
$\Ytableau{2'&1'}$\ytableausetup{nosmalltableaux} cancel, and we are
left with\\[-1.3\baselineskip]
\begin{multline*}\ytableausetup{smalltableaux}
  (h_1 - g_1)\Omega\acts(B,S,W)
  =
    \sum_{\substack{S_1\in\Orbii(\dgr\sqcup S)\\
  S_1\not\ni\,\mbox{\tiny$\begin{ytableau}2'\\1'
    \end{ytableau}, \begin{ytableau}2'&1'
    \end{ytableau}$}}}
  \ffrac{1}{\kappa^2}\,z^{\phantom{y}}_{S,S_1}
  \biggl(
  \eta\bigl(\q^{\HookDistance(\Pos{S_1}{2'},\,\Pos{S_1}{1'})}\bigr)
  (B,S_1,W_1)_{(n,n-1)}
  \\[-12pt]
  -\zeta\bigl(\q^{\HookDistance(\Pos{S_1}{2'},\,\Pos{S_1}{1'})}\bigr)
    (B,S_1,W_1)_{(1',2')}\biggr).
  \end{multline*}
  Here, in view of the structure of the $1'2'$-orbit without the
  $\Ytableau{2'\\1'}$ and $\Ytableau{2'&1'}$ configurations, the terms
  can be collected pairwise, and hence
\begin{gather*}
  (h_1 - g_1)\Omega\acts(B,S,W) =
  \sum_{\substack{S_1\in\Orbii(\dgr\sqcup S)\\
  S_1\not\ni\,\mbox{\tiny$\begin{ytableau}2'\\1'
    \end{ytableau}, \begin{ytableau}2'&1'
    \end{ytableau}$}}}
  \ffrac{1}{\kappa^2}\,z^{\phantom{y}}_{S,S_1} y^{\phantom{y}}_{S,S_1}
  (B,S_1,W_1)_{(1',2')}.
\end{gather*}\ytableausetup{nosmalltableaux}
The factors that ``compare'' the $z_{S,S_1}$ coefficients in front of
similar terms are given by
\begin{multline*}
  y^{\phantom{y}}_{S,T}=
  -\zeta\bigl(\q^{\HookDistance(\Pos{T}{2'},\,\Pos{T}{1'})}\bigl)
  \\
  {}+ 
  \mfrac{\theFF\bigl(\q^{\HookDistance(\Pos{S}{2'},\, \Pos{T}{2'})}\bigl)}{
    \theFF\bigl(\q^{\HookDistance(\Pos{S}{2'},\, \Pos{T}{1'})}\bigl)}
  \mfrac{\widetilde{\eta}\bigl(\q^{\HookDistance(\Pos{S}{2'},\,
      \Pos{T}{1'})}\bigr)
    \widetilde{\zeta}\bigl(\q^{\HookDistance(\Pos{S}{2'},\, \Pos{T}{1'})}\bigr)}
  {\widetilde{\eta}\bigl(\q^{\HookDistance(\Pos{S}{2'},\, \Pos{T}{2'})}\bigr)
    \widetilde{\zeta}\bigl(\q^{\HookDistance(\Pos{S}{2'},\,
      \Pos{T}{2'})}\bigr)}\,
  \eta\bigl(\q^{-\HookDistance(\Pos{T}{2'},\,\Pos{T}{1'})}\bigr),
\end{multline*}
where, to avoid proliferating cases in the formula for
$y^{\phantom{y}}_{S,T}$, we define
\begin{gather*}
  \widetilde{\eta}(x) =
  \begin{cases}
    1,&x=1,\\
    \eta(x)&\text{otherwise},
  \end{cases}
  \ \
  \widetilde{\zeta}(x) =
  \begin{cases}
    1,&x=1,\\
    \zeta(x)&\text{otherwise},
  \end{cases}
  \ \
  \theFF(x) =
  \begin{cases}
    -\fffrac{1}{\alpha\beta},&x=1,\\
    \theF(x)&\text{otherwise}.
  \end{cases}
\end{gather*}

Using~\eqref{eta-inv} and~\eqref{zeta-inv}, we now conclude that $(g_1
- h_1) \Omega\acts(B,S,W) = 0$ for all standard triples if and only
if~\eqref{etazeta-consistency} holds.  This completes the proof
of~\bref{Thm:action}.

\subsubsection{Remarks}
\begin{enumerate}
\item Defined for generic parameter values (at which the algebra is
  semisimple), the seminormal representation $\SN{m,n,\dgri,\dgrii}$
  is irreducible, and is then readily identified with the irreducible
  $\qwB_{m,n}$ representation defined by two Young diagrams
  $(\dgri,\dgrii)$, of dimension $\|\dgri\|\|\dgrii\|\mbox{%
    \scriptsize$\displaystyle\binom{m}{f}\binom{n}{f} $}f!$, \
  $f=m-|\dgri|=n-|\dgrii|$ (where $\|\lambda\|$ is the number of
  standard tableaux built on a Young diagram $\lambda$).  In
  particular, this dimension is the number of standard triples in
  $\HyperT_{m,n}(\dgri, \dgrii)$.

\item The matrices for all $\qwB$ generators acquire denominators in
  the seminormal basis, which is the technical reason why seminormal
  representations cease to exist at certain parameter values.
  Studying these denominators (which are very explicit in our
  formulation) can be rather informative (cf., e.g.,~\cite{[Steen]} in
  a different, but not unrelated~context).
\end{enumerate}

\subsection{Scalar product}\label{sec:scalar-p}
We next construct an invariant scalar product on each seminormal
representation $\SN{m,n,\dgri,\dgrii}$.

\subsubsection{Entropy} We define the entropy of a pair $(\tbl,w)$,
where $\tbl$ is a Young tableau filled with $1$,\dots,$n$ and $w$ is
an $n$-dimensional vector with nonzero components.  We write
$\tbl=\sigma\tbl^{\text{s}}$, where $\tbl^{\text{s}}$ is the
(row-reading) superstandard tableau and $\sigma\in\mathbb{S}_n$
(assumed to be in reduced form) acts by permuting entries.  We then
define an $\mathbb{S}_n$ action on $w$ by specifying how the
elementary transpositions~act:
\begin{equation*}
  \sigma_{i,i+1}(w_1,\dots,w_n)
  = \ffrac{w_i}{w_{i+1}}
  (w_1,\dots,w_{i+1},w_i,\dots,w_n)
\end{equation*}
(which \textit{is} a representation of the symmetric group).  Then the
entropy $\Ent(\tbl,w)$ is the collection (unordered list) of all
factors thus obtained in acting with $\sigma$ on~$w$ (\textit{not} the
product of factors).  The inherent nonuniqueness does not affect the
final result in view of the nature of the functions applied to the
entropy in what follows.

If $w$ is an $n$-dimensional vector, but a standard tableau $\tbl$ is
filled not with $1,\dots,n$ but with some numbers $a_{i_1}< \dots <
a_{i_n}$, then there is a unique monotonic map
$\phi:(a_{i_1},\dots,a_{i_n})\to(1,\dots,n)$, and we set
$\Ent(\tbl,w)=\Ent(\phi\tbl,w)$.  For example,
\ytableausetup{smalltableaux}%
$\phi:\Ytableau{2&6&7\\4&9}\mapsto\Ytableau{1&3&4\\2&5}$.%
\ytableausetup{nosmalltableaux}

If, further, $\tbl$ is an \textit{anti}standard tableau (and, as in
our case, is built on a skew shape), then we define $\sigma$ by
comparing with the super-antistandard row-reading tableau based on the
same skew shape, and then apply $\sigma$ to $w$ by the same rule as
above.

We apply this construction to standard triples as follows. Given a
standard triple $(B,S,W)\in\HyperT_{m,n}(\dgri, \dgrii)$ and the
corresponding weight
$(\mu_1,\dots,\linebreak[0]\mu_{m+n-1})=\Wt(B,S,W)$, we set
\begin{equation*}
  \Ent(W)=\Ent(W,(1,\mu_1,\dots,\mu_{n-1}))
\end{equation*}
(with $1$ prepended to the first $n-1$ components of the weight so as
to make an $n$-compo\-nent vector).

Quite similarly,
\begin{equation*}
  \Ent(B)=\Ent(B,[\mu_{n},\dots,\mu_{n+m-1}]_{B}),
\end{equation*}
where the vector comprises those components among $\mu_{n}$, \dots,
$\mu_{n+m-1}$ that correspond to the entries of $B$ (not proportional
to~$\theta$; see~\bref{sec:weights}); the corresponding $\phi$ map is
understood here wherever needed.

Finally, to define the entropy of $S$, we extract $\sigma$ from
$S=\sigma S^{\text{as}}$, where $S^{\text{as}}$ is a
super-antistandard tableau as noted above.  With this ``antistandard
option'' indicated by a prime, we set
\begin{equation*}
  \Ent(S)=\Ent'(S,[\mu_{n+m-1},\dots,\mu_{n}]_{S}),
\end{equation*}
where this time only the components corresponding to the entries of
$S$ (those proportional to~$\theta$) are selected and placed in
reverse order, in accordance with their association with the entries
of~$S$.

\subsubsection{Examples} For the second standard triple
in~\bref{hyper-example}, we have
\begin{equation*}
  W=\Ytableaui{
    1&3&4\\
    2&5\\
    6
  }\,,\qquad
  w=(1,\mu_1,\dots,\mu_5)=
  (1, \q, \q^{-1}, \q^{-2}, 1, \q^2).
\end{equation*}
Then $\sigma=\sigma_{3,4}\sigma_{2,3}$, and we have
\begin{align*}
  \sigma_{2,3}w&=\q^2(1, \q^{-1}, \q, \q^{-2}, 1, \q^2),
  \\
  \sigma_{3,4}(1, \q^{-1}, \q, \q^{-2}, 1,\q^2) &=
  \q^3(1, \q^{-1}, \q^{-2}, \q, 1, \q^2),
\end{align*}
whence the entropy $\Ent(W)=\{q^3, q^2\}$.

For the antistandard tableau entering the same standard triple
in~\bref{hyper-example}, we have
\begin{equation*}
  S=\Ytableaui{
    \none[\cdot]&\none[\cdot]&4'\\
    6'&3'\\
    1'
  }\,,\qquad
  w= (\mu_6,\mu_4,\mu_3,\mu_1)
  = (-\theta \q^{-1},-\q^2 \theta,-\theta, -\theta \q^{-2}).
\end{equation*}
The ordering in $4',6',3',1'$ differs from a super-antistandard by a
single transposition, and hence the entropy is
$\Ent(S)=\{\mu_6/\mu_4\}=\{\q^{-3}\}$.

As another example, we consider the standard triple
\begin{equation*}
  (B,S,W)=
  \Threetableaux{\Ytableaui{2'& 6'\\ 5'}}{\Ytableaui{\none[\cdot]& \none[\cdot]& \none[\cdot]& 1'\\ 7'& 4'\\ 3'}}{\Ytableaui{1& 3& 4& 7\\ 2& 5\\ 6}}
  \in
  \HyperT_{7,7}\Bigl(\Ytableaui{{}& {}\\ {}}, \Ytableaui{{}&{}&{}}\Bigr)
\end{equation*}
and the associated weight
\begin{equation*}
  \mu
  =\bigl(\q, \q^{-1}, \q^{-2}, 1, \q^2, \q^{-3}; -\theta \q^3 , 1, -\theta \q^{-2}, -\theta, \q, \q^{-1}, -\theta \q^{-1}\bigr).
\end{equation*}
The corresponding entropies are \ytableausetup{smalltableaux}%
\begin{align*}
  \Ent(W)&=\Ent(W,(1,\q, \q^{-1}, \q^{-2}, 1, \q^2, \q^{-3}))
           =\{\q^4, \q^3, \q^2, \q^3, \q^5\},
  \\
  \Ent'(S)&=\Ent\bigl(\Ytableaui{\none[\cdot]& \none[\cdot]& \none[\cdot]&1\\
  4& 3\\ 2}\,,\,(-\theta \q^{-1}, -\theta, -\theta \q^{-2}, -\theta \q^3)\bigr)
  =\{\q^{-4},\q^{-3},\q^{-5}\},
  \\
  \Ent(B)&=\Ent\bigl(\Ytableaui{1& 3\\ 2}\,,\,(1, \q, \q^{-1})\bigr)
  = \{\q^2\}.
\end{align*}
\ytableausetup{nosmalltableaux}%

\subsubsection{Mutual entropy} Given a standard triple
$(B,S,W)\in\HyperT_{m,n}(\dgri, \dgrii)$, we define the mutual entropy
of $B$ and $S$ as the following (unordered) set of $\q^k$:
\begin{equation*}
  \Ent(B,S)=\bigsqcup_{i'\in S}\bigsqcup_{\substack{j'\in B\\ j'<i'}}
  \q^{\HookDistance(S[i'])-\HookDistance(B[j'])}.
\end{equation*}

This can be interpreted as follows. From the weight
$(\mu_1,\dots,\mu_{m+n-1})$, we select the last $m$ components
$\mu_{n},\dots,\mu_{m+n-1}$; among these, there are $f$ components
$\mu_{i_1},\dots,\mu_{i_f}$ corresponding to the entries of $S$.  For
each such $\mu_{i_r}$, we take all ratios $\mu_{i_r}/\mu'$, where
$\mu'$ ranges over components \textit{to the left of $\mu_{i_r}$} that
correspond to the entries of $B$.  The right-hand side of the last
formula is the collection of all such ratios.

\subsubsection{Shape factor}\label{sec:shape}
The scalar product that we define on standard triples in what follows
is made of the entropies, which depend on the \textit{tableaux}, and
of a factor that depends only on the shape of the diagrams.  We now
define the shape factors for each $\SN{m,n,\dgri,\dgrii}$.

For a standard triple $(B,S,W)\in\HyperT_{m,n}(\dgri, \dgrii)$, the
shapes of $S$ and $W$ are related as
$s\equiv\shape{S}=\shape{W}/\dgr$, and, as previously, we prefer
dropping the shape $\shape{W}$ from the notation, replacing it with
$\dgr\sqcup s$.  The shapes in $\HyperT_{m,n}(\dgri, \dgrii)$ are
\begin{equation}\label{shapes}
  \dgr\sqcup s
  =(\dgr_1+f_1,\dots,\dgr_k+f_k; f_{k+1},\dots,f_{\ell}),
\end{equation}
where $f_1+\dots + f_{\ell}=f\equiv n-|\dgrii|$, $\ell\geq k$, with
$f_1\geq 0$, \dots, $f_k\geq 0$ and $f_{k+1}\geq \dots\geq
f_{\ell}\geq 1$ (we use a semicolon to separate rows $i$ with
$\lambda_i>0$).

We select a reference shape where $s$ is a single row (of $f$ boxes)
that extends the top row of $\dgr$:
\begin{equation*}
  \dgr\sqcup s^{\text{ref}}
  =(\dgr_1+f,\dgr_2,\dots,\dgr_k).
\end{equation*}
%% (we recall that $f=m-|\dgri|=n-|\dgr|$).
We set $\theD(\dgr\sqcup s^{\text{ref}})=1$, and define the shape
factor $\theD(\cdot)$ for all other shapes recursively.

For $\dgr\sqcup s$ in~\eqref{shapes}, choose a corner $\Circ_J$ of
$\dgr\sqcup s$ in any of the rows $2\leq J\leq\ell$ such that
$f_{J}>0$, i.e., a corner not belonging to $\dgr$. (If there is no
such corner, then we already have $s=s^{\text{ref}}$.)  \ We let
$\Phi_J(\dgr\sqcup s)$ denote the diagram where the chosen corner is
moved to the end of the top row; depending on where the corner was
chosen, $\Phi_J(\dgr\sqcup s)$ is one of the diagrams
\begin{align*}
  &(\dgr_1+f_1+1,\dgr_2+f_2,\dots,\dgr_J + f_J -1,\dots,\dgr_k+f_k; f_{k+1},\dots,f_{\ell})
  \\
  \intertext{or}
  &(\dgr_1+f_1+1,\dgr_2+f_2,\dots,\dgr_k+f_k; f_{k+1},\dots,f_J-1,\dots,f_{\ell}).
\end{align*}
We then set
\begin{equation}\label{shape-factor}
  \theD(\dgr\sqcup s)
  =\mfrac{\Mass_1\bigl(\Phi_J(\dgr\sqcup s)\bigr)}{\Mass_J(\dgr\sqcup s)}
  \mfrac{\Massi_1\bigl(\Phi_J(\dgr\sqcup s),\dgr\bigr)}{\Massi_J(\dgr\sqcup s,\dgr)}\,
  \theD\bigl(\Phi_J(\dgr\sqcup s)\bigr),
\end{equation}
where, for a Young diagram of $K$ rows,
\begin{alignat*}{2}
  \Mass_J(\Dgr)
  &=\prod_{\Box\,\in\Dgr(J+1,K)}
  \theH(\q^{\HookDistance(\Box,\,\Circ_J)}),&\qquad
  \theH(x)&=-\mfrac{\eta(x)^2}{\alpha \beta \theF(x)},
  \\
  \Massi_J(\Dgr,\dgr)
  &=\prod_{\substack{\Box\,\in\Dgr(J+1,K)\\
      \Box\,\notin\dgr%%%_{J+1}\cup\dots\cup\dgr_k
    }}\theZ(\q^{\HookDistance(\Box,\,\Circ_J)}),
  &\qquad \theZ(x)&=-\mfrac{\zeta(x)^2}{\alpha \beta \theF(x)}.
\end{alignat*}
Here, for a Young diagram $\Dgr$ and $i\leq j$, we let
$\Dgr(i,j)=\Dgr_i\sqcup\dots\sqcup\Dgr_j$ denote the part of the
diagram made of the rows between (and including) the $i$th and $j$th
rows.  The products in the last two formulas are over boxes
\textit{strictly below the chosen corner} in the $J$th row: the boxes
of $\Dgr$ in the first case and the boxes of $\Dgr/\dgr$ in the
second.

We note that from~\eqref{etazeta-consistency}, we have
\begin{equation}\label{HZ-consistency}
  \theH(x) \theZ(x) \theH(y) \theZ(y)=
  \theH(x y) \theZ(x y),
\end{equation}
which is a consistency condition for the construction of
$\theD(\dgr\sqcup s)$.

Formula~\eqref{shape-factor} defines $\theD(\cdot)$ for all shapes
encountered in $\HyperT_{m,n}(\dgri, \dgrii)$, by gradually moving all
boxes not belonging to $\dgr$ to the top row.

\subsubsection{Example}As an example of the calculation of the $\theD$
factor, we take $\lambda=(2,2,1)$ and choose a skew shape $s$ such
that $\dgr\sqcup s = (4,2,2,1)$. \
Then\ytableausetup{smalltableaux}%
\begin{align*}
  \theD\Bigl(
  \Ytableaui{\none[\cdot]&\none[\cdot]&{}&{}\\
  \none[\cdot]&\none[\cdot]\\
  \none[\cdot]&{}\\
  {}
  }\Bigr)
  &=\theH(\q^4)\theH(\q^5)^2\theH(\q^6)\theZ(\q^5)\,\theD\Bigl(
  \Ytableaui{\none[\cdot]&\none[\cdot]&{}&{}&{}\\
    \none[\cdot]&\none[\cdot]\\
    \none[\cdot]&{}
    }\Bigr)
  \\
  &=\theH(\q^4)\theH(\q^5)^2\theH(\q^6)\theZ(\q^5)\cdot
    \theH(\q^5)\theH(\q^6)\theH(\q^7)\,
    \theD\Bigl(
  \Ytableaui{\none[\cdot]&\none[\cdot]&{}&{}&{}&{}\\
    \none[\cdot]&\none[\cdot]\\
    \none[\cdot]
  }\Bigr)
  \\
  &=\theH(\q^4)\theH(\q^5)^3\theH(\q^6)^2\theH(\q^7) \theZ(\q^5).
\end{align*}
\ytableausetup{nosmalltableaux}%
Moving the boxes around starting from another corner gives
$\fffrac{\theZ(\q^7) \theH(\q^4) \theH(\q^5)^2 \theH(\q^6)^2
  \theH(\q^7)^2}{\theZ(\q^2) \theH(\q^2)}$,
which is the same in view of~\eqref{HZ-consistency}.

In the next theorem, we speak of a \textit{diagonal} scalar product of
standard triples, i.e., such that $((B, S,W), (B', S',W')) = 0$ unless
$B = B'$, $S = S'$, and $W = W'$.

\begin{thm}\label{thm:scalar}
  The diagonal scalar product $(\cdot,\cdot)$ on standard triples in
  $\HyperT_{m,n}(\dgri, \dgrii)$ \textup{(}the seminormal basis of
  $\SN{m,n,\dgri,\dgrii}$\textup{)} defined by
  \begin{equation*}
    \bigl((B,S,W),(B,S,W)\bigr)=
    \mfrac{\theD(\dgr\sqcup\shape{S})}{\qDim(\dgr\sqcup\shape{S},\q)}
    \prod_{\Box\,\in S}\Bigl(\!1 + \ffrac{\theta}{\q^{\HookDistance(\Box)}}\Bigr)
    \mfrac{\displaystyle
      \prod_{t\in\Ent(W)}\theH(t)\;
      \prod_{v\in\Ent(B)}\theZ(v)
    }{\displaystyle
      \prod_{u\in\Ent'(S)}\theZ(u)
      \prod_{w\in\Ent(B,S)}\theZ(w)
    }
  \end{equation*}
  is invariant under the action of $\qwB_{m,n}$ generators:
%%   defined in \bref{h-action}, \bref{g-action}, and~\bref{E-action}:
  $(A\acts\hyp_1,\hyp_2)=(\hyp_1,A\acts\hyp_2)$ for any standard
  triples $\hyp_1$ and $\hyp_2$, and $A$ any of the $\qwB_{m,n}$
  generators $g_j$, $\EE$, or~$h_i$.
\end{thm}
The quantum dimension $\qDim$ of a Young diagram is defined
in~\bref{sec:qdim}.

The proof is by direct verification.  Showing the invariance under
$g_j$ and $h_i$ amounts to a standard analysis of cases (which are
somewhat more numerous for the $g_j$).  As regards the action of
$\EE$, we consider two skew tableaux $S_1$ and $S_2$ belonging to the
orbit of the mobile element; they differ by the position of a single
box.  For the coefficients $c^{(1)}_{S_i}$ in~\eqref{cSS1}, their
ratio $\nofrac{c^{(1)}_{S_1}}{c^{(1)}_{S_2}}$ is evidently reproduced
from the ratio of the products $\prod_{\Box\,\in S}$ in the formula
for the scalar product.  Moreover, for the coefficients
$c^{(2)}_{S_i}$, we have
\begin{equation*}
  \mfrac{c^{(2)}_{S_1}}{c^{(2)}_{S_2}}=
  \mfrac{\qDim(\dgr\sqcup\shape{S_1},\q)}{
    \qDim(\dgr\sqcup\shape{S_2},\q)},
  \qquad
  \q=-\fffrac{\beta}{\alpha},
\end{equation*}
leading to the desired result.

\begin{rem}
  The formula for the scalar product considerably simplifies for the
  totally symmetric choice in~\eqref{fully-symmetric}: then $\theH$,
  $\theZ$, and $\theD$ are identically equal to~$1$, and
  \begin{equation*}
    \bigl((B,S,W),(B,S,W)\bigr)=
    \mfrac{1}{\qDim(\shape{W},\q)}
    \prod_{\Box\,\in S}
    \Bigl(\!1 + \ffrac{\theta}{\q^{\HookDistance(\Box)}}\Bigr)
  \end{equation*}
  (where, of course, the shape of the ``white'' tableau is
  $\shape{W}=\dgr\sqcup\shape{S}$).
\end{rem}

\subsubsection{Example}
For the $16$-dimensional seminormal representation with $m = 2$,
$n = 4$, and \ytableausetup{smalltableaux}%
$ (\dgr',\dgr)= \bigl(\Ytableaui{{} }, \Ytableaui{
  {}&\\
  {} }\bigr) $,\ytableausetup{nosmalltableaux}
we list some (a half) of its basis vectors---standard
triples~$\hyp$---and their scalar squares in the format
$\hyp\to(\hyp,\hyp)\to(\hyp,\hyp)\bigr|_{\text{sym}}$, where the last
term is the form taken by the scalar product for the totally symmetric
choice~\eqref{fully-symmetric}.  For $\hyp=(B,S,W)$ such that
$1'\in B$, we have \ytableausetup{smalltableaux}
\begin{alignat*}{2}
%%   \Bigl(\Ytableaui{
%%     1' \\
%%   }
%%   \,,\; 
%%   \Ytableaui{
%%     \none[\cdot] & \none[\cdot] \\
%%     \none[\cdot] & 2' \\
%%   }
%%   \,,\; 
%%   \Ytableaui{
%%     1 & 2 \\
%%     3 & 4 \\
%%   }\Bigr)&\to\ffrac{(\frac{\theta}{\q}+1)(1-\q^3)^2 (\q \theta +1)
%%     \eta(\q^3)^2}{\q^2 (1-\q^4)^2 (\theta +1) \zeta(-\theta)^2}
%% &&\to\ffrac{(\theta +1) (1-\q^2)}{\q^2
%%    (1-\q^4)},
%% \\
  \Bigl(\Ytableaui{
    1' \\
  }
  \,,\; 
  \Ytableaui{
    \none[\cdot] & \none[\cdot] \\
    \none[\cdot] & 2' \\
  }
  \,,\; 
  \Ytableaui{
    1 & 3 \\
    2 & 4 \\
  }\Bigr)&\to
  \ffrac{(\frac{\theta }{\q}+1) (\q \theta +1)(1-\q^2)^2 (1-\q^3)   
   \eta(\q^2)^2 \eta(\q^3)^2}{(1-\q) \q^3
   (1-\q^4)^2 \alpha^2 (\theta +1)
   \zeta(-\theta)^2}
&&\to\ffrac{(\theta +1)
   (1-\q^2)}{\q^2
   (1-\q^4)},
\\
%%   \Bigl(\Ytableaui{
%%     1' \\
%%   }
%%   \,,\; 
%%   \Ytableaui{
%%     \none[\cdot] & \none[\cdot] & 2' \\
%%     \none[\cdot]   \\
%%   }
%%   \,,\; 
%%   \Ytableaui{
%%     1 & 2 & 3 \\
%%     4   \\
%%   }\Bigr)&\to\ffrac{(\q \theta +1) (\theta \q^3+1)(1-\q) \alpha^2}{(1-\q^3)
%%    (\theta \q^2+1) \zeta(-\q^2 \theta)^2}
%% &&\to\ffrac{(\theta \q^2+1)(1-\q)}{\q (1-\q^3)},
%% \\
%%   \Bigl(\Ytableaui{
%%     1' \\
%%   }
%%   \,,\; 
%%   \Ytableaui{
%%     \none[\cdot] & \none[\cdot] & 2' \\
%%     \none[\cdot]   \\
%%   }
%%   \,,\; 
%%   \Ytableaui{
%%     1 & 2 & 4 \\
%%     3   \\
%%   }\Bigr)&\to\ffrac{(\q \theta +1) (\theta  \q^3+1)(1-\q) (1-\q^3) \eta(\q^3)^2}{\q (1-\q^2)
%%    (1-\q^4) (\theta  \q^2+1)
%%    \zeta(-\q^2 \theta)^2}
%% &&\to\ffrac{(\theta  \q^2+1)(1-\q)}{\q (1-\q^3)},
%% \\
  \Bigl(\Ytableaui{
    1' \\
  }
  \,,\; 
  \Ytableaui{
    \none[\cdot] & \none[\cdot] & 2' \\
    \none[\cdot]   \\
  }
  \,,\; 
  \Ytableaui{
    1 & 3 & 4 \\
    2   \\
  }\Bigr)&\to\ffrac{(\q \theta +1) (\theta  \q^3+1)
    (1-\q^2) \eta(\q^2)^2 \eta(\q^3)^2}{\q^2 (1-\q^4)
   \alpha^2 (\theta  \q^2+1) \zeta(-\q^2 \theta)^2}
&&\to\ffrac{(\theta \q^2+1)(1-\q)}{\q (1-\q^3)},
\\
%%   \Bigl(\Ytableaui{
%%     1' \\
%%   }
%%   \,,\; 
%%   \Ytableaui{
%%     \none[\cdot] & \none[\cdot] \\
%%     \none[\cdot]  \\
%%     2'  \\
%%   }
%%   \,,\; 
%%   \Ytableaui{
%%     1 & 2 \\
%%     3  \\
%%     4  \\
%%   }\Bigr)&\to\ffrac{(\frac{\theta}{\q^3}+1) (\frac{\theta }{\q}+1)
%%     (1-\q) (1-\q^3) 
%%    \eta(\q^3)^2}{\q^3
%%    (1-\q^2) (1-\q^4)
%%    (\frac{\theta }{\q^2}+1) \zeta(-\frac{\theta}{\q^2})^2}
%% &&\to\ffrac{(\frac{\theta }{\q^2}+1)(1-\q)}{\q^3 (1-\q^3)},
%% \\
    \Bigl(\Ytableaui{
    1' \\
  }
  \,,\; 
  \Ytableaui{
    \none[\cdot] & \none[\cdot] \\
    \none[\cdot]  \\
    2'  \\
  }
  \,,\; 
  \Ytableaui{
    1 & 3 \\
    2  \\
    4  \\
  }\Bigr)&\to\ffrac{(\frac{\theta }{\q^3}+1)
    (\frac{\theta }{\q}+1)
    (1-\q^2) \eta(\q^2)^2 \eta(\q^3)^2}{\q^4 (1-\q^4)
   \alpha^2 (\frac{\theta }{\q^2}+1)
   \zeta(-\frac{\theta}{\q^2})^2}
&&\to\ffrac{(\frac{\theta}{\q^2}+1)(1-\q)}{\q^3
   (1-\q^3)},
\\
  \Bigl(\Ytableaui{
    1' \\
  }
  \,,\; 
  \Ytableaui{
    \none[\cdot] & \none[\cdot] \\
    \none[\cdot]  \\
    2'  \\
  }
  \,,\; 
  \Ytableaui{
    1 & 4 \\
    2  \\
    3  \\
  }\Bigr)&\to\ffrac{(\frac{\theta }{\q^3}+1)
    (\frac{\theta }{\q}+1)
    (1-\q^3)^2
    \eta(\q^2)^2 \eta(\q^3)^4}{\q^5 (1-\q^4)^2
   \alpha^4 (\frac{\theta }{\q^2}+1)
   \zeta(-\frac{\theta}{\q^2})^2}
&&\to\ffrac{(\frac{\theta}{\q^2}+1)(1-\q) }{\q^3 (1-\q^3)},
\end{alignat*}
and for those with $1'\in S$,
\begin{alignat*}{2}
%% \Bigl(\Ytableaui{
%%     2' \\
%%   }
%%   \,,\; 
%%   \Ytableaui{
%%     \none[\cdot] & \none[\cdot] \\
%%     \none[\cdot] & 1' \\
%%   }
%%   \,,\; 
%%   \Ytableaui{
%%     1 & 2 \\
%%     3 & 4 \\
%%   }\Bigr)&\to\ffrac{(\theta +1)(1-\q^3)^2 \eta(\q^3)^2}{
%%     \q^3 (1-\q^4)^2 \alpha^2}
%% &&\to\ffrac{(\theta +1) (1-\q^2)}{\q^2
%%    (1-\q^4)},
%% \\
  \Bigl(\Ytableaui{
    2' \\
  }
  \,,\; 
  \Ytableaui{
    \none[\cdot] & \none[\cdot] \\
    \none[\cdot] & 1' \\
  }
  \,,\; 
  \Ytableaui{
    1 & 3 \\
    2 & 4 \\
  }\Bigr)&\to\ffrac{(\theta +1)(1-\q^2)^2 (1-\q^3)
    \eta(\q^2)^2
   \eta(\q^3)^2}{(1-\q) \q^4
   (1-\q^4)^2 \alpha^4}
&&\to\ffrac{(\theta +1) (1-\q^2)}{\q^2 (1-\q^4)},
\\
  \Bigl(\Ytableaui{
    2' \\
  }
  \,,\; 
  \Ytableaui{
    \none[\cdot] & \none[\cdot] & 1' \\
    \none[\cdot]   \\
  }
  \,,\; 
  \Ytableaui{
    1 & 2 & 3 \\
    4   \\
  }\Bigr)&\to\ffrac{(\theta  \q^2+1)(1-\q)}{\q
   (1-\q^3)}
&&\to\ffrac{(\theta \q^2+1)(1-\q)}{\q (1-\q^3)},
\\
%%   \Bigl(\Ytableaui{
%%     2' \\
%%   }
%%   \,,\; 
%%   \Ytableaui{
%%     \none[\cdot] & \none[\cdot] & 1' \\
%%     \none[\cdot]   \\
%%   }
%%   \,,\; 
%%   \Ytableaui{
%%     1 & 2 & 4 \\
%%     3   \\
%%   }\Bigr)&\to \ffrac{(\theta  \q^2+1)(1-\q) (1-\q^3)
%%     \eta(\q^3)^2}{\q^2 (1-\q^2)
%%    (1-\q^4) \alpha^2}
%% &&\to\ffrac{(\theta\q^2+1)(1-\q)}{\q
%%    (1-\q^3)},
%% \\
  \Bigl(\Ytableaui{
    2' \\
  }
  \,,\; 
  \Ytableaui{
    \none[\cdot] & \none[\cdot] & 1' \\
    \none[\cdot]   \\
  }
  \,,\; 
  \Ytableaui{
    1 & 3 & 4 \\
    2   \\
  }\Bigr)&\to \ffrac{(\theta\q^2+1)(1-\q^2) \eta(\q^2)^2
   \eta(\q^3)^2}{\q^3
   (1-\q^4) \alpha^4}
&&\to\ffrac{(\theta \q^2+1)(1-\q)}{\q (1-\q^3)},
\\
%%   \Bigl(\Ytableaui{
%%     2' \\
%%   }
%%   \,,\; 
%%   \Ytableaui{
%%     \none[\cdot] & \none[\cdot] \\
%%     \none[\cdot]  \\
%%     1'  \\
%%   }
%%   \,,\; 
%%   \Ytableaui{
%%     1 & 2 \\
%%     3  \\
%%     4  \\
%%   }\Bigr)&\to \ffrac{(\frac{\theta}{\q^2}+1)(1-\q)
%%    (1-\q^3)  \eta(\q^3)^2}{\q^4 (1-\q^2)
%%    (1-\q^4) \alpha^2}
%% &&\to\ffrac{(\frac{\theta }{\q^2}+1)(1-\q)}{\q^3
%%    (1-\q^3)},
%% \\
%%   \Bigl(\Ytableaui{
%%     2' \\
%%   }
%%   \,,\; 
%%   \Ytableaui{
%%     \none[\cdot] & \none[\cdot] \\
%%     \none[\cdot]  \\
%%     1'  \\
%%   }
%%   \,,\; 
%%   \Ytableaui{
%%     1 & 3 \\
%%     2  \\
%%     4  \\
%%   }\Bigr)&\to \ffrac{(\frac{\theta}{\q^2}+1)(1-\q^2) \eta(\q^2)^2
%%    \eta(\q^3)^2}{\q^5
%%    (1-\q^4) \alpha^4}
%% &&\to\ffrac{(\frac{\theta}{\q^2}+1)(1-\q)}{\q^3
%%    (1-\q^3)},
%% \\
  \Bigl(\Ytableaui{
    2' \\
  }
  \,,\; 
  \Ytableaui{
    \none[\cdot] & \none[\cdot] \\
    \none[\cdot]  \\
    1'  \\
  }
  \,,\; 
  \Ytableaui{
    1 & 4 \\
    2  \\
    3  \\
  }\Bigr)&\to \ffrac{(\frac{\theta}{\q^2}+1)(1-\q^3)^2 \eta(\q^2)^2
   \eta(\q^3)^4}{\q^6
   (1-\q^4)^2 \alpha^6}
 &&\to\ffrac{(\frac{\theta}{\q^2}+1)(1-\q)}{\q^3 (1-\q^3)}.
\end{alignat*}\ytableausetup{nosmalltableaux}%

\begin{Thm}\label{Thm:diagonalize}
  In a seminormal representation $\SN{m,n,\dgri,\dgrii}$, the
  Jucys--Murphy elements $\JMBare(n)_2$, \dots, $\JMBare(n)_{m+n}$
  \textup{(}see~\bref{sec:JM}\textup{)} act on the seminormal basis
  elements
  %% ---standard triples $\hyp\in\HyperT_{m,n}(\dgri, \dgrii)$---
  $\hyp=(B,S,W)$ as
  \begin{equation*}
    \JMBare(n)_j\acts\hyp = \Wt(\hyp)_{j-1}\;\hyp,\qquad j=2,\dots,m+n,
  \end{equation*}
  where the weight of a standard triple is defined
  in~\bref{sec:weights}.
\end{Thm}

\subsection{Proof}
\subsubsection{}The first $n-1$ Jucys--Murphy elements
$\JMBare(n)_{2}$, \dots, $\JMBare(n)_{n}$ are Jucys--Murphy elements
of the Hecke subalgebra $\Hecke_n\subset\qwB_{m,n}$, and the assertion
is well known~\cite{[Mu-81]}; we recall that it can be proved by
induction on $i$ in $\JMBare(n)_{i}$, based on
definition~\eqref{JM-white}, which can be equivalently rewritten as
\begin{equation*}
    \Bigl(-\ffrac{1}{\alpha\beta}\,h_{n+1-i} + \ffrac{\alpha+\beta}{\alpha\beta}\Bigr)
  \JMBare(n)_{i}
  = -\ffrac{1}{\alpha \beta}\JMBare(n)_{i-1}h_{n+1-i},\quad
    2\leq i\leq n,
\end{equation*}
and the fact that when $h_{n-i}$ acts nondiagonally, it gives rise to
a new weight that differs from the original weight $\mu$ by the
transposition of two neighboring components, $\mu_{i-1}$ and~$\mu_{i}$.

For the remaining Jucys--Murphy elements $\JMBare(n)_{n+1}$, \dots,
$\JMBare(n)_{n+m}$, the definition also implies the identities
\begin{equation*}
  \Bigl(-\ffrac{1}{\alpha\beta}g_{j-1} +
  \ffrac{\alpha+\beta}{\alpha\beta}\Bigr)\JMBare(n)_{n+j}
  =-\ffrac{1}{\alpha\beta}\,\JMBare(n)_{n+j-1} g_{j-1},\quad
  1\leq j\leq m,
\end{equation*}
and it is also the case that whenever the action of $g_{j'}$ produces
a new standard triple $((B, S)_{(j',j'+1)}, W)$, its weight differs
from the weight of $(B, S, W)$ by the transposition of $\mu_{n+j'-1}$
and $\mu_{n+j'}$.  Hence, by the same argument, the statement of the
theorem holds for $\JMBare(n)_{n+2}$, \dots, $\JMBare(n)_{n+m}$ as
soon as it holds for $\JMBare(n)_{n+1}$. \ It therefore remains to
establish the claim for~$\JMBare(n)_{n+1}$, i.e.,
\begin{equation}\label{remains}
  \JMBare(n)_{n+1}\acts\hyp = \Wt(\hyp)_{n}\;\hyp.
\end{equation}
The relevant component of the weight is determined by the position
of~$1'$ in the standard triple.

\subsubsection{}
The Jucys--Murphy element $\JMBare(n)_{n + 1}$ is not related to the
``lower'' ones by a simple formula, and we instead use its explicit
form found in~\bref{J-explicit}:
\begin{gather*}
  \JMBare(n)_{n + 1} =
  1 - (\alpha+\beta)\kappa\sum_{s=1}^{n}(-\alpha\beta)^{s-1}
  \hdowninv_{s-1,1}\,\EE\,\hupinv_{1,s-1}.
\end{gather*}
It readily follows that $\JMBare(n)_{n + 1}$ commutes with $h_1$,
\dots, $h_{n-1}$.  For $h_2$, \dots, $h_{n-1}$, which commute with
$\EE$, this is entirely a Hecke-algebra statement, and the
commutativity for $h_1$ is also immediate because, concentrating on
the generators that do not commute with $h_1$, we have
\begin{align*}
  h_1 \cdot h_2 h_1\,\EE\,h_1 h_2
  &= h_2 h_1 h_2\,\EE\,h_1 h_2\\
  &= h_2 h_1\,\EE\, h_2 h_1 h_2 \\
  &= h_2 h_1\,\EE\, h_1 h_2 \cdot h_1.
\end{align*}
Therefore, $\JMBare(n)_{n + 1}$ acts by an eigenvalue in each
irreducible representation of $\Hecke_{n}$.  Because
$\SN{m,n,\dgri,\dgrii}$ decomposes into a direct sum of irreducible
$\Hecke_{n}$ representations, it remains to find these eigenvalues.

\subsubsection{}All calculations for $\JMBare(n)_{n + 1}$ can be done
in $\qwB_{1,n}$, and, hence, in seminormal representations
$\SN{1,n,\emptyset,\dgrii}$ (with $|\dgr|=n-1$) and
$\SN{1,n,\Box,\dgrii}$ (with $|\dgr|=n$). \ The second case is
immediate, because $\EE$ then acts trivially, and therefore
$\JMBare(n)_{n + 1}$ acts as identity; but for all standard triples
$(\mbox{\tiny$\begin{ytableau}
    1'
  \end{ytableau}$},S,W)
\in\HyperT_{1,n}(\Box,\dgrii)$,
the weight component in~\eqref{remains} is indeed~$1$.

We are therefore left with the first case, i.e., finding the
eigenvalues of $\JMBare(n)_{n + 1}$ acting on standard triples
\begin{equation*}
  \hyp=(\emptyset,S,W),
\end{equation*}
where the skew shape $S$ is a single box, attached to $\lambda$ in one
of the possible positions; the sought eigenvalue depends on that
position (and the shape~$\dgrii$).

By the invariance property of the scalar product, the
sought eigenvalue is
\begin{equation}\label{eigen-work}
  \mfrac{(\hyp,\JMBare(n)_{n + 1}\hyp)}{(\hyp,\hyp)} =
  1 - (\alpha+\beta)\kappa\sum_{s=1}^{n}(-\alpha\beta)^{s-1}
  \mfrac{(\hupinv_{1,s-1}\hyp,\,\EE\,\hupinv_{1,s-1}\hyp)}{(\hyp,\hyp)}.
\end{equation}
For $\hyp$ of the above form, we choose $W\remove n$ to be the
row-reading superstandard tableau filled with $1$, \dots, $n-1$; its
shape, we recall, is $\lambda$.  Because the single box of $S$
carries~$1'$, this standard triple has a mobile element.  The
nondiagonal part of the action of $h_1^{-1}$ destroys the mobile
element, and hence (as many times in the foregoing) $h_1^{-1}$
effectively acts by an eigenvalue.  Next, each $h^{-1}_i$ in
$h^{-1}_2\dots h^{-1}_s$ acts by an eigenvalue whenever $n-i$ and
$n-i+1$ are in the same row.

The left-hand side of~\eqref{eigen-work} is therefore expressed as a
sum over the rows of $W$. \ More precisely, let $n$ stand in the $K$th
row of $W$, and let also $\Box_i$ be the last box in the $i$th row
(thus, the shape of $W$ is $\shape{W}=\dgr\sqcup\Box_K$).  The sum
over $i=1,\dots, K-1$ in the next formula is over the rows above the
$K$th one, to which the contribution of the rest of the diagram is
added.  In addition to the notation $\Dgr(i,j)$ introduced
in~\bref{sec:shape}, we let $\Dgr({\geq}j)$ and $\Dgr({>}j)$ denote
the parts a Young diagram made of rows nonstrictly and strictly below
a $j$th row.  Then
\begin{align*}
  \mfrac{(\hyp,\JMBare(n)_{n + 1}\hyp)}{(\hyp,\hyp)}
  \kern-30pt
  &\kern30pt=
  1 - \Bigl(\!1+\ffrac{\theta}{\q^{\HookDistance(\Box_K)}}\Bigr)
  \mfrac{\displaystyle
    \prod_{\delta\in\corners{\dgrii\,}}\bigl(1-\q^{\HookDistance(\delta,\,\Box_K)}
    \bigr)
  }{\displaystyle
    \prod_{\substack{\Star\in\cocorners{\dgrii\,}\\ \Star\neq\Box_K}}
    \bigl(1-\q^{\HookDistance(\Star,\,\Box_K)} \bigr)
  }
  \\
  &\times
  \Biggl(\sum_{i=1}^{K-1}(1-\q^{-1})
   \q^{\HookDistance(\Box_i,\,\Box_K)}\,
  \mfrac{\displaystyle
    \prod_{\Star\in\cocornersii{\dgrii(i+1,K-1)}}
    \bigl(1-\q^{\HookDistance(\Star,\,\Box_K)}\bigr)
  }{\displaystyle
    \prod_{\delta\in\corners{\dgrii(i,K-1)}}
    \bigl(1-\q^{\HookDistance(\delta,\,\Box_K)}\bigr)
  }
  +\mfrac{\displaystyle
    \prod_{\Star\in\cocornersi{\dgr({>K})}}
    \bigl(1-\q^{\HookDistance(\Star,\,\Box_K)}\bigr)
  }{\displaystyle
    \prod_{\delta\in\corners{\dgr({\geq}K)}}
    \bigl(1-\q^{\HookDistance(\delta,\,\Box_K)}\bigr)
  }
    \Biggr)
  \\
  &\kern30pt=1 - \Bigl(\!1+\ffrac{\theta}{\q^{\HookDistance(\Box_K)}}\Bigr)
    = -\ffrac{\theta}{\q^{\HookDistance(\Box_K)}}
    = -\ffrac{\theta}{\q^{\HookDistance(\Pos{S}{1'})}},
\end{align*}
which is the $n$th component of the weight defined
in~\bref{sec:weights}. \ This shows~\eqref{remains} and
hence~\bref{Thm:diagonalize}.

\section{Outlook}
We have discussed the quantum walled Brauer algebras $\qwB_{m,n}$
starting with the endomorphism algebras of mixed tensor products.
%% $X^*{}^{\otimes m}\otimes X^{\otimes n}$. \ 
We constructed the link-state basis in $\qwB$ Specht modules, a
Baxterization of the algebra (more precisely, of morphisms in an
``ambient'' category), and seminormal $\qwB$ representations for
generic parameters of the algebra, which allowed us to find the
spectrum of a family of Jucys--Murphy elements.

The results can be developed in various directions.  Among these, we
note finding a generalized seminormal basis for the special parameter
values $\theta=-\bigl(-\frac{\beta}{\alpha}\bigr)^r$, and
investigating lattice models$/$spin chains that can be constructed
from the monodromy matrix obtained by Baxterization.

At $\theta=-\bigl(-\frac{\beta}{\alpha}\bigr)^r$, $r\in\oZ$, suitable
quotients of the $\qwB$ algebra centralize the action of $q$-deformed
general linear Lie superalgebras on tensor products of their natural
representations; the $\qwB$ algebra becomes nonsemisimple, and the
Jucys--Murphy elements acquire root vectors in Specht modules.  A
generalization of the seminormal basis can then be defined as a basis
in which Jucys--Murphy elements take the standard Jordan form.  The
common Jordan structure of Jucys--Murphy elements then depends on $r$
and is a subject to be investigated.  A generalized seminormal basis
is important, in particular, in finding the bimodule structure of the
mixed tensor products of $U_qg\ell(M|N)$ representations.

The ``universal monodromy matrix'' resulting from the proposed
Baxterization relates to a ``universal spin chain,'' which yields
specific, true spin chains (corresponding to the $U_q g\ell(M|N)$
series) at special parameter values. It is of interest to develop an
appropriate version of the Bethe-ansatz approach and to trace how the
step-by-step degeneration descends from the universal model to a
specific spin chain$/$lattice model with the chosen $U_qg\ell(M|N)$
symmetry.
%% Among the particular models, the standard XXZ spin chain with
%% $U_qs\ell(2)$ symmetry and the $\qwB$ algebra degenerated into the
%% Temperley--Lieb algebra occur at~$r=2$.
Deeper insights are to be gained from the root-of-unity case
(cf.~\cite{[GN]}), where the centralizer of $U_qg\ell(M|N)$ is
expected to be a lattice $W$-algebra---a discretization
(cf.~\cite{[GST]}) of a $W$-algebra defined in two-dimensional
conformal field theory in terms of the intersection of kernels of the
screening operators corresponding to~$U_qg\ell(M|N)$.

\subsection*{Acknowledgments}We thank A.~Davydov, B.~Feigin,
M.~Finkelberg, A.~Gainutdinov, A.~Ki\-se\-lev, G.~Kufryk, S.~Lentner,
I.~Runkel, Y.~Saint-Aubin, and H.~Saleur for the useful discussions
and suggestions.  Special thanks, for the hospitality, go to
D.~Adamovic for the conference ``Representation Theory XIV'' and to
D.~Ridout and S.~Wood for ``The Mathematics of Conformal Field
Theory,'' where the above results were reported.  Advice from
Y.~Saint-Aubin is greatfully appreciated.  Very useful remarks,
suggestions, and corrections by the referee are gratefully
appreciated.  This paper was supported in part by the RFBR grant
13-01-00386.  The work of IYuT was supported in part by the ERC
Advanced Grant NuQFT.

\appendix

\section{Notation and conventions}

\subsection{}
We let $|\lambda|$ denote the number of boxes in a Young diagram or a
skew shape (or in fact in a tableau built on any of these).

\subsection{}\label{positions} By the \textit{position} of a box in a
Young diagram or a tableau or a skew shape, we mean the coordinates of
the box in a quadrant of $\oZ^2$, assigned in accordance with the
pattern
\begin{equation*}
  \begin{ytableau}
    \scriptscriptstyle\pos{1,1}& \scriptscriptstyle\pos{1,2} & \scriptscriptstyle\pos{1,3}  & \none[\dots] \\
    \scriptscriptstyle\pos{2,1}& \scriptscriptstyle\pos{2,2} \\
    \scriptscriptstyle\pos{3,1}&\none&\none[\ddots]\\
    \none[\vdots]
  \end{ytableau}
  \ .
\end{equation*}

For a tableau $\tbl$ containing a number $k$, we let
$\Pos{\tbl}{k}$ denote the position of $k$ in $\tbl$.
Clearly, $\Pos{\tbl}{1}=\pos{1,1}$ for any nonempty standard
tableau $\tbl$.

\subsection{}\label{defs:HookD}
Given two positions $\pos{i_1,j_1}$ and $\pos{i_2,j_2}$, we define
their \textit{hook distance} as
\begin{equation*}
  \HookDistance(\pos{i_1,j_1},\pos{i_2,j_2})
  = i_1 - i_2 - j_1 + j_2.\pagebreak[3]
\end{equation*}\pagebreak[3]%
For a single position $(i,j)$, we set
\begin{equation*}
  \HookDistance(\pos{i,j}) = i - j.
\end{equation*}
For a box $\Box$ in position $\pos{i,j}$, we set
$\HookDistance(\Box)=\HookDistance(\pos{i,j})$.

%% For a given Young tableau $\tbl$, $|\tbl|=n$, filled with numbers $1$,
%% \dots, $n$, and for two numbers $i$ and $j$, $1\leq i,j\leq n$, we
%% also write
%% \begin{equation*}
%%   \HookDistance(i,j)
%%   = \HookDistance_{\tbl}(i,j)
%%   = \HookDistance(\Pos{\tbl}{i},\Pos{\tbl}{j})
%% \end{equation*}
%% and 
%% \begin{equation*}
%%   \HookDistance(i)
%%   = \HookDistance_{\tbl}(i)
%%   = \HookDistance(\Pos{\tbl}{i}).
%% \end{equation*}

\subsection{Quantum dimension}\label{sec:qdim}For a Young diagram
$\dgr=(\dgr_1,\dots,\dgr_k)$, we define its quantum dimension as
\begin{equation*}
  \qDim(\dgr,q)
  =\mfrac{\displaystyle\prod_{j=1}^kq^{(j-1)\dgr_j}
    \prod_{i=1}^{|\dgr|}(1 - q^i)
  }{\displaystyle
    \prod_{\Box\,\in\dgr}(1 - q^{h(\Box)})},
\end{equation*}
where $h(\Box)$ is the length of the hook passing through a chosen box
(and $|\dgr|=\dgr_1+\dots+\dgr_k$).

For example,
\ytableausetup{smalltableaux}%
\begin{equation*}
  \qDim\Bigl(\Ytableaui{{}&{}&{}\\
    {}&{}\\
    {}&{}\\
    {}}\,,
  q\Bigr)
  =
  \mfrac{q^9 (1 - q) (1 - q^2) (1 - q^3) (1 - q^4) (1 - q^5) (1 - q^6) (1 - q^7) (1 - q^8)}{(1 - q)^3 (1 - q^2) (1 - q^3) (1 - q^4)^2 (1 - q^6)}.
\end{equation*}
\ytableausetup{nosmalltableaux}%

\section{Hecke algebras}\label{app:Hecke}
\subsection{}
By the Hecke algebra $\Hecke_{n}=\Hecke_{n}(\alpha,\beta)$, we mean
the Iwahori--Hecke algebra of type $A_n$ over~$\oC$. It is the
quotient of the braid group on $n$ strands, with generators $h_1$,
\dots, $h_{n-1}$, by the relations
\begin{equation}\label{Hecke-rel}
  (h_i - \alpha)(h_i - \beta) = 0,\qquad 1\leq i\leq n-1,
\end{equation}
where $\alpha$ and $\beta$ are two complex numbers (typically, such
that $\alpha+\beta\neq 0$ and $\alpha\beta\neq 0$).  The algebra
actually depends not on two but on one parameter, because $\alpha$,
$\beta$, and $h_i$ can be rescaled simultaneously.
%% \footnote{Popular choices of fixing this freedom of rescaling are to
%%   set $(\alpha,\beta)$ equal to $(1,-q^2)$, $(q,-q^{-1})$, or
%%   variations thereof; each has to produce the two eigenvalues of a
%%   transposition, $(+1,-1)$, as $\Hecke_{n}$ is evolved into the
%%   symmetric group on $n$ elements, which in the cases just listed
%%   means taking $q\to 1$.}

Thinking of the $\Hecke_{n}$ generators as coming from the braid
group, we use the braid-group diagram notation for them:
\begin{equation*}
  h_1 = \
  \begin{tangles}{l}
    \hx\step\id\step\id\step[1.5]\object{\raisebox{8pt}{$\dots$}}
  \end{tangles}\ \ ,
  \qquad
  h_2 = \
  \begin{tangles}{l}
    \id\step\hx\step\id\step[1.5]\object{\raisebox{8pt}{$\dots$}}
  \end{tangles}\ \ ,
  \qquad\dots.
\end{equation*}
Then Hecke relations~\eqref{Hecke-rel} take the graphic form
\begin{equation*}
  \begin{tangles}{l}
    \vstr{67}\hx\\
    \vstr{67}\hx
  \end{tangles}\ \
  = -\alpha\beta\ \ 
  \begin{tangles}{l}
    \vstr{133}\id\step\id
  \end{tangles}\ \ 
  + (\alpha+\beta)\ \
  \begin{tangles}{l}
    \vstr{133}\hx
  \end{tangles}.
\end{equation*}

\subsection{Jucys--Murphy elements}\label{JM-Hecke} In $\Hecke_{n}$,
we define a commuting family of Jucys--Murphy elements $\JMBare_i$,
$1\leq i\leq n$, as follows:
\begin{equation}\label{JM-h}
  \JMBare_1=1,\qquad
  \JMBare_{i}= (-\alpha \beta)^{-1}h_{n+1-i}\JMBare_{i-1}h_{n+1-i},
\end{equation}
$2\leq i\leq n$.  In terms of braid diagrams,
\begin{equation*}
  \JMBare_2 = (-\alpha\beta)^{-2}\ \
  \begin{tangles}{l}
    \id\step\step\id\step\id\step\id\step\id\\
    \id\step\object{\dots}\step\id\step\id\step\hx\\
    \id\step[2]\id\step\id\step\hx\\
    \id\step[2]\id\step\id\step\id\step\id
  \end{tangles}\ \ ,
  \quad
  \JMBare_3 = (-\alpha\beta)^{-2}\ \
  \begin{tangles}{l}
    \id\step\step\id\step\hx\step\id\\
    \id\step\object{\dots}\step\id\step\id\step\hx\\
    \id\step[2]\id\step\id\step\hx\\
    \id\step[2]\id\step\hx\step\id
  \end{tangles}\ \ ,
  \quad
  \JMBare_4 =  (-\alpha\beta)^{-3}\ \
  \begin{tangles}{l}
    \vstr{66}\id\step\step\hx\step\id\step\id\\
    \vstr{66}\id\step\step\id\step\hx\step\id\\
    \vstr{66}\id\step\object{\dots}\step\id\step\id\step\hx\\
    \vstr{66}\id\step\step\id\step\id\step\hx\\
    \vstr{66}\id\step[2]\id\step\hx\step\id\\
    \vstr{66}\id\step[2]\hx\step\id\step\id
  \end{tangles}\ \ ,
\end{equation*}
and so on.  Diagram manipulations immediately show that the
$\JMBare_i$ pairwise commute.  The actual choice of the $\JMBare_i$
family (which is not unique) and the labeling reflect our preferences
in the main body of the paper.

\subsection{Specht modules of $\Hecke_{n}$}\label{sec:ytr}
We essentially follow~\cite{[DJ-87]} in describing the action of
$\Hecke_{n}$ on its Specht modules.

A Specht module $S^{\dgr}$ of $\Hecke_{n}$ is associated with each
Young diagram $\dgr$, $|\dgr|=n$, and is defined for any values of
$\alpha$ and $\beta$. \ It has a basis labeled by all standard Young
tableaux of shape $\dgr$.  The $\Hecke_{n}$ generators $h_k$ act on
$S^{\dgr}$ by first mapping into a larger space $W^{\dgr}$ and then
taking the quotient by a set of relations $R$ such that
$W^{\dgr}/R=S^{\dgr}$:
\begin{equation}\label{two-arrow-action}
  \xymatrix@R=12pt@C=20pt{
  S^{\dgr}\ar[r]\ar@/_12pt/_{h_k}[rr]&
  W^{\dgr}\ar[r]^{\pi}&
  S^{\dgr},}
\end{equation}
where $\pi$ is the canonical projection.  The space $W^{\dgr}$ is the
linear span of all (not necessarily standard) Young tableaux obtained
by filling $\dgr$ with $1,\dots,n$, and $R$ are the Garnir
relations~\cite{[Welsh]}, which we describe below.

\subsubsection{}
The first short arrow in~\eqref{two-arrow-action} is defined as
follows.  We recall that every tableau $\tbl$ can be obtained by
applying an element $\sigma\in\oS_n$ to a reference tableau
$\tbl^{\text{s}}$ (which we choose as the row-reading superstandard
tableau), $\tbl^\sigma=\sigma\tbl^{\text{s}}$, where $\sigma$ acts
just by permuting the entries.  We write $\sigma$ as a reduced
(minimal-length) representation
$\sigma=\sigma_{i_1}\sigma_{i_2}\cdots\sigma_{i_k}$ in terms of
elementary transpositions $\sigma_i$, $i=1,\dots,n-1$, and,
accordingly,
\begin{equation*}
  \tbl^\sigma=\tbl^{i_1i_2\dots i_k}
  =\sigma_{i_1}\sigma_{i_2}\cdots\sigma_{i_k}\tbl^{\text{s}}.
\end{equation*} 
We 
%% call $k$ the length of the tableau and
write $k=\ell(\tbl^\sigma)$ (and set $\ell(\tbl^{\text{s}})=0$).  Then
the first arrow in~\eqref{two-arrow-action} is
\begin{equation}\label{eq:h-action}
  h_k\tbl^\sigma=
    \begin{cases}
      \sigma_k\tbl^\sigma,&
      \ell(\sigma_k\tbl^\sigma)=\ell(\tbl^\sigma)+1,
      \\
      (\alpha+\beta)\tbl^\sigma-\alpha\beta\sigma_k\tbl^\sigma&
      \text{otherwise}.
    \end{cases}
\end{equation}
Clearly, this gives a nonstandard Young tableau in general.  But
modulo the Garnir relations $R$, any nonstandard Young tableau can be
expressed as a linear combination of standard Young tableaux.

\subsubsection{Garnir relations}The set of Garnir relations $R$
consists of two subsets, $R=R_\alpha\cup R_\beta$.

The relations in $R_\alpha$ are those that make the rows of $\tbl$
standard: for any row $k_1$, \dots, $k_g$ with a ``disorder''
$k_i>k_{i+1}$, we order the offending numbers at the expense of the
factor $\alpha$ appearing in front of the tableau.  Hence, if $\tbl$
is a tableau with the total of $K$ instances of disorder in its rows,
then the corresponding relation in $R_{\alpha}$ is
\begin{equation*}
  \tbl = \alpha^K \widetilde{\tbl},
\end{equation*}
where $\widetilde{\tbl}$ is the corresponding row-standard
tableau.

The relations in $R_{\beta}$ allow linearly expressing any tableau
with nonstandard columns as linear combinations of column-standard
tableaux.  For any tableau with a transposition in a column,
\begin{equation*}
  \qquad\qquad\tbl\ \ =
  \vcenter{\Ytableaui{
      \cdot&\none[\dots]&\none[\dots]&\none[\dots]&\none[\dots]&\none[\dots]&\none[\dots]&\none[\dots]&\none[\dots]&\cdot
      \\
      \none&\none&\none&\none[\vdots]
      \\
      \cdot&\none[\dots]&\cdot&x_1&x_2&\none[\dots]&\none[\dots]&\none[\dots]&x_a
      \\
      y_1&\none[\dots]&y_{b'}&y_b&\cdot&\none[\dots]&\cdot
      \\
      \none&\none&\none&\none[\vdots]
      \\
      \cdot&\none[\dots]&\none[\dots]&\none[\dots]&\none[\dots]&\cdot
    }} \kern-200pt (b'=b-1),
\end{equation*}
where $x_1>y_b$, let $\Sh(x,y)$ be the set of all permutations of the
form
\begin{equation}\label{eq:sigma}
  \sigma=
  \begin{pmatrix}
    x_{j_1}&x_{j_2}&\dots&x_{j_r}\\
    y_{i_1} &y_{i_2}&\dots&y_{i_r}
  \end{pmatrix},\qquad 1\leq r\leq\min(a,b), 
\end{equation}
where $(x_{j_1},x_{j_2},\dots,x_{j_r})$ and
$(y_{i_1},y_{i_2},\dots,y_{i_r})$ are ordered subsets
of~$(x_1,\dots,x_a)$ and $y=(y_1,\dots,y_b)$. \ We set
\begin{align*}
  L(\sigma)&=r,
  \\
  w(\sigma)&=\sum_{k=1}^r\HookDistance(x_{j_k},y_{i_k}),
\end{align*} 
where the hook distance $\HookDistance({}\cdot{},{}\cdot{})$ is
defined in~\bref{defs:HookD}. \ We then have the relation
\begin{equation*}
  \tbl=-\sum_{\sigma\in\Sh(x,y)}\alpha^{L(\sigma)^2}
  \left(-\ffrac{\beta}{\alpha}\right)^{-w(\sigma)}%%%%%%%%
  \sigma\tbl,
\end{equation*}
where the $\sigma$ act by permuting the entries.  The set $R_\beta$
contains all such relations.

This defines $\pi$ and hence the action in~\eqref{two-arrow-action}:
by the repeated use of Garnir relations, every tableau in the
right-hand side of~\eqref{eq:h-action} is expressed as a linear
combination of standard Young tableaux.

\begin{example}
  We consider the tableau
  \begin{equation*}
    \tbl = \ \Ytableaui{
        1&4&5&7\\
        2&3&6
      }
  \end{equation*}
  with disorder in the second column.  The corresponding permutations
  are then given by {\scriptsize$\begin{pmatrix}
    5& 7\\
    2& 3
  \end{pmatrix}$, $\begin{pmatrix}
    4& 7\\
    2& 3
  \end{pmatrix}$, $\begin{pmatrix}
    4& 5\\
    2& 3
  \end{pmatrix}$, $\begin{pmatrix}
    7\\
    2
  \end{pmatrix}$, $\begin{pmatrix}
    7\\
    3
  \end{pmatrix}$, $\begin{pmatrix}
    5\\
    2
  \end{pmatrix}$, $\begin{pmatrix}
    5\\
    3
  \end{pmatrix}$, $\begin{pmatrix}
    4\\
    2
  \end{pmatrix}$}, and {\scriptsize$\begin{pmatrix}
    4\\
    3
  \end{pmatrix}$}, and applying the $R_\beta$ relations
yields\ytableausetup{smalltableaux}
\begin{align*}
    \tbl&=
    -\fffrac{\beta^6}{\alpha^2}~\Ytableaui{
        1&4&2&3\\
        5&7&6
      }
    +\fffrac{\beta^5}{\alpha}~\Ytableaui{
        1&2&5&3\\
        4&7&6
      }
    -\beta^4~\Ytableaui{
       1&2&3&7\\
       4&5&6
      }
    -\fffrac{\beta^4}{\alpha^3}~\Ytableaui{
        1&4&5&2\\
        7&3&6
      }
    \\
    &\quad{}
    +\fffrac{\beta^3}{\alpha^2}~\Ytableaui{
        1&4&5&3\\
        2&7&6
      }
    +\fffrac{\beta^3}{\alpha^2}\!~\Ytableaui{
        1&4&2&7\\
        5&3&6
      }
    -\fffrac{\beta^2}{\alpha}~\Ytableaui{
        1&4&3&7\\
        2&5&6
      }
    -\fffrac{\beta^2}{\alpha}~\Ytableaui{
        1&2&5&7\\
        4&3&6
      }
    +\beta~\Ytableaui{
        1&3&5&7\\
        2&4&6
      }\,.
  \end{align*}
  After applying the $R_\alpha$ relations to order the elements in
  each row, we obtain
  \begin{align*}
    \tbl&=
      -\alpha\beta^6\;\Ytableaui{
        1&2&3&4\\
        5&6&7
      }
    +\alpha\beta^5\;\Ytableaui{
        1&2&3&5\\
        4&6&7
      }
   -\beta^4\;\Ytableaui{
        1&2&3&7\\
        4&5&6
      }
    -\alpha\beta^4\;\Ytableaui{
        1&2&4&5\\
        3&6&7
      }
    \\
    &\quad{}
    +\alpha\beta^3\;\Ytableaui{
        1&3&4&5\\
        2&6&7
      }
    +\beta^3\;\Ytableaui{
        1&2&4&7\\
        3&5&6
      }
    -\beta^2\;\Ytableaui{
        1&3&4&7\\
        2&5&6
      }
    -\beta^2\;\Ytableaui{
        1&2&5&7\\
        3&4&6
      }
    +\beta\;\Ytableaui{
        1&3&5&7\\
        2&4&6
      }\,.
  \end{align*}\ytableausetup{nosmalltableaux}%
\end{example}

\subsection{Seminormal representations of $\Hecke_{n}(\alpha,\beta)$}
We briefly review the facts that we need about the seminormal
representations of Hecke algebras~\cite{[We-88]} (also see
\cite{[Ram-semi]} and the references therein).

\subsubsection{}\label{weights-f-white}
For a standard $n$-box tableau $\tbl$ filled with $1$, \dots,
$n$, we define its weight
\begin{equation*}
  \wt(\tbl)=(\wt_1(\tbl),\dots,\wt_{n-1}(\tbl))\in\oC^{n-1},\qquad
  \wt_i(\tbl) = \Bigl(-\ffrac{\beta}{\alpha}\Bigr)^{\HookDistance(i+1)}
\end{equation*}
For example,\ytableausetup{smalltableaux}%
\begin{equation*}
  \tbl = \ \ \Ytableaui{
    1&3&6\\
    2&4&8\\
    5\\
    7
  }\ \ \Longrightarrow \wt(\tbl)=
  (-\ffrac{\beta}{\alpha},-\ffrac{\alpha}{\beta},
  1,
  \ffrac{\beta^2}{\alpha^2},
  \ffrac{\alpha^2}{\beta^2},
  -\ffrac{\beta^3}{\alpha^3},
  -\ffrac{\alpha}{\beta})
\end{equation*}\ytableausetup{nosmalltableaux}%

\subsubsection{}\label{sec:h-acts}
For generic $\alpha$ and $\beta$, given a Young diagram $\lambda$, the
seminormal representation $\modV^{\lambda}$ is defined by specifying
the $\Hecke_n(\alpha,\beta)$ action on a basis of standard tableaux of
shape $\lambda$. \ The generators act on such a tableau $\tbl$
as\footnote{That the left-hand side involves $h_{n-i}$ rather than the
  more standard $h_i$ means that a Hecke algebra automorphism is
  applied to a standard construction, which is done for our purposes
  in the main text.}
\begin{equation}\label{h-acts-first}
  h_{n-i}\acts \tbl=
  \begin{cases}
    \alpha\tbl,& \text{$i$ and $i+1$ are in the same row},\\
    \beta\tbl,& \text{$i$ and $i+1$ are in the same column},\\
    \ffrac{\alpha + \beta}{1 - \ffrac{w_{i-1}}{w_i}}\tbl
    +\eta\bigl(\ffrac{w_{i-1}}{w_i}\bigr)\tbl_{(i,i+1)}& \text{otherwise},
  \end{cases}
\end{equation}
where $w=\wt(\tbl)$, and $\tbl\mapsto\tbl_{(i,i+1)}$ is the
transposition of $i$ and $i+1$.\footnote{In the third line
  in~\eqref{h-acts-first}, $\tbl_{(i,i+1)}$ is again a standard
  tableau. The third line cannot occur for $i=1$ just because $1$ and
  $2$ are necessarily in the same row or column.}\pagebreak[3]
%% \footnote{For example, \ytableausetup{smalltableaux}%
%%   \Ytableaui{
%%     1 & 2 & 5 \\
%%     3 & 4 & 6 } \ is sent into a standard tableau by $(2,3)$ and
%%   $(4,5)$, and is not by $(1,2)$, $(3,4)$, and
%%   $(5,6)$.\ytableausetup{nosmalltableaux}}\;
The function $\eta(\cdot)$ is a convenient way to reflect the
arbitrariness of rescaling the basis elements.  For the Hecke-algebra
property~\eqref{Hecke-rel} to hold for~\eqref{h-acts-first}, as is
easy to see, $\eta$ must satisfy the relation
\begin{equation}\label{eta-inv}
  \eta(x)\eta\bigl(\ffrac{1}{x}\bigr) = -\alpha \beta
  \,\theF(x),
\end{equation}
where
\begin{equation}\label{theF}
  \theF(x)=\mfrac{
    \Bigl(\!1+\fffrac{x\beta}{\alpha}\Bigr)
    \Bigl(\!1+\fffrac{x\alpha}{\beta}\Bigr)}{(1-x)^2}.
\end{equation}

The same relation~\eqref{eta-inv} suffices for the Yang--Baxter
equation $h_i h_{i+1} h_i = h_{i+1} h_{i} h_{i+1}$ to hold for the
action in~\eqref{h-acts-first}. \ To see this, note that
$\widetilde{w}=\wt(\tbl_{(i,i+1)})$ differs from $w$ \textit{by the
  transposition of the $(i-1)$th and $i$th components}.  Verifying the
Yang--Baxter equation then amounts to a straightforward analysis of
the possible cases.  For example, if both $h_{n-i}$ and $h_{n-i-1}$
act in accordance with the third line in~\eqref{h-acts-first}, then
\begin{multline*}
  h_{n-i} h_{n-i-1} h_{n-i}\tbl - h_{n-i-1} h_{n-i} h_{n-i-1}\tbl=
  (\alpha +\beta)\times
  \\
  \biggl(\mfrac{(\alpha +\beta)^2 w_{i} w_{i+1} 
    (w_{i-1} w_{i+1}-w_{i}^2)}{(w_{i-1}-w_{i})^2 (w_{i}-w_{i+1})^2}  
  +\ffrac{ w_{i+1}
    \Bigl(\!\eta\bigl(\ffrac{w_{i}}{w_{i+1}}\bigr)\eta\bigl(\ffrac{w_{i+1}}{w_{i}}\bigr)-\eta\bigl(\ffrac{w_{i-1}}{w_{i}}\bigr)
    \eta\bigl(\ffrac{w_{i}}{w_{i-1}}\bigr)\!\Bigr)}{w_{i-1}-w_{i+1}}
  \biggr)\tbl,
\end{multline*}
which does vanish by virtue of~\eqref{eta-inv}. \ The other cases
(where one of the generators acts by an eigenvalue) are verified
similarly.

A ``symmetric'' choice for $\eta(\cdot)$ satisfying~\eqref{eta-inv} is
\begin{equation*}
  \eta(x)=\sqrt{-\alpha\beta\theF(x)}.
\end{equation*}

\subsubsection{} In the seminormal representation $\modV^{\lambda}$,
Jucys--Murphy elements~\eqref{JM-h} are diagonalized,
\begin{equation*}
  J_i \acts \tbl  = \wt_{i-1}(\tbl)\tbl,
  \qquad
  2\leq i\leq n,
\end{equation*}
for each standard tableau of the given shape $\lambda$.

\section{Two-variate identities}\label{app:2-variate}
We prove identities~\eqref{id1} and~\eqref{id2} by actually proving an
apparently stronger statement in~\bref{lemma:identity}, which implies
that the identities in fact hold in a generalized, two-variate form,
in which they border with generalized Pieri rules (see~\cite{[GHXZ]}
and the references therein).

Instead of a single indeterminate $\q=-\beta/\alpha$ in
Sec.~\ref{sec:seminormal}, we introduce two variables $x$ and $t$ and
associate them with vertical and horizontal distances $\Deltav$ and
$\Deltah$ between boxes in a Young diagram.  To facilitate the
comparison with the Pieri rule formulas in the literature (see,
e.g.,~\cite{[GHXZ]}), we define the weight of a box in terms of its
coleg and coarm:
\begin{equation*}
  \weight(\Box)
  = x^{l'(\Box)} t^{a'(\Box)}.
\end{equation*}

\begin{Lemma}\label{lemma:identity}
  Let $\dgr$ be a Young diagram with $\ell$ corners \ $\Circ_i$, \
  $1\leq i\leq\ell$, and the addable boxes \ $\Star_k$, \
  $0\leq k\leq\ell$. \ Set $\bx_i= x\,t\,\weight(\Circ_i)$ and
  $\bu_k=\weight(\Star_k)$. \ Then
  \begin{gather*}
    \sum_{k=0}^{\ell}
    \ffrac{1}{\bu_k^j}\,
    \mfrac{1}{\displaystyle
      \prod_{\substack{i=0\\ i\neq k}}^{\ell}
      \bigl(1 - \frac{\bu_i}{\bu_k}\bigr)}
    =
    \begin{cases}
      1,&j=0,\\
      0,&1\leq j\leq \ell,\\
      \ffrac{(-1)^\ell}{\bx_1\dots \bx_{\ell}},
      &j=\ell+1.
    \end{cases}
  \end{gather*}
\end{Lemma}
To prove this, we simply note that for \textit{any} pairwise distinct
complex numbers $u_0$, $u_1$, \dots, $u_{\ell}$, the function
\begin{equation*}
  F_j(z)= \mfrac{z^j}{\displaystyle\prod_{i=0}^{\ell}(z - u_i)}
\end{equation*}
has the property that
\begin{equation*}
  \sum_{i=0}^{\ell}\res_{u_i}F_j(z) =
  \sum_{i=0}^{\ell}
  \mfrac{u_{i}^j}{\displaystyle
    \prod_{\substack{k=0\\ k\neq i}}^{\ell}(u_{i} - u_{k})},
\end{equation*}
where the sum goes over all residues at finite $z$, but at the same
time the residue at infinity is
\begin{equation*}
  \res_{\infty}F_j(z) =
  \begin{cases}
    0,&1\leq j\leq \ell-1,\\
    -1,&j=\ell,\\
    -h_{j-\ell}(u_0,\dots,u_{\ell}),&j\geq \ell,
  \end{cases}
\end{equation*}
where $h_n$ are the complete homogeneous symmetric polynomials
\begin{gather*}
  h_n(x)=\sum_{i_1\leq\dots\leq i_n}x_{i_1}\dots x_{i_n}.
\end{gather*}
In applying this with $u_i=\bu_i$ as defined above, it remains to
recall (see, e.g.,~\cite{[GT-qt]}) that
\begin{gather*}
  \bu_0\dots\bu_{\ell} = \bx_1\dots \bx_{\ell}.
\end{gather*}

\subsection{}
It follows from the lemma, in particular, that
\begin{align}\label{ID1}
  \sum_{\Star\,\in\,\cocorners{\dgr}}
  \mfrac{\displaystyle
  \prod_{\Circ\in\corners{\dgr\,}}\Bigl(\!1-x^{\Deltav(\Star,\,\Circ)}t^{\Deltah(\Star,\,\Circ)}\Bigr)
  }{\displaystyle
  \prod_{\Star'\neq\Star}\Bigl(\!1-x^{\Deltav(\Star,\Star')}t^{\Deltah(\Star,\Star')}\Bigr)}
  &=1
  \\
  \intertext{and}
  \label{ID2}
  \sum_{\Star\,\in\,\cocorners{\dgr}}
  \mfrac{\displaystyle
  \prod_{\Circ\in\corners{\dgr\,}\;,\;\; {\Circ\neq\Ast}}
  \Bigl(\!1-x^{\Deltav(\Star,\,\Circ)}t^{\Deltah(\Star,\,\Circ)}\Bigr)
  }{\displaystyle
  \prod_{\Star'\neq\Star}
  \Bigl(\!1-x^{\Deltav(\Star,\Star')}t^{\Deltah(\Star,\Star')}\Bigr)}
  &=1,
\end{align}
which are two-variate forms of~\eqref{id1} and~\eqref{id2}.  (As
before, $\Corners{\lambda}$ are the corners of a Young diagram
$\lambda$ and $\Cocorners{\lambda}$ are the boxes addable to it, and
$\Ast$ is a fixed corner of $\dgr$).

\subsection{Examples} We illustrate~\eqref{ID1}
with
%% \ytableausetup{smalltableaux}
%% $\dgr=\ \ \Ytableaui{
%%   {}&&&\\
%%   {}&\\
%%   {} }$\ytableausetup{nosmalltableaux}\,,
$\lambda=(4,2,1)$,
which has the following addable boxes $\Star$ and corners $\Circ$:
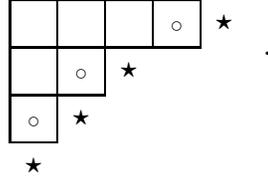
\begin{equation*}
  \Ytableaui{
    {}&&&\Circ&\none[\Star]\\
    {}&\Circ&\none[\Star]\\
    \Circ&\none[\Star]\\
    \none[\Star]
  }\ \ .
\end{equation*}
Identity~\eqref{ID1} then takes the form
\begin{multline}\label{ide-ex-1}
  \mfrac{(1-x) (1-x^2 t^{-1})
    (1-x^3t^{-3})}{(1-x t^{-1}) (1-x^2 t^{-2})
    (1-x^3t^{-4})}
  +\mfrac{(1-t) (1-x)
    (1-x^2 t^{-2})}{(1-x^{-1} t) (1-x t^{-1})
    (1-x^2 t^{-3})}
  \\
  +\mfrac{(1-t)
    (1 - x^{-1}t^2)
    (1-x t^{-1})}{(1 - x^{-2}t^2) (1-x^{-1} t)
    (1-x t^{-2})}
  +\mfrac{(1-t) (1 - x^{-2} t^4)
   (1-x^{-1} t^3)}{(1 - x^{-3}t^4) (1-x^{-2}t^3)
   (1- x^{-1}t^2)} = 1.
\end{multline}
With a corner $\Ast$ selected as\ytableausetup{smalltableaux}
  $(\dgr,\Ast)=\ \
  \Ytableaui{
    {}&&&\\
    {}&\Ast\\
    {}
  }$\ytableausetup{nosmalltableaux}\,,
the corresponding identity~\eqref{ID2} becomes the following
``thinning'' of~\eqref{ide-ex-1}:
\begin{multline*}
  \mfrac{(1-x) (1-x^3 t^{-3})}{(1-x t^{-1}) (1-x^2 t^{-2})
    (1-x^3 t^{-4})}
  +\mfrac{(1-t)
    (1-x^2 t^{-2})}{(1-x^{-1}t) (1-x t^{-1})
    (1-x^2 t^{-3})}
  \\
  +\mfrac{(1-x^{-1}t^2)
    (1-x t^{-1})}{(1-x^{-2}t^2) (1-x^{-1}t)
    (1-x t^{-2})}
  + \mfrac{(1-t)
    (1- x^{-2}t^4)}{(1-x^{-3}t^4)
    (1- x^{-2}t^3)
    (1- x^{-1}t^2)}
 = 1.
\end{multline*}
Another, yet ``thinner,'' identity was needed in the proof
of~\bref{lemma:1234-zero}. 

\begin{Rem}
  Lemma~\bref{lemma:identity} implies a (simple) part of the
  identities discussed in~\cite{[GHXZ],[GT-qt]} (also see the
  references therein). \ If we set
  \begin{gather*}
    d_{\dgr^{(k)},\dgr}(t,x)=\ffrac{1}{\bu_k}\,
    \mfrac{\displaystyle
      \prod_{i=1}^{\ell}\Bigl(\!1 - \frac{\bx_i}{\bu_k}\Bigr)
    }{\displaystyle
      \prod_{\substack{i=0\\ i\neq k}}^{\ell}
      \Bigl(\!1 - \frac{\bu_i}{\bu_k}\Bigr)} ,
  \end{gather*}
  then the identity %%%(see, e.g.,~\cite{[GHXZ]})
  \begin{align*}
    \sum_{k=0}^{\ell} d_{\dgr^{(k)},\dgr}(t,x)
    &
    = 1
    \\
    \intertext{follows from~\bref{lemma:identity} by expanding the
    brackets in the numerator, after which we are left with only the
    top-degree term in the $\bu_i$. \ Similarly, we also have}
    \sum_{k=0}^{\ell}\bu_k d_{\dgr^{(k)},\dgr}(t,x)&= 1,
  \end{align*}
  where just the constant term in the numerator contributes to produce
  the $1$ in the right-hand side.
  
%%   We finally remark that from~\cite{[GT-qt]} (also see~\cite{[GHXZ]})
%%   we know that
%%   \begin{equation*}
%%     \sum_{k=0}^{\ell} (\bu_k)^j d_{\dgr^{(k)},\dgr}(t,x)=
%%     (-1)^{j-1} e_{j-1}[D_{\dgr}(t,x)],
%%     \qquad j\geq 1,
%%   \end{equation*}
%%   where
%%   \begin{align*}
%%   D_{\dgr}(t,x) &= \sum_{i=1}^{\ell}\bx_i - \sum_{i=1}^{\ell+1} \bu_i,
%%   \end{align*}
%%   the elementary symmetric functions are
%%   \begin{gather*}
%%     e_n(x)=\sum_{i_1<\dots<i_n}x_{i_1}\dots x_{i_n},
%%   \end{gather*}
%%   and $F[G]$ denotes plethystic substitution (see,
%%   e.g.,~\cite{[expose],[random]}).
\end{Rem}

\parindent0pt

\end{document}